\documentclass[12pt]{amsart}
\usepackage{amssymb}
\usepackage[all]{xy}
\usepackage{hyperref}

\newtheorem{theorem}{Theorem}[section]
\newtheorem{prop}[theorem]{Proposition}

\newtheorem{coro}[theorem]{Corollary}
\newtheorem{lemma}[theorem]{Lemma}
\newtheorem{claim}[theorem]{Claim}
\newtheorem{rem}[theorem]{Remark}

\def\wh #1{{\widehat{#1}}}
\def\w #1{{\widetilde{#1}}}
\def\Id{\operatorname{Id}}
\DeclareMathOperator{\Nm}{{Nm}}

\DeclareMathOperator{\id}{{id}}
\DeclareMathOperator{\pic}{{Pic}}
\DeclareMathOperator{\Aut}{{Aut}}
\DeclareMathOperator{\spp}{{ps}}

\allowdisplaybreaks
\newcount\hbadness \let\vbadness=\hbadness

\begin{document}
\title[Prym Varieties and fourfold covers]{Prym varieties and fourfold covers}
\author{Sev\'\i n Recillas}
\address{Instituto de Matem\'aticas,
UNAM Campus Morelia, Morelia, Mich. 58190, M\'exico }
 \email{sevin@unam.mx}
\address{CIMAT, Callej\'on Jalisco s/n,
Valenciana, Guanajuato, M\'exico}
\email{sevin@fractal.cimat.mx}
\author{Rub\'\i\ E. Rodr\'\i guez}
\address{Facultad de Matem\'aticas, Pontificia Universidad Cat\'olica
de Chile, Casilla 306, Correo 22, Santiago, Chile}
\email{rubi@mat.puc.cl}
\thanks{This research was partially supported by CONACYT--CONICYT Exchange
Program, by grants CONACYT 27962--E, Fondecyt 1000623 and
Presidential Science Chair on Geometry.}

\maketitle

\tableofcontents

\section{Introduction}

The problem of finding an isogenous decomposition of the Jacobian
$J\widetilde{C}$ of a smooth connected complete curve
$\widetilde{C}$ has been studied by many people, as well as its
generalization for any abelian variety.

The first example of finding subtori $X$ and $Y$ of
$J\widetilde{C}$ such that $X \times Y$ is isogenous to
$J\widetilde{C}$ is due to Wirtinger \cite{w}; the situation is a
double cover of curves $f : \widetilde{C} \to C$, where
$J\widetilde{C}$ is found to be isogenous to the product of the
abelian subvariety $f^{*}(JC)$ and its natural complement
$P(\widetilde{C}/C)$, later called the Prym variety of the cover
$f$, by Mumford in \cite{mumprym}.

In this last work one finds a careful discussion of the kernels of
the isogeny and of $f^{*}$, which is then used to describe the
type of polarization that $P(\widetilde{C}/C)$ inherits from
$J\widetilde{C}$.

These ideas were generalized by Ries in \cite{ri} where he studied
the Prym variety associated to a cyclic unramified cover of a
hyperelliptic Riemann surface, and in \cite{gg} for the case of a
curve with an automorphism of prime order.

In a more general context, their results may be viewed as
particular cases of the isotypical decomposition of an abelian
variety with a finite group acting on it. Such decompositions have
been studied by several people in different settings (Donagi,
Kanev, M\'erindol, Lange and Recillas, among others), with
applications to the theory of integrable systems and the moduli
spaces of principal bundles of curves.

In another related direction, there are two constructions that
have proved very useful in the study of relations between
Jacobians and Prym varieties. The Recillas trigonal construction
(c.f. \cite{r}) shows that the Jacobian of a tetragonal curve is
isomorphic to the Prym variety of a double cover of a trigonal
curve, and Pantazis  showed in \cite{p} that the Prym varieties
associated to the bigonal construction of Donagi (c.f. \cite{d})
are dual to each other.

The aim of this work is the study of the structure of the Jacobian
of the Galois extension $W \to T$ of a fourfold cover of smooth
connected complete curves $X \to T$. We give explicit descriptions
of the rational simple components of this Jacobian as generalized
Pryms of intermediate covers of $W \to T$ and analyze the
polarization types and kernels of the isogenies involved.

Even though the isotypical decomposition of $JW$ may be obtained
faster from general methods (see, for instance,  \cite{d2},
\cite{kanev} and \cite{merindol}), our calculations for the
kernels of the natural isogenies that appear allow us to obtain
new proofs of the bigonal and trigonal constructions  in purely
group-theoretical terms, by considering $T = \mathbb{P}^1$ and
particular values of the ramification of $X \to T$ (see Remark
\ref{rem:orient} and Section \ref{subsec:big}, Remark
\ref{rem:classical} and Theorem \ref{thm:classical}).

Furthermore, by applying appropriate restrictions for the
ramification data of $X \to T$, we also obtain families of
Jacobian varieties isogenous to products of Jacobians and families
of Prym varieties isogenous to the product of elliptic curves.

When working with the Jacobian of a curve, we use additive or
tensor notation depending on whether the points of the Jacobian
are considered as points of an abelian variety or as line bundles
on the curve.

If $Y$ is a subvariety of an abelian variety, we denote by $Y^0$
its connected component of the origin, and by $Y[d]$ the subgroup
of $Y$ consisting of its points of order $d$.

If $f : \w C \to C$ is an unramified cyclic cover of curves of
degree $d$, we will denote by $\eta_f$ a point in $JC[d]$ which
determines the cover $f$; equivalently, a generator for the kernel
of $f^*$.

If $A$ is a polarized abelian variety with polarization $\lambda :
A \to \wh A$, where $\wh A$ is the dual abelian variety, then
$K(\lambda)$ will denote the kernel of $\lambda$. If $\lambda$ is
of type $(d_1, \ldots, d_g)$ with $d_i$ dividing $d_{i+1}$ for $1
\leq i \leq g-1$, we will denote by $\lambda_{\wh A} : \wh A \to
\widehat{\widehat{A}} = A$ the polarization on $\wh A$ such that
$\lambda_{\wh A} \circ \lambda$ is multiplication by $d_g$ on $A$.

If $f : A \to B$ is a morphism between abelian varieties, the dual
map will be denoted by $\wh f : \wh B \to \wh A$; if $B$ has a
polarization $\lambda_B$, then $A$ inherits the induced
polarization $\lambda_A$ defined by $\lambda_A = \wh f \circ
\lambda_B \circ f$ (the pullback of $\lambda_B$ via $f$).

\section{Prym varieties for covers of curves}

>From now on,  $\w C$ and $C$ are smooth connected complete curves
over the complex numbers ${\mathbb C}$, of respective genera $\w
g$ and $g$ , and $f : \w C \to C$ is a cover of degree $d$.

Let $( J, \lambda )$ and $(\w J, \w \lambda )$ denote the
Jacobians of $C$ and $\w C$, respectively. The cover $f$ induces
two homomorphisms between Jacobians, denoted by $ f^* : J
\longrightarrow  \w J$ and $\Nm f : \w J \longrightarrow  J$
respectively. Just as in \cite{mumprym} for the case $d=2$, they
are related as follows.

\begin{equation}
\label{nm}
  \Nm f \circ  f^* =  d \, \cdot \,  \Id \vert _{J} \quad \text{ and }
  \quad   \widehat{\Nm f} = \w \lambda \circ f^* \circ \lambda^{-1}
\end{equation}

Also note that $(\ker \Nm f)^0=(\ker(f^* {\circ} \Nm f))^0$.

\medskip

We will use the following characterization given in \cite[p.
337]{lange:cav92}.

\begin{prop}
\label{prop:lange}
Given a  cover of curves $f : \widetilde{C} \to
C$, the induced homomorphism $f^* : JC \to J\widetilde{C}$ is not
injective if and only if $f$  factors via a cyclic \'etale
covering $f'$ of degree $\geq 2\; $ as is shown in the following
diagram.
 $$ \xymatrix{ \widetilde{C} \ar[rr]^{f} \ar[dr]_{f''}  &
 &  C \\
                                         & C' \ar[ur]_{f'} & \ \ \ .
 }
 $$
\end{prop}

\begin{rem}
\label{rem:lange}
Note that by the proof of the above proposition,
$\ker f^{'*}$ is cyclic of order equal to the degree of $f'$;
hence if $f^{''*}$ is injective, then $\ker f^{*}$ is cyclic of
order equal to the degree of $f'$.
\end{rem}

\medskip

The following result is well known.

\begin{lemma}
\label{lem:kernm}
Let $f : \w C \to C$ be a  cover of curves of
degree $d$.

Then the number of connected components of $\ker \Nm f$ is the
cardinality of $\ker f^{*}$.

In particular, if $d$ is prime and $f$ is ramified or if $d$ is
prime and $f$ is not a cyclic cover then $\ker \Nm f$ is
connected.
\end{lemma}

\medskip

We now recall the definition of the Prym variety of a subtorus of
a principally polarized abelian variety (p.p.a.v.) and some of its
properties (c.f. \cite{ri} and \cite{mumprym}).

Let $(A, \lambda)$ be a p.p.a.v., where $\lambda : A \to \w A$ is
the principal polarization and let $X \overset i \hookrightarrow
A$ be a subtorus. Then  the \textit{Prym variety} $P = P(A,
\lambda , X)$ is defined by $P=\ker (\wh i \, {\circ} \, \lambda )
= \lambda^{-1}(\ker \wh i)$. Observe that $P$ is connected (since
$\ker \wh i$ is injective) and its dimension is the codimension of
$X$ in $A$.

If we denote by $j : P \to A $ the inclusion and by $\lambda _X$
and $\lambda _P$ the induced polarizations on $X$ and $P$,
respectively, we obtain the following commutative diagram:

\begin{equation}
\label{dia:prym}
\xymatrix{
  X \ar[d]_{\lambda_X} \ar[r]^{i} & A \ar[d]_{\lambda}  & P \ar[l]_{j} \ar[d]^{\lambda_P}\\
  \wh X   & \wh A  \ar[l]^{\wh i} \ar[r]_{\wh j} & \wh P
  }
\end{equation}

The following results appear in \cite{ri}.

\begin{prop}
 \label{prop:ries}
 Let $(A, \lambda)$ be a p.p.a.v., let $X$ be a
subtorus of $A$ and let $P = P(A, \lambda , X)$ be its Prym
variety.
Then
\begin{enumerate}

\item[i)] The sum homomorphism
\begin{align*}
 \sigma: X \times P & \longrightarrow A \\
 (x,y) & \longmapsto x + y
\end{align*}
is an isogeny with kernel given by
 $$ \ker \sigma= \{\ (x, -x)\ : \ x \in X \cap P\ \} .
 $$

\item[ii)] If we denote by $\lambda _{X\times P}$ the
pullback of the polarization $\lambda$ to $X\times P$ via
$\sigma$, then
$$\lambda _{X\times P} = \wh \sigma {\circ} \lambda {\circ} \sigma
=   \left ( \begin{matrix} \lambda_X & 0\\ 0 & \lambda
_P\end{matrix} \right )$$

\item[iii)] The kernels of the
induced polarizations are related as follows. $$
K(\lambda_{X\times P})=K(\lambda_X)\times K(\lambda_P)
 $$

Furthermore, $\ker \sigma$ is a maximal isotropic subgroup.

\vskip12pt

\item[iv)] $K(\lambda_X) = X \cap P = K(\lambda_P)$.
\end{enumerate}
\end{prop}

\medskip

We now apply these results to our cover of curves $f : \w C \to
C$. In this case the \textit{Prym Variety of $f$} is defined by
 $$
 P = P(f) = P(\w C/C)  = (\ker ( \wh {f^{\ast}} {\circ}
 \w{\lambda}))^0
 $$
and we can prove the following results.

\begin{theorem}
\label{prop1}
 Let $f : \w C \to C$ be a  cover of curves of degree
$d$. Denote by ${\tilde \lambda}_{f^{*}(J)}$ and by ${\tilde
\lambda}_{P}$ the polarizations induced by $\w \lambda$ on
$f^{*}(J)$ and on $P$, respectively.
Then
\begin{enumerate}
\item[i)] $P = (\ker \Nm f)^0$;
\item[ii)] $P = P(\w J , \w \lambda , f^{*}(J))$; i .e., $P$ is the
Prym variety of the subtorus $f^{*}(J)$ of the p.p.a.v. $(\w J ,
\w \lambda )$;
\item[iii)] the polarization on $J$ given by the pullback of
${\tilde \lambda}_{f^{*}(J)}$ via $f^{*}$ is $d \lambda$ and
$H_0:=\ker f^{*} \subset J[d]$ is isotropic with respect to the
Weil form associated to $d \lambda$;
\item[iv)] $f^{*}$ induces an isomorphism
 $$ H_0^{\bot }/H_0 \longrightarrow \ker (\lambda_{f^*J}) \, = \, \ker
(\lambda_P)  = \, f^*J \cap P \, \subseteq P[d]
 $$
 where orthogonality is with respect to the skew-symmetric multiplicative
form   $e_{d \lambda}$ on $J[d]$, and
\begin{equation*}
\label{eq:natural4}
 \bigl| f^*J \cap P \bigl| = \dfrac{\bigl| J[d] \bigl|}{\bigl| \ker f^* \bigl|^2}
 \ .
\end{equation*}
\item[v)] Consider the isogeny
\begin{equation*}
\phi \; : J \times P \to \w{J}  \  , \ \ \phi(c,\tilde{c}) =
f^*(c) + \tilde{c}
\end{equation*}
and the projection onto the first factor $\pi_1 : J \times P \to
J$.

Set $H_1 = \pi_1(\ker \phi)$. Then
 \begin{align*}
 \ker \phi = \{ (c, -f^*(c)) : c \in J[d]  \text{ and } f^*(c) \in P
 \} & \longrightarrow H_1 \\
 (c, -f^*(c)) & \dashrightarrow c
 \end{align*}
 is an isomorphism.

Furthermore, $H_1 = H_0^{\bot }$ and
\begin{equation*}
\label{eq:natural3}
 \bigl| \ker \phi \bigl| = \bigl| H_1 \bigl| = \bigl| H_0^{\bot } \bigl|=
 \dfrac{\bigl| J[d] \bigl|}{\bigl| \ker f^* \bigl|} \ .
\end{equation*}
\end{enumerate}
\end{theorem}
\begin{proof}
i) follows from $\lambda^{-1} {\circ} \wh {f^{\ast}} {\circ}
\w{\lambda} = \Nm f$.

ii) holds since  $f^{*}$ has finite kernel.

iii) is just a computation involving the relations given in
(\ref{nm}): $$ \wh{f^{*}} \circ {\w \lambda} \circ f^{*} = \lambda
\circ \Nm f \circ {\w \lambda}^{-1} \circ {\w \lambda} \circ f^{*}
= \lambda \circ \Nm f \circ f^{*} = d \,
 \lambda
 $$

As for iv), it follows from descent theory that $f^{*} : J \to
f^{*}J$ induces an isomorphism from $H_0^{\bot }/H_0$ to $\ker
(\lambda_{f^*J})$; therefore,
 $$
 H_0^{\bot} = (f^{*})^{-1}(\ker \lambda_{f^*J}) \, .
 $$

But we already know from Proposition \ref{prop:ries} that
 $$
 \ker (\lambda_{f^*J}) = \ker (\lambda_{P}) = f^*J \cap P \, .
 $$
Observing that $y \in f^*J \cap P$ if and only if $y$ is in $P$
and  $y = f^{*}(x)$ for some $x \in J$, we obtain $0 = \Nm f (y) =
\Nm f (f^{*}(x)) = dx$, and therefore $dy =0$.

The equality
\begin{equation}
\label{eq:lattices} \bigl| H_0 \bigl| = \dfrac{\bigl| J[d]
\bigl|}{\bigl| H_0^{\bot} \bigl|}
\end{equation}
will complete the proof of iv).

To prove it, consider $J = V/L$ and $f^{*}J = V/M$, with $L
\subseteq M$ lattices in $V$, and recall that $H_0 = M/L$,
$H_0^{\bot} = M^{\bot}/L$, $J[d] = K(d \, \lambda) = L^{\bot}/L$,
and that the index of $L$ in $M$ equals the index of $ M^{\bot}$
in $L^{\bot}$.

Then from $L \subseteq M \subseteq M^{\bot} \subseteq L^{\bot}$ it
follows that $[ L^{\bot} : L ] = [ L^{\bot} : M^{\bot} ] \cdot [
M^{\bot} : L ]$; since $[ L^{\bot} : M^{\bot} ] = [ M : L ]$, we
obtain (\ref{eq:lattices}) in the form $[ L^{\bot} : L ] = [ M : L
] \cdot [ M^{\bot} : L ]$.

As for v), it is clear that $\ker \phi \subseteq J[d] \times P[d]$
and therefore $\phi$ is an isogeny. It is also clear that $s : H_1
 \to \ker \phi$ defined by $c \mapsto
(c,-f^{*}(c))$ is a section of $\pi_1 : \ker \phi \to H_1$, so it
follows that $\pi_1$ is an isomorphism.

But also $H_1 = \{c \in J : f^{*}(c) \in P  \} =
(f^{*})^{-1}(f^{*}J \cap P) = H_0^{\bot }$. Therefore $ \bigl|
\ker \phi  \bigl| =  \bigl| H_0^{\bot } \bigl|$, which together
with equality (\ref{eq:lattices}) complete the proof of v) and of
the Theorem.
\end{proof}

Our next result gives the relation between the Prym variety of the
composition of two covers of curves and the Prym varieties of the
intermediate covers.

\begin{prop}
\label{prop:comp}
Let $f : X \to Y$ and $g : Y \to Z$ be two
covers of curves. Set $h = g \circ f : X \to Z$.

Then there are two isogenies given as follows.
 $$\spp : P(Y/Z) \times P(X/Y) \to P(X/Z) \, , \ \ \spp (y,x) = x + f^*y
 $$
and
 $$\psi : JZ \times P(Y/Z) \times P(X/Y) \to JX \, , \ \ \psi
(z, y, x) = x + f^*y + h^*z  \, .
 $$

Furthermore,
\begin{itemize}
\item[i)] The kernel of $\spp$ is contained in
$P(Y/Z)[\deg f] \times P(X/Y)[\deg f]$ and its cardinality is
equal to
 $$ \bigl| \ker \spp \bigl| = \bigl| P(Y/Z)[\deg f] \bigl|
\; \dfrac{\bigl| g^*(JZ) \cap \ker f^* \bigl|}{\bigl| \ker f^*
\bigl|} \, . $$
\item[ii)] The kernel of $\psi$ has cardinality equal to
$$ \bigl| \ker \psi \bigl| = \dfrac{\bigl| JY[\deg f]\bigl|}
{\bigl| \ker f^* \bigl|} \dfrac{\bigl| JZ[\deg g] \bigl|}{\bigl|
\ker g^* \bigl|} = \dfrac{\bigl| JZ[\deg h] \bigl| \bigl|
P(Y/Z)[\deg f] \bigl|}{\bigl| \ker f^* \bigl| \bigl|\ker g^*
\bigl|}  \, . $$
\end{itemize}
\end{prop}

\begin{proof}
By Theorem~\ref{prop1} we have isogenies
 $$ \alpha : JY \times P(X/Y) \to JX \  , \ \ (y, x) \to f^{*}y + x \ ,
 $$
 $$ \beta : JZ \times P(Y/Z) \to JY \  , \ \ (z, y) \to g^{*} z + y \ ,
 $$
and
 $$ \gamma : JZ \times P(X/Z) \to JX \  , \ \ (z, x) \to h^{*} z + x \
 $$
\noindent whose  kernels have sizes $\bigl| \ker \alpha \bigl| =
\dfrac{\bigl| JY[\deg f] \bigl|}{\bigl| \ker f^{*} \bigl|} \; $,
$\bigl| \ker \beta \bigl| = \dfrac{\bigl| JZ[\deg g]
\bigl|}{\bigl| \ker g^{*} \bigl|} \;$ and  $\bigl| \ker \gamma
\bigl| = \dfrac{\bigl| JZ[\deg h] \bigl|}{\bigl| \ker h^{*}
\bigl|} \,$ respectively.

 \vskip12pt

A short computation and connectedness show that the image of
$\spp$ is contained in $P(X/Z)$.

Furthermore if $(y, x) \in \ker \spp$, we have that $0 = f^* y
+x$. But then $0 = \Nm f(f^*y) + \Nm f(x) = (\deg f) \; y$; hence
$y \in P(Y/Z)[\deg \; f]$ and $x = -f^* y$; therefore $\ker \spp
\subseteq \{ (y, -f^*y) : y \in P(Y/Z)[\deg \; f] \}$, which
proves the first part of i).

>From the commutative diagram
 $$ \xymatrix{ JZ \times P(Y/Z) \times P(X/Y)
 \ar[rr]^{(\beta , \, \text{id}_{P(X/Y)})} \ar[d]_{(\text{id}_{JZ} , \, \spp )}
\ar[drr]^{\psi} & & JY \times P(X/Y) \ar[d]^{\alpha} \\
 JZ \times P(X/Z) \ar[rr]^{\gamma} & &  JX}
 $$
\noindent we obtain that $\psi$ is an isogeny and also that
 $$  \bigl| \ker \psi \bigl| =  \bigl| \ker \spp \bigl|   \bigl| \ker \gamma \bigl|
 = \bigl| \ker \alpha \bigl|  \bigl| \ker \beta \bigl| \ .
 $$

Therefore $\bigl| \ker \psi \bigl| = \bigl| \ker \spp \bigl|
 \dfrac{\bigl| JZ[\deg h] \bigl| }{\bigl| \ker h^* \bigl|} =
  \dfrac{\bigl| JY[\deg f] \bigl| }{\bigl| \ker f^* \bigl|} \dfrac{\bigl| JZ[\deg g]
\bigl| }{\bigl| \ker g^* \bigl|} $ and in order to finish the
proof we only need to compute $\bigl| \ker \spp \bigl| $.

Now the last equality implies that
 $$ \bigl| \ker \spp \bigl| =
\dfrac{\bigl| JY[\deg f]  \bigl| \, \bigl| JZ[\deg g] \bigl|
}{\bigl| JZ[\deg h] \bigl|}
 \dfrac{\bigl| \ker h^* \bigl|}{\bigl| \ker f^* \bigl| \, \bigl| \ker g^*
 \bigl|} \ .
 $$

But it is clear that
 $$
 \dfrac{\bigl| JY[\deg f]  \bigl| \, \bigl|
 JZ[\deg g] \bigl| }{\bigl| JZ[\deg h] \bigl|} =
 \bigl| P(Y/Z)[\deg f] \bigl|
 $$
and that
 $$\ker h^{*} = {g^{*}}^{-1} (g^*(JZ) \cap \ker f^*)
 $$
from where it follows that
 $$ \bigl| \ker h^{*} \bigl| =
 \bigl| \ker g^{*} \bigl| \bigl| g^*(JZ) \cap \ker f^* \bigl| \, .
 $$
The proof is now complete.
\end{proof}

\section{Galois covers}

Assume that the cover $f : \w C \to C$ is a Galois cover; i.e.,
that there exists a subgroup $G$ of  the automorphism group of $\w
C$ such that $C = \w C/G$ and $f$ is the canonical quotient map.
For any $g \in G$ we denote by the same symbol $g$ the
automorphism induced by $g$ on $\w J$, by $\langle g \rangle$ the
subgroup of $G$ generated by $g$,  and by ${\w J}^{\, H}$ the set
of fixed points of $H$ in $\w J$, for any subgroup $H$ of $G$. The
following results are immediate.

\begin{prop}
\label{prop:galois}
 Let $f : \w C \to C = \w C/G$ be a Galois cover
of degree $d$ and let $P = P(f)$ denote the corresponding Prym
variety.

Then
\begin{enumerate}
\item[i)] $f^*(J) = (\w J^G)^0$ and
${\w{J}}^{\, G} = f^*(J) + P_0$, where $P_0 = P \cap {\w{J}}^{\,
G} \subseteq P[d]  \,$.
\item[ii)] If we define $\Nm G : \w J \to \w J$ by $\Nm G = \displaystyle{\sum_{g
\in G}} \, g$, then $\Nm G = f^* \circ \Nm f$.
\item[iii)] $P = (\ker \Nm G )^0$.
\end{enumerate}
\end{prop}

\begin{coro}
\label{coro:double}
If $f : \w C \to C$ denotes a double cover
given by the involution $\sigma : \w C \to \w C$ and $P$ denotes
the corresponding Prym variety, then
 $$ {\w J}^{\, \langle \, \sigma \rangle} = f^{*}J + P[2]
  $$
\end{coro}

\begin{proof}

>From Proposition \ref{prop:galois} i) we know that
 ${\w J}^{\, \langle \, \sigma \rangle} =  f^* J + P_0$,
where $P_0 = {\w J}^{\, \langle \, \sigma \rangle} \cap P
\subseteq P[2]$.

Since $\sigma$ is an involution, it is clear that $P[2] \subseteq
{\w J}^{\, \langle \, \sigma \rangle}$ and the result follows.
\end{proof}

\begin{prop}
Let $f : \w C \to C = \w C/\langle \alpha \rangle $ be a cyclic
unramified cover.

Then ${\w J}^{\, \langle \alpha \rangle}$ is connected; i.e., ${\w
J}^{\, \langle \alpha \rangle} = f^* J$.
\end{prop}

\begin{proof}
Assume $f$ is of degree $d$; then ${\w J}^{\, \langle \alpha
\rangle} = \ker ( 1 - \alpha )$.

If ${\mathcal L}$ is in $\ker ( 1 - \alpha )$, then there is an
isomorphism $\phi : {\mathcal L} \to \alpha ({\mathcal L})$. But
then $\alpha^{d-1} (\phi) \circ \ldots \alpha(\phi) \circ \phi$ is
an isomorphism from $ {\mathcal L}$ to itself, and therefore equal
to a nonzero complex constant $c$; adjusting the constant, we
obtain an isomorphism $\phi_1 =
\displaystyle{\frac{1}{\sqrt[d]{c}}} \phi$ of $ {\mathcal L}$ of
order $d$.

Therefore $G$ acts on ${\mathcal L}$; i.e., $G$ is linearizable
and it follows from \cite{mumav} that then ${\mathcal L}$ is in
$f^* J$.
\end{proof}

The next result appears in \cite{s3}; its corollary will be very
useful later on.

\begin{lemma}
\label{lemma:norm}
 Let $f : \w C \to C = \w C/G$ be a Galois cover
of degree $d$, let $H$ be a subgroup of $G$ and denote by $h : \w
C \to \w C/H$ and $\varphi : \w C/H \to C$ the corresponding
covers, as in the following diagram.
$$ \xymatrix{
              & \w C \ar[dl]_{h} \ar[dd]^{f}  \\
 \w C/H  \ar[dr]^{\varphi} & \\
                          & C = \w C/G}
$$

Then   $f^{*}(\Nm \varphi (z)) = g_1 h^{*} z + g_2 h^{*} z +
\cdots + g_r h^{*} z$ for all $z \in J(\w C/H)$, where $\{g_1 ,
g_2 , \ldots , g_r \}$ is a complete set of representatives for
$G/H$.
\end{lemma}

\vskip12pt

\begin{coro}
\label{coro:norm}
 Under the hypothesis of Lemma \ref{lemma:norm},
 $$h^{*}(P(\varphi)) = \{ x \in {\w{J}}^{\, H} : \
 \displaystyle\sum_{i=1}^r g_i (x) = 0 \}^{\circ} \ . $$
\end{coro}

\begin{proof}
Let
 $$
 A = \{ x \in {\w{J}}^{\, H} : \  \displaystyle\sum_{i=1}^r g_i (x) = 0 \}^{\circ} \, .
 $$

If $z \in P(\varphi)$ then $\Nm f (h^{*}z) = |H| \Nm \varphi(z) =
0$; hence $h^{*}(P(\varphi)) \subseteq (\ker \Nm f)^{\circ} =
(\ker \Nm G)^{\circ}$ and $\sum_{g \in G} g(h^{*}z) = 0$.

But clearly $h^{*}(P(\varphi)) \subseteq {\w J}^{H}$, so the last
equation can be written as $\lvert H \rvert
\displaystyle\sum_{i=1}^r g_i (h^{*}z) = 0$, and we have proven
that $h^{*}(P(\varphi)) \subseteq A$.

Conversely, let $x \in A$; then
 $$
 \Nm G (x) = \displaystyle\sum_i g_i \left( \displaystyle\sum_{k \in H}
 k \right)(x) = |H| \displaystyle\sum_i g_i (x) = 0
 $$
and we obtain that $A \subseteq (\ker \Nm G)^{\circ} = (\ker \Nm
f)^{\circ}$.

Therefore $\Nm h (A) \subseteq (\ker \Nm \varphi)^{\circ} =
P(\varphi)$; then
 $$|H| A = \displaystyle\sum_{k \in H} k (A) = h^{*}(\Nm h(A))
  \subseteq h^{*}(P(\varphi))
 $$
and the result follows.
\end{proof}

\section{Degree two covers}
\label{section:deg2}

 Consider a degree two morphism of curves $f :
\widetilde{C} \to C$ with total ramification degree $\omega$.

Set $g = g_C$ and $P = P(\widetilde{C}/C)$. Then
$g_{\widetilde{C}} = 2g-1+\frac{\omega}{2}$, $\dim P =
g-1+\frac{\omega}{2}$, and $|P[2]| = 2^{2g-2+\omega}$.

In this section we describe $P[2]$ explicitly in terms of the
ramification of $f$.

By  Theorem~\ref{prop1}  $\bigl| {(\ker f^{*})}^{\bot}\bigl| =
\dfrac{\bigl| J[2] \bigl|}{\bigl| \ker f^{*} \bigl|}$
 and $f^*$ induces an isomorphism
 $$
 {(\ker f^{*})}^{\bot}/\ker f^* \cong f^*J \cap P \subseteq P[2] \,
 .
 $$

If $\omega = 0$, then $|\ker f^*| = 2$; hence $|{(\ker
f^{*})}^{\bot}/\ker f^*| = 2^{2g-2}$ and $f^*J \cap P = P[2]$.

If $\omega \neq 0$, then $|\ker f^*| = 0$; hence
 $f^*$ induces an isomorphism from
$J[2]$ to $f^*(J[2]) \cap P \subseteq P[2]$. Therefore $f^*(J[2])
\subseteq P[2]$ and  $|P[2]/f^*(J[2])| = 2^{\omega - 2}$.

In particular if $\omega = 2 \,$, then we have that $J[2]$ is
isomorphic to $f^*(J[2]) = P[2]$.

Also note that it follows from Lemma \ref{lem:kernm} and from
Proposition \ref{prop:galois} that if $\omega \neq 0$, then $\ker
(\Nm f) = \ker (1 + \sigma)$ is connected, where $\sigma$ is the
involution of $\widetilde{C}$ associated to the cover $f$.

The cases $\omega = 0$ and $\omega = 2$ of above are the \lq\lq
Classical  Pryms\rq\rq \, described in Mumford's beautiful paper
\cite{mumprym}.

We now consider $\omega > 2$. Let  $ P_1 , \ldots , P_{\omega} \in
\widetilde{C}$ be the ramification points of $f$, and for $i \in
\{ 1, \ldots, \omega \}$, let $Q_i  = f(P_i) \in C$.

For $j \in \{ 2, \ldots, \omega  \}$, choose $m_j \in \pic^{0}(C)$
such that $m_j^{\otimes 2} = {\mathcal O}_C (Q_1 -Q_j)$ and
consider
 $$
 {\mathcal F}_j =   {\mathcal O}_{\widetilde{C}} (P_j
 -P_1) \otimes f^*(m_j) \, .
 $$

Then each ${\mathcal F}_j$ is in $P[2]$ and we have the following
result.

\begin{theorem}
\label{prop:case2}
 Let $f : \w C \to C$ be a degree two cover with
total ramification degree $\omega > 2$.

Then for any ${\mathcal L} \in P[2]$ there exist $m \in J[2]$ and
unique $\nu_j \in \{ 0,1 \}$ such that
 $$ {\mathcal L} = {\mathcal F}_2^{\nu_2} \otimes
 \ldots {\mathcal F}_{\omega - 1}^{\nu_{\omega - 1}} \otimes f^* (m) \ .
 $$

In other words,
 $$ P[2] = f^*(J[2]) \oplus_{j=2}^{\omega-1} {\mathcal F}_j
 \ {\mathbb Z}/2{\mathbb Z} \ .
 $$
\end{theorem}

\begin{proof}
Let ${\mathcal D} \in \pic^{\frac{\omega}{2}}(C)$ be the divisor
class defining the cover $f : \widetilde{C} \to C$; that is,
 \begin{align*}
 {\mathcal D}^{\otimes 2} & = {\mathcal O}_C (Q_1 + \ldots + Q_{\omega})
 & \text{ and } & \
 f^* ({\mathcal D}) & = {\mathcal O}_{\widetilde{C}} (P_1 + \ldots + P_{\omega})
 \, .
 \end{align*}

Consider  the Abel-Jacobi map
\begin{align*}
 \widetilde{C}^{(\widetilde{g})} & \stackrel{\int}\longrightarrow \widetilde{J} \\
 D & \longmapsto {\mathcal O}_{\widetilde{C}}(D-{\widetilde{g}} P_1) \ .
\end{align*}

Then, given ${\mathcal L} \in P[2]$, $\sigma$ acts on $\int
^{-1}(\mathcal L)=|D|$; hence there exists $D \in |D|$ such that
$\sigma D= D$, and therefore $D$ must be of the form $D= f^*(E) +
\sum_{i=1}^r P_{l_i} \,$, for some effective divisor $E$ on $C$.

It follows that ${\mathcal L} = {\mathcal
O}_{\widetilde{C}}(\sum_{i=1}^r (P_{l_i}-P_1)) \otimes f^*(n) $,
for some $n$ in $\pic^0 (C)$. But  ${\mathcal L}^{\otimes 2} =
{\mathcal O}_{\widetilde{C}}$ and, since $f^*$ is injective, it
follows that $n^{\otimes 2} = {\mathcal O}_{C}(\sum_{i=1}^r (Q_1 -
Q_{l_i}))$.

Therefore $$ {\mathcal L} = {\mathcal F}_{l_1} \otimes \ldots
{\mathcal F}_{l_r} \otimes f^* (m) \, , $$ where $m = n \otimes
m_{l_1}^{-1} \otimes \ldots \otimes m_{l_r}^{-1} \in J[2]$.

The following result completes the proof.
\end{proof}

\begin{lemma}
\label{lemma:case2}
 Under the conditions of Theorem~\ref{prop:case2}, we
have
 $${\mathcal F}_{i_1} \otimes \ldots \otimes {\mathcal F}_{i_k} \in
 f^*(J[2]), \text{ with } \,  2 \leq i_1 < \ldots < i_k \leq \omega
 $$
if and only if
 $$
 \{ i_1 , \ldots ,i_k \} = \{ 2, \ldots ,\omega \} \, .
 $$
\end{lemma}

\begin{proof}
Let ${\mathcal D} \in \pic^{\frac{\omega}{2}}(C)$ be as in the
proof of Theorem~\ref{prop:case2} and assume $\{ i_1 , \ldots ,i_k
\} = \{ 2, \ldots ,\omega \}$.

Then
\begin{align*}
 {{\mathcal F}_2 \otimes \ldots
\otimes {\mathcal F}_{\omega}} &
= {\mathcal O}_{\widetilde{C}}(P_2 + \ldots + P_{\omega} -
(\omega - 1)P_1) \otimes f^*(m_2 \otimes \ldots \otimes m_{\omega}) \\
& = f^*({\mathcal D}(-{\frac{\omega}{2}} Q_1)
\otimes m_2 \otimes \ldots m_{\omega}) \, .
\end{align*}

But
 $$({\mathcal D}(-{\frac{\omega}{2}} Q_1) \otimes m_2 \otimes
 \ldots m_{\omega})^{\otimes 2} = {\mathcal O}_C
 $$
i.e.,
 $${\mathcal F}_2 \otimes \ldots \otimes {\mathcal F}_{\omega}
  \in f^*(J[2]) \, .
 $$

To show the other implication we may assume, without loss of
generality, that $\{ i_1, \ldots , i_k \} = \{ 2, 3, \ldots , r
\}$. That is, it suffices to show that if $ {\mathcal F}_2 \otimes
\ldots \otimes {\mathcal F}_r = f^*(n)$, with $n \in J[2]$, then
$r = \omega$.

Indeed, $ {\mathcal F}_2 \otimes \ldots \otimes {\mathcal F}_r =
f^*(n)$ with $n \in J[2]$, is equivalent to
 $$
 {\mathcal O}_{\widetilde{C}}(\sum _{i=2}^r P_i-(r-1)P_i) = f^*(n \otimes
 m_2^{-1} \otimes \ldots \otimes m_r^{-1}) \, .
 $$

Now consider the case $r$ even and define
 $$
 {\mathcal F} = n \otimes m_2^{-1} \otimes \ldots \otimes m_r^{-1} \otimes
 {\mathcal O}_C(\frac{r}{2} Q_1) \, .
 $$

Then ${\mathcal F}^{\otimes 2} = {\mathcal O}_C(\sum_{i=1}^r Q_i)$
and $f^*({\mathcal F}) = {\mathcal O}_{\widetilde{C}}(\sum_{i=1}^r
P_i) \,$. But this is possible by the definition of ${\mathcal D}$
only if $r = \omega$.

If $r$ were odd, then we would define
 $$
 {\mathcal F} = n \otimes m_2^{-1} \otimes \ldots \otimes
 m_r^{-1} \otimes {\mathcal O}_C(\frac{r-1}{2} Q_1)
 $$
and have that ${\mathcal F}^{\otimes 2} = {\mathcal
O}_C(\sum_{i=2}^r Q_i)$ and that $f^*({\mathcal F}) = {\mathcal
O}_{\widetilde{C}}(\sum_{i=2}^r P_i) \,$, which is impossible
since $r < \omega$.
\end{proof}

We now give a different set of generators for $P[2]/f^*(J[2])$,
which will be useful later on.

\begin{coro}
\label{coro:p2}
Under the conditions of Theorem~\ref{prop:case2},
choose $n_i$ and $s_j$ in $\pic^{0}(C)$ such that $n_i^{\otimes 2}
= {\mathcal O}_C (Q_{2i-1} -Q_{2i})$ and that $s_j^{\otimes 2} =
{\mathcal O}_C (Q_{2j} -Q_{2j+1})$. Also let
 $${\mathcal L}_i = {\mathcal O}_{\widetilde{C}}
  (P_{2i} -P_{2i-1}) \otimes f^*(n_i), \text{ for } 1 \leq i \leq
\omega/2
 $$
and let
 $${\mathcal G}_j = {\mathcal O}_{\widetilde{C}} (P_{2j+1}
-P_{2j}) \otimes f^*(s_j), \text{ for } 1 \leq j \leq \omega/2-1
 $$

\noindent respectively.

Then $$ P[2] = f^*(J[2]) \displaystyle{\oplus_{i=1}^{\omega/2-1}}
{\mathcal L}_i {\mathbb Z}/2{\mathbb Z}
\displaystyle{\oplus_{i=1}^{\omega/2-1}} {\mathcal G}_i {\mathbb
Z}/2{\mathbb Z} \ . $$
\end{coro}

\begin{proof}
Modulo $f^*(J[2])$ we may write
$ {\mathcal F}_{2i} = {\mathcal L}_1 \otimes
{\mathcal G}_1 \otimes \ldots \otimes
{\mathcal G}_{i-1} \otimes {\mathcal L}_i
$
and
$ {\mathcal F}_{2i+1} = {\mathcal L}_1 \otimes
{\mathcal G}_1 \otimes \ldots \otimes
{\mathcal G}_{i-1} \otimes {\mathcal L}_i
\otimes {\mathcal G}_i
$, for $1 \leq i \leq \omega/2-1$.
\end{proof}

\begin{rem}
\label{rem:rel}
Note that ${\mathcal L}_1 \otimes {\mathcal L}_2 \otimes \ldots {\mathcal
L}_{\omega/2} = {\mathcal F}_2 \otimes {\mathcal F}_3 \otimes \ldots {\mathcal
F}_{\omega} \in f^*(J[2])$.
\end{rem}

\vskip12pt

\section{Covers of degree three}
For degree three covers, the following results appear in \cite{s3}.

\begin{theorem}
\label{theo:degree3}
 Let $u : Z \to X$ be a non-Galois cover of degree three.

Denote by $g$ the genus of $X$, by $\alpha$  the number of total
ramification points and by $\beta$  the number of simple
ramification points.

Then there is a curve $W$ admitting the symmetric group on three
letters $S_3 = \langle \tau, \sigma : \tau^3 = \sigma^2 = (\tau
\sigma)^2 =1 \rangle$ as group of automorphisms and a commutative
diagram of curves and covers as follows
\begin{equation}
 \xymatrix{
&W\ar[dl]_{\ell}\ar[dd]_{\gamma}\ar[dr]^{\psi}&\\
Z\ar[dr]_{u}&&Y\ar[dl]^{v}\\ &X&} \label{dia:s3}
\end{equation}
\noindent where $Y = W/\langle \tau \rangle$, $Z = W/\langle
\sigma \rangle$ and $X = W/S_3$.

Furthermore,  there is an $S_3$-equivariant isogeny
\begin{equation}
 R : JX \times P(Y/X) \times P(Z/X) \times P(Z/X)  \longrightarrow JW
\label{bigiso}
\end{equation}
\noindent given by $R(x,y,z_1,z_2) = \gamma ^* x + \psi ^*y + \ell
^*z_1 + \tau \ell ^*z_2$ where the action of $S_3$ is trivial on the
first factor of the left hand side, the alternating action on the
second factor, and the unique irreducible complex representation of
degree two of $S_3$ on the product of the last two factors.

The cardinality of its kernel is given by $$ |\ker R \, | =
\begin{cases}
2^{2g -1} \cdot 3^{6g -3 + \alpha} &\text{if $\beta = 0$} \\
 2^{2g} \cdot 3^{6g -3 + \alpha + \beta} &\text{if $\beta \neq 0$} \ .
\end{cases}
$$
\end{theorem}

\begin{rem}
Observe that there are some natural restrictions on the data.
Since $\beta = |B(Y \to X)|$, we must have $\beta \equiv 0 \ (2)$.
Also, if $g_X = 0$, then $Y$ and $Z$ connected imply $\beta \geq
2$ and $2\alpha + beta \geq 4$.
\end{rem}

The calculation of $ |\ker R \, |$ mentioned above depends upon
the following lemma (c.f. \cite[p. 136]{s3}), similar in spirit to
the description of the points of order two in a Prym variety for a
degree two cover given in Section \ref{prop:case2}; we include it
here for the sake of completeness.

\begin{lemma}
\label{kerl}
Let $W$ denote a curve with an $S_3$ action and
associated Diagram~(\ref{dia:s3}).

Let
 $$
 L = \{ z \in P(Z/X)[3] : \ell^* z = \tau \ell^* z \}
 $$
for $\tau $ an element of order three in $S_3$, and let $\{ P_1
,\ldots , P_{\alpha } \} \subseteq Z$ denote the set of total
ramification points of $u : Z \to X$.

Also, if $\alpha$ is greater than or equal to two,  for each $j$
in $\{ 2, \ldots , \alpha\}$ choose ${\mathfrak m}_j \in
\pic^{\circ}(X)$ such that ${\mathfrak m}_j ^{\otimes 3} =
{\mathcal O}_X (u(P_1) - u(P_j))$ and consider the element of $L$
given by
 $${\mathfrak F}_j = {\mathcal O}_Z (P_j - P_1) \otimes
u^*({\mathfrak m}_j) \, .
 $$

Then $u^{*} : JX \to JZ$ is injective and

\begin{enumerate}
\item[i)] $u^*(JX[3]) \subseteq L$;
\item[ii)] for $\alpha$ equal to zero or one, $u^*(JX[3]) = L$. Thus
in this case  $L$ is isomorphic to $JX[3]$;
\item[iii)] for $\alpha$ greater than one,
$L = u^*(JX[3]) \bigoplus_{j=2}^{\alpha} {\mathfrak F}_j \, {\mathbb
Z}/3{\mathbb Z}$.
\end{enumerate}
\end{lemma}

\begin{rem} One cand find in  \cite{ri} the particular instance
of Theorem \ref{theo:degree3} corresponding to the case $v : Y \to
X = {\mathbb P}^1$ is a hyperelliptic cover and $\psi : W \to Y$
is a cyclic unramified cover of odd prime degree (in our notation,
this is equivalent to $\alpha = 0$ and $X = {\mathbb P}^1$).

There he obtains an $S_3-$equivariant isomorphism
 $$
 JZ \times JZ \to P(W/Y)
 $$

>From our point of view, this situation is part of the more general
$S_3-$equivariant isogeny $P(Z/X) \times P(Z/X) \to P(W/Y)$
described in \cite{s3};  together with the result that the group
$L$ is isomorphic to the kernel of the isogeny
 $$
 r : P(Z/X) \times P(Z/X) \to P(W/Y)
 $$
given by
 $$
 r(z_1 , z_2) = \ell^{*} z_1 + \tau \ell^{*} z_2
 $$
where $\tau$ is any element of order three in $S_3$, also in
\cite{s3}, and under the above assumptions, we also obtain that
the isogeny $r$ is an isomorphism in this case.

In \cite{recillas:jeg94} one can find the case when $u : Z \to X =
{\mathbb P}^1$ is a simple trigonal cover (in our notation, this
is equivalent to $\beta = 0$ and $X = {\mathbb P}^1$).

There, an $S_3-$equivariant isogeny is obtained
 $$
 JY \times JZ \times JZ \to JW \, ,
 $$
and statement ii) of the lemma is shown for the case $\alpha = 0$
(and $X = {\mathbb P}^1$).
\end{rem}

 \vskip12pt

\section{Covers of degree four}

In this section we study covers $f : X \to T$ of degree four.

There are five possibilities for $f$: if it is a Galois cover,
then it may be either cyclic or given by the action of the Klein
group on $X$; if $f$ is non Galois, then its corresponding Galois
cover may be given either by the action of the dihedral group of
order eight or by the alternating group on four letters or by the
symmetric group on four letters.

We will see that the first three cases give the bigonal
construction as a particular instance, whereas the last two imply
the trigonal construction.

For each possibility of the Galois cover $W \to T$ with Galois
group $G$, associated to $f : X \to T$, we will give a geometric
decomposition of $JW$. Such decomposition is a $G-$equivariant
isogeny between $JW$ and a product of appropriate Jacobians and
Pryms of intermediate covers.

Moreover, each of the complex irreducible representations of $G$
corresponds to the action of $G$ on precisely one of the factors
in the product.

We will also compute the kernel of each isogeny.

\subsection{The cyclic case}

If $X$ is a curve such that
 $$
 {\mathbb Z}/4{\mathbb Z} = \langle g: g^4 = 1 \rangle \subseteq \Aut(X)
 $$
then $T$ will denote the quotient $X/\langle g \rangle$ and $F$
will denote  the quotient $X/ \langle g^2 \rangle$; the
corresponding maps will be denoted by $a : X \to F$ and  $b : F
\to T$.

We then have the following diagram of curves and covers.
\begin{equation}
\xymatrix{
    X \ar[d]_{a}   \\
    F \ar[d]_{b}   \\
 T
 }
 \label{cyclic}
 \end{equation}

\begin{theorem}
Let $X$ be  a curve such that ${\mathbb Z}/4{\mathbb Z} = \langle
g : g^4 = 1 \rangle \subseteq \Aut(X)$ with associated
Diagram~(\ref{cyclic}).

Denote by  $\delta $ the number of fixed points of $g$ in $X$
(i.e., the number of total ramification points of the cover $X \to
T$), and by $\gamma $ the number of fixed points of $g^2$ in $X$
not fixed by $g$.

Then:

\begin{enumerate}
\item[i)] $g_F = 2g_T -1 + \frac{\delta}{2}$ and
$g_X = 2g_F-1+\frac{\gamma + \delta}{2} = 4g_T-3+\frac{\gamma +
3\delta}{2}$.

Furthermore, the cardinality of the ramification locus is: $\delta
+ \gamma$, for $a : X \to F$,  and $\delta$, for $b : F \to T$.

In particular, the signature type of $T$ is $(g_T ; \overbrace{2,
\ldots ,2}^{\gamma/2} , \overbrace{4, \ldots ,4}^{\delta})$.

\vskip12pt

\item[ii)] There is a ${\mathbb Z}/4{\mathbb Z}-$equivariant
isogeny
 $$ \phi_{{\mathbb Z}/4{\mathbb Z}} : JT \times P(F/T)
 \times P(X/F) \to JW
 $$
defined by
 $$
 \phi_{{\mathbb Z}/4{\mathbb Z}} (t, f, x) = a^*(b^*(t)) + a^*(f) + x
 $$
where the action of ${\mathbb Z}/4{\mathbb Z}$ is: the trivial
action on $JT$, the action given by the representation $g \to -1$
on $P(F/T)$, and, on $P(X/F)$, the sum of the other two
irreducible representations of ${\mathbb Z}/4{\mathbb Z}$ given by
$g \to I = \sqrt{-1}$ and $g \to -I$, respectively.

\vskip12pt

\item[iii)] The cardinality of the kernel of
$\phi_{{\mathbb Z}/4{\mathbb Z}}$ is given by $$ \bigl| \ker
\phi_{{\mathbb Z}/4{\mathbb Z}} \bigl| =
\begin{cases}
        2^{6g_T-2 + \delta} , &\text{if $\delta > 0$;}   \\
        2^{6g_T-3 } , &\text{if $\delta = 0$ and
        $\gamma > 0$;}   \\
        2^{6g_T-4 }, &\text{if $\delta = \gamma = 0$.}
\end{cases}
$$
\end{enumerate}
\end{theorem}

Before proving the theorem, we state an elementary remark which
will be useful in the sequel.

\begin{rem}
Let us observe that i) of the theorem imposes the conditions
 $$ \gamma \equiv \delta \equiv  0 \ (2) \, .
 $$

Moreover, if $g_T = 0$ we must have $\delta \geq 2$ in order that
the corresponding covers be connected.
\end{rem}

\begin{proof}
For the first statement, the cardinality of the ramification loci
is clear, thus the formulae for the genera follow from the
Riemann-Hurwitz formula.

For the second and third statements, note that by
Proposition~\ref{prop1} and Remark~\ref{rem:lange} we have that if
$\delta > 0$, all induced homomorphisms between Jacobians are
injective, that if $\delta = 0$ and $\gamma > 0$, then $a^*$ is
injective and $|\ker (b \circ a)^*| = |\ker b^*| = 2$, finally
that if $\delta = \gamma = 0$, then $|\ker a^*| = |\ker b^*| = 2$,
and then apply Proposition~\ref{prop:comp}.
\end{proof}

\subsection{The Klein case}

Let $X$ be a curve such that the Klein group
 $${\mathcal K} =
\langle \sigma , \tau : \sigma^2 = 1, \tau^2 = 1, (\sigma \tau)^2
= 1   \rangle
 $$
is contained in $\Aut(X)$.

We let $T$ denote the quotient $ X/{\mathcal K}$ and, for $k \in
{\mathcal K} $, $k \neq 1$, we will let $X_{k}$ denote the
quotient $X/ \langle k \rangle$; the corresponding quotient maps
will be denoted by  $a_k : X \to X_k$, \quad   $b_k : X_k \to T $,
\, and $f : X \to T$.

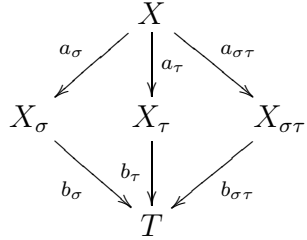
 \begin{equation}
 \xymatrix{
 &  X \ar[dl]_{a_{\sigma}} \ar[d]^{a_{\tau}}  \ar[dr]^{a_{\sigma \tau}} \\
 X_{\sigma} \ar[dr]_{b_{\sigma}} & X_{\tau}  \ar[d]_{b_{\tau}}
    & X_{\sigma \tau} \ar[dl]^{b_{\sigma \tau}}   \\
 & T
 }
 \label{klein}
 \end{equation}

\begin{theorem}
 \label{T:klein}
Let $X$ be a curve such that ${\mathcal K} \subseteq \Aut(X)$ with
associated Diagram~(\ref{klein}).

We also let $2r$ (resp. $2s$, $2t$) denote the number of fixed
points in $X$ of $\sigma \tau$ (resp. $\sigma$, $\tau$) (or
equivalently, the cardinality of the ramification locus of
$a_{\sigma \tau}$, $a_{\sigma }$, $a_{\tau }$, respectively).

Then:

\begin{enumerate}
\item[i)]
If $g$ denotes the genus of $T$, the genera of the intermediate
covers and the cardinality of the corresponding ramification loci
are given in the following table; in particular, the signature
type of $T$ is $(g ; \overbrace{2, \ldots ,2}^{s} , \overbrace{2,
\ldots ,2}^{t} , \overbrace{2, \ldots ,2}^{r})$.
 {\small{
\begin{center}
\begin{tabular}{|l|l|} \hline
\label{gt:klein}
 &  \\ genus & order of ramification $|B|$ \\ &
\\ \hline &      \\ $g_X = 4g-3 + s + t + r$  &   \\
       &                       \\ \hline
       &                       \\
$g_{X_{\sigma}} = 2g - 1 + \dfrac{r + t}{2}$
       & $|B(X \to X_{\sigma})| =  2s$\\
       &                       \\ \hline
       &                       \\
$g_{X_{\tau}} = 2g - 1 + \dfrac{s + r}{2}$
       & $|B(X \to X_{\tau})| =  2t$ \\
&   \\ \hline &   \\ $g_{X_{\sigma \tau}} = 2g - 1 + \dfrac{s +
t}{2}$
       & $|B(X \to X_{\sigma \tau})|  = 2r $ \\
&  \\ \hline &  \\
       & $|B(X_{\sigma} \to T)| =  t + r$ \\
       & $|B(X_{\tau} \to T)| =  s + r$   \\
       & $|B(X_{\sigma \tau} \to T)| =  s + t$\\
&  \\ \hline
\end{tabular}
\end{center}
}}
 \normalsize

\vskip12pt

\item[ii)] For $j,k,l \in \{ \sigma , \tau , \sigma \tau \}$ all different,
there are respective isogenies
 $$ \phi_{j} : P(X_{k}/T) \times P(X_{l}/T) \to P(X/X_{j})
 $$
defined by
 $$
 \phi_{j} (x_1, x_2) = a_{k}^*(x_1) + a_{l}^*(x_2) \, .
 $$

Furthermore, $\ker \phi_{j} \subseteq P(X_{k}/T)[2] \times
P(X_{l}/T)[2]$ and  its cardinality  is given as follows, where
$B_j = |B(X \to X_{j})|$.

\begin{equation*}
|\ker \phi_{j}| =
\begin{cases}
        2^{2g_T-2},    &\text{if $B_j = B_k = B_l = 0$;}   \\
        2^{2g_T-1},    &\text{if $B_j = 0$ and exactly
                         one of $B_k , B_l$ is zero;}    \\
        2^{2g_T},  &\text{if $B_j = 0$ and $ B_k B_l > 0$;} \\
        2^{2g_T-1+B_j/2},  &\text{if $B_j > 0$.}
\end{cases}
\end{equation*}

\item[iii)] There is an isogeny
 $$
 \varphi : P(X_{\sigma \tau}/T) \times P(X_{\sigma}/T) \times
 P(X_{\tau}/T) \to P(X/T)
 $$
defined by
 $$\varphi (x_1, x_2, x_3) = a_{\sigma \tau}^*(x_1) +
 a_{\sigma}^*(x_2) + a_{\tau}^*(x_3) \, .
 $$

Moreover, $\ker \varphi \subseteq P(X_{\sigma \tau}/T)[2] \times
P(X_{\sigma}/T)[2] \times P(X_{\tau}/T)[2] $ and its cardinality
is given by

\begin{equation*}
|\ker \varphi | =
\begin{cases}
        2^{4g_T-4}, &\text{if all covers are unramified;}   \\
        2^{4g_T-3+r+s+t}, &\text{otherwise.}
\end{cases}
\end{equation*}

\item[iv)] There is a ${\mathcal K}-$equivariant isogeny
 $$ \phi_{\mathcal K} : JT  \times P(X_{\sigma \tau}/T) \times
 P(X_{\sigma}/T) \times P(X_{\tau}/T) \to JX
 $$
defined by
 $$
 \phi_{\mathcal K} (t, x_1, x_2, x_3) = f^*(t) +
 a_{\sigma \tau}^*(x_1) + a_{\sigma}^*(x_2) + a_{\tau}^*(x_3)
 $$
where the action of ${\mathcal K}$ is:  the trivial one on $JT$,
the action of the irreducible representation of ${\mathcal K}$
given by $\sigma \to -1$ and $\tau \to -1$ on $P(X_{\sigma
\tau}/T)$, the action of the irreducible representation of $G$
given by $\sigma \to 1$ and $\tau \to -1$ on $ P(X_{\sigma}/T)$,
and, on $P(X_{\tau}/T)$, the action of the irreducible
representation of ${\mathcal K}$ given by $\sigma \to -1$ and
$\tau \to 1$.

The cardinality of the kernel of $\phi_{\mathcal K}$ is given by
\begin{equation*}
|\ker \phi_{\mathcal K}| =
\begin{cases}
        2^{8g_T-6}, &\text{if all covers are unramified;}   \\
        2^{8g_T-4+r+s+t}, &\text{if exactly two of $r,s,t$ are equal to zero;} \\
        2^{8g_T-3+r+s+t}, &\text{if exactly one or if none of $r,s,t$} \\
                          & \text{is equal to zero.}
\end{cases}
\end{equation*}
\end{enumerate}
\end{theorem}

\begin{rem}
The restrictions on the data in this case are
 $$
 r \equiv s \equiv t \equiv 0 \ (2)
 $$
and if $g_T = 0$ then $r+s$, $r+t$ and $s+t \geq 2$.
\end{rem}

\begin{proof}
For any involution $k \in {\mathcal K}$, denote $P_k = P(X_k/T)$.

Note that $f^*(t) \in JX^{{\mathcal K}}$, for any $t \in JT$, and
that if $k$ is any involution in ${\mathcal K}$, then $a_k^*(x)
\in JX^{\langle k \rangle}$, for any $x \in JX_k$.

Moreover if $h$ is any of the other two involutions in ${\mathcal
K}$, then $h$ induces the involution $\tilde{h}$ on $X_k$ which
gives the cover $b_k : X_k \to T$. By Corollary ~\ref{coro:norm}
we have that $P_k$ is the connected component of the identity of
$\ker (1 + \tilde{h})$. Therefore $h(a_k^*(x)) =
a_k^*(\tilde{h}(x)) = -a_k^*(x)$, for every $x \in P_k$. This
remark proves that the homomorphism $\phi_{\mathcal K}$ is
${\mathcal K}-$equivariant.

Consider the following homomorphisms
 $$
 \spp = \spp_{\sigma \tau} : P_{\sigma \tau} \times
 P(X/X_{\sigma  \tau}) \to P(X/T) \ , \, \spp (y, x) = a_{\sigma \tau}^*(y) + x
 $$
and
 $$
 \phi : JT \times P(X/T) \to JX  \ , \, \phi (t, x) = f^*(t) + x
 \, .
 $$

By Propositions~\ref{prop:comp} and \ref{prop1} we have that
 $\spp$ and $\phi$ are isogenies with kernels of cardinalities
 $$
 |\ker \spp | = |P_{\sigma \tau }[2]|
 \dfrac{| b_{\sigma \tau }^* JT \cap \ker a_{\sigma \tau }^*|}
 {|\ker a_{\sigma \tau }^*|}
 $$
and
 $$
 |\ker \phi \ | = \dfrac{|JT[4]|}{|\ker f^*|} =
 \dfrac{4^{2g_T}}{|\ker f^*|}
 $$
respectively.

By Remark~\ref{rem:lange} we have
 $$ |\ker f^*| = \begin{cases}
        4, &\text{if all covers are unramified;}   \\
        2, &\text{if exactly two of $r,s,t$ are equal to zero;} \\
          1, &\text{if exactly one or if none of $r,s,t$ is equal to zero.}
                 \end{cases}
$$

Now the homomorphism $\varphi$ may be factored as $\varphi = \spp
\circ \left( \id_{ P_{\sigma \tau}}, \phi_{\sigma \tau} \right)$,
and the homomorphism $\phi_{\mathcal K}$ as $\phi_{\mathcal K} =
\phi \circ \left( \id_{JT}, \varphi \right)$.

Hence if we show that $ \phi_{\sigma \tau} $ is an isogeny, it
will follow that $\varphi$ and $\phi_{\mathcal K}$ are isogenies.

Assuming this, it also follows that
 $$|\ker \varphi | = |\ker \spp | \cdot |\ker \phi_{\sigma \tau}|
 = |\ker \phi_{\sigma \tau}| |P_{\sigma \tau }[2]|
 \dfrac{| b_{\sigma \tau }^* JT \cap \ker a_{\sigma \tau}^*|}{|\ker a_{\sigma \tau }^*|}
 $$
and that
 $$ |\ker \phi_{\mathcal K}| = |\ker \phi | \cdot |\ker \varphi | =
 \dfrac{4^{2g_T}}{|\ker f^*|} |\ker \phi_{\sigma \tau}|
 |P_{\sigma \tau}[2]| \dfrac{| b_{\sigma \tau }^* JT \cap \ker a_{\sigma \tau}^*|}
 {|\ker a_{\sigma \tau }^*|} \ .
 $$

To complete the proof we therefore have to show that the
homomorphism $ \phi_{\sigma \tau} : P(X_{\sigma}/T) \times
P(X_{\tau}/T) \to JX$ is an isogeny onto $P(X/X_{\sigma \tau})$
and compute $|\ker \phi_{\sigma \tau} |$, $|P_{\sigma \tau }[2]|$
and $ \dfrac{| b_{\sigma \tau }^* JT \cap \ker a_{\sigma \tau
}^*|}{|\ker a_{\sigma \tau }^*|}$.

To show that $\phi_{\sigma \tau} (P(X_{\sigma}/T) \times
P(X_{\tau}/T))$ is contained in $P(X/X_{\sigma \tau})$ it suffices
to prove that $(1 + a_{\sigma \tau}) (\phi_{\sigma \tau}(x_1 ,
x_2)) = 0$ for all $(x_1 , x_2)$ in $P(X_{\sigma}/T) \times
P(X_{\tau}/T)$, since $P(X_{\sigma}/T) \times P(X_{\tau}/T)$ is
connected and since $P(X/X_{\sigma \tau}) = (\ker (1 + a_{\sigma
\tau}))^{\circ}$ by Corollary~\ref{coro:norm}.

But it follows from the remark at the beginning of the  proof of
this theorem that
 $$
 (1 + a_{\sigma \tau})(\phi_{\sigma \tau}(x_1 , x_2)) =
 (1 + a_{\sigma \tau})(a_{\sigma}^*(x_1) + a_{\tau}^*(x_2))= 0 \, .
 $$

Hence we have proven that $\phi_{\sigma \tau}$ is a homomorphism
from $P(X_{\sigma}/T) \times P(X_{\tau}/T)$ to $P(X/X_{\sigma
\tau})$. But these two varieties have the same dimension by i);
therefore to show that $\phi_{\sigma \tau}$ is an isogeny, it
suffices to show that its kernel is finite, which is what we prove
next.

The fact that if $(y,x) \in \ker \spp $, then $y \in P_{\sigma
\tau}[2]$ follows from the proof of Proposition~\ref{prop:comp}.

It follows that if $(x_1 , x_2 , x_3) \in \ker \varphi $, then we
have $x_1 \in P_{\sigma \tau}[2]$ by the factorization for
$\varphi$ given above. By symmetry we obtain the inclusion
\begin{equation}
\label{eq:cont} \ker \varphi \subseteq P_{\sigma \tau}[2] \times
P_{\sigma}[2] \times P_{\tau}[2] \ .
\end{equation}

Hence we may write
\begin{align*}
 \ker (\phi_{\sigma \tau})  &=  \{ (x_1 , x_2 ) \in
P_{\sigma}[2] \times P_{\tau}[2] : a_{\sigma}^* (x_1) = a_{\tau}^*
(x_2)    \} \\ & = (a_{\sigma}^* , a_{\tau}^*)^{-1} \{ (x, x) : x
\in a_{\sigma}^*(P_{\sigma}[2]) \cap a_{\tau}^*(P_{\tau}[2])  \}
\, \,
\end{align*}
which says it is finite, since
 $$
 \bigl| \ker (\phi_{\sigma \tau} ) \bigl|  =
 \deg a_{\sigma}^*\bigl|_{P_{\sigma}} \cdot
 \deg a_{\tau}^*\bigl|_{P_{\tau}} \cdot \bigl|
 a_{\sigma}^*(P_{\sigma}[2]) \cap a_{\tau}^*(P_{\tau}[2]) \bigl| \,
 \, .
 $$

The numbers appearing on the right side of the last equality and
also the number $\ker \spp_{\sigma \tau} =|P_{\sigma \tau }[2]|
\dfrac{| b_{\sigma \tau }^* JT \cap \ker a_{\sigma \tau
}^*|}{|\ker a_{\sigma \tau }^*|}$ are computed in the Appendix,
thus the proof is now complete.
\end{proof}

\begin{rem}
\label{rem:orient}
Given a double cover of a double cover $X \to Y
\to T$ which corresponds to a Galois four-fold cover, we have seen
that in both possible cases - a cyclic cover or a Klein group
action cover - there is another such object $X' \to Y' \to T$ (a
double cover of a double cover of $T$) which is naturally
associated to the original one: the same one for the cyclic case
and any of the other two in Diagram (\ref{klein}) for the Klein
case.

In other words, the cyclic and the Klein constructions together
give the bigonal construction (see \cite{d}) for the orientable
cover case.

In another section (see \ref{subsec:big}) we will complete the
bigonal construction to include the non-orientable case.
\end{rem}

\newpage

\section{The dihedral case}

Let $W$ be a curve such that the dihedral group of order eight
${\mathcal D}_4 = \langle r , s : r^4 = 1, s^2 = 1, (r s)^2 = 1
\rangle$ is contained in $\Aut(W)$.

We let $T$ denote the quotient $W/{\mathcal D}_4$, and for $d \in
{\mathcal D}_4$, $d \neq 1$, we will let $W_{d}$ denote the
quotient $W/ \langle d \rangle$.

Let $K_s$ (resp. $K_{rs}$) denote the Klein subgroup of ${\mathcal
D}_4$ generated by $r^2$ and $s$ (resp. $r^2$ and $rs$), and let
$W_{K_s}$ (resp. $W_{K_{rs}}$) denote the quotient $W/K_s$ (resp.
$W/K_{rs}$).

The corresponding quotient maps will be denoted by: $\gamma : W
\to T \,$; for $n \in \{ s, r^2 s, r^2 , rs, r^3 s \}$, $a_n : W
\to W_n \,$; for $n \in \{ s, r^2 s, r^2 \}$, $b_{n} : W_n \to
W_{K_s}$; for $n \in \{ r s, r^3 s, r^2 \}$, $c_{n} : W_n \to
W_{K_{rs}}$; for $n \in \{ K_s , K_{r s}, r^2 \}$, $d_{n} : W_n
\to T$; and $e : W_{r^2} \to W_r$.

Then we have the following diagram of curves and covers.

 \begin{equation}
 \xymatrix{
    & &  W\ar[dl] \ar[d]_{a_{r^2}} \ar[dr] \ar[dll]_{a_{r^2s}} \ar[drr]^{a_{r^3s}} & & \\
 W_{r^2s} \ar[dr] & W_s \ar[d]_{b_s} &  W_{r^2} \ar[dr] \ar[dl] \ar[d]_e &
    W_{rs} \ar[d]^{c_{rs}} & W_{r^3s}\ar[dl] \\
 &W_{K_s} \ar[dr]_{d_{K_s}}  & W_r \ar[d] & W_{K_{rs}}\ar[dl] \\
 & & T
 }
 \label{dihedral}
 \end{equation}

\begin{theorem}
\label{T:dih}
 Let $W$ be a curve such that ${\mathcal D}_4 \subseteq
\Aut(W)$ with associated Diagram~(\ref{dihedral}).

We let $2\delta$, $2\alpha$, $2\gamma_1$, $4\gamma_2$  denote the
number of fixed points in $W$ of $r$, $s$ (or $r^2 s$), $rs$ (or
$r^3s$), $r^2$ not fixed by $r$, respectively.

Then:

\begin{enumerate}
\item[i)]
If $g$ denotes the genus of $T$, the genera of the intermediate
covers and the cardinality of the corresponding ramification loci
are given in the following table.

In particular, the signature type of $T$ is
 $$(g ; \overbrace{4, \ldots ,4}^{\delta} , \overbrace{2, \ldots
 ,2}^{\alpha} , \overbrace{2, \ldots ,2}^{\gamma_1}, \overbrace{2,
 \ldots ,2}^{\gamma_2}) \ .
 $$

\small{
\hspace*{-2cm}\
\begin{tabular}{|l|l|} \hline
\label{gt:d4} &  \\ genus & order of ramification $|B|$ \\ &  \\
\hline &      \\ $g_W = 8g-7 + 2\alpha + 2\gamma_1 + 2\gamma_2 +
3\delta$  &   \\
       &                       \\ \hline
       &                       \\
$g_{W_s} =  4g-3+ \dfrac{\alpha+3\delta}{2}
            +\gamma_1+\gamma_2 =$ & $|B(W \to W_s)| =  2\alpha =$\\
$g_{W_{r^2 s}}$ &  $|B(W \to W_{r^2 s})| $  \\
       &                       \\ \hline
       &                       \\
$g_{W_{rs}} =  4g-3+ \dfrac{\gamma_1+3\delta}{2}
            +\alpha + \gamma_2$ = & $|B(W \to W_{rs})| = 2\gamma_1 =$ \\
$g_{W_{r^3 s}} $       & $|B(W \to W_{r^3 s})|    $ \\
       &   \\ \hline &   \\
$g_{W_{r^2}} = 4g-3 + \alpha + \gamma_1 + \delta$
       & $|B(W \to W_{r^2})|  = 4\gamma_2+2\delta  $ \\
&  \\ \hline &  \\ $g_{W_{K_s}} = 2g- 1 + \dfrac{\gamma_1 +
\delta}{2}$
       & $|B(W_s \to W_{K_s})| = \alpha + 2\gamma_2 + \delta$ \\
       & $|B(W_{r^2} \to W_{K_s})| = 2\alpha$   \\
&  \\ \hline &   \\ $g_{W_{K_{rs}}} = 2g-1 + \dfrac{\alpha +
\delta}{2}$
       & $|B(W_{rs} \to W_{K_{rs}})| = \gamma_1 + 2\gamma_2 + \delta$ \\
       & $|B(W_{r^2} \to W_{K_{rs}})| = 2\gamma_1$ \\
&  \\ \hline &  \\ $g_{W_r} = 2g - 1 + \dfrac{\alpha +
\gamma_1}{2}$
       & $|B(W_{r^2} \to W_r )| = 2\delta$ \\
&  \\ \hline & \\
       & $|B(W_{K_s} \to T)| = \gamma_1 + \delta$ \\
       & $|B(W_{r} \to T)| = \alpha + \gamma_1$ \\
       & $|B(W_{K_{rs}} \to T)| = \alpha + \delta$ \\
& \\ \hline
\end{tabular}
}

\normalsize

\vskip12pt

\item[ii)] There is a ${\mathcal D}_4-$equivariant
isogeny
 $$ \phi_{{\mathcal D}_4} : JT \times P(W_r/T) \times P(W_{K_s}/T)
 \times P(W_{K_{rs}}/T) \times 2P(W_s/W_{K_s}) \to JW
 $$
given by
 \begin{multline*}
 \phi_{{\mathcal D}_4} (t, w_1 , w_2 , w_3 , y_1 , y_2) )
 = \gamma^* (t) + (e \circ a_{r^2})^*(w_1) \\
 + (b_{s} \circ \,  a_s)^*(w_2) + (c_{rs} \circ \, a_{rs})^*(w_3)
 + a_s^*(y_1) + r a_{s}^*(y_2)
 \end{multline*}
where the action of ${\mathcal D}_4$ is: the trivial one on $JT$,
the action of the irreducible representation of ${\mathcal D}_4$
given by $r \to 1$ and $s \to -1$  on $ P(W_r/T) $, the action of
the irreducible representation of ${\mathcal D}_4$ given by $r \to
-1$ and $s \to 1$  on $ P(W_{K_s}/T) $, the action of the
irreducible representation of ${\mathcal D}_4$ given by $r \to -1$
and $s \to -1$ on $ P(W_{K_{rs}}/T) $, and, on $2P(W_s/W_{K_s})$,
the action of the unique irreducible representation of degree two
of ${\mathcal D}_4$.

The kernel of $\phi_{{\mathcal D}_4}$ has cardinality as follows.

\small{
\begin{center}
\begin{tabular}{|c|c|} \hline
\label{eq:d4}
&  \\ $\bigl| \ker \phi_{{\mathcal D}_4}   \,
\bigl|$ & Case \\
                   &  \\ \hline
                   &      \\
 $2^{20g-17}$    &  if $\alpha = \gamma_1 = \gamma_2 = \delta = 0$ \\
       &                       \\ \hline
      &  \\
 $2^{20g-15+ 4\gamma_1}$  & if $\alpha = \gamma_2 = \delta = 0$ and $\gamma_1 > 0$  \\
 &  \\ \hline
 &  \\
 $2^{20g-15}$  & if $\alpha = \gamma_1 = \delta = 0$ and $\gamma_2 > 0$  \\
 &  \\ \hline
 &  \\
 $2^{20g-14+ 3\alpha}$  & if $\gamma_1 = \gamma_2 = \delta = 0$ and $\alpha > 0$  \\
 &  \\ \hline
 &  \\
 $2^{20g-13+ 4\delta}$  & if $\alpha = \gamma_1 = 0$ and
     $\delta > 0$, for any $\gamma_2$  \\
 &  \\ \hline
 &  \\
 $2^{20g-13+ 4\gamma_1}$ & if $\alpha = \delta = 0$ and
   $\gamma_1 \cdot \gamma_2 > 0$  \\
 &  \\ \hline
 &  \\
 $2^{20g-12+ 4\delta + 4\gamma_1}$  & if $\alpha  = 0$ and
     $\delta \cdot \gamma_1 > 0$, for any $\gamma_2$  \\
 &  \\ \hline
 &  \\
 $2^{20g-13+ 3\alpha + 4\gamma_1 + 2\gamma_2}$  & if $\delta  = 0$,
     $\alpha  > 0$,   \\
 & and  $\gamma_1  = 0 ,  \gamma_2 >0$ or $\gamma_2 = 0, \gamma_1 > 0$\\
 &  \\ \hline
 &  \\
 $2^{20g-12+ 3\alpha + 4\gamma_1 + 2\gamma_2+5\delta}$  & if
     $\alpha \cdot \gamma_1 \cdot \gamma_2   > 0$, for any $\delta$   \\
 &  \\ \hline
 \end{tabular}
\end{center}
}

\vskip12pt

\item[iii)] There is a natural isogeny (the bigonal construction)
 $$
 \Nm a_{rs} \circ {a_s^*}_{|_{P(W_{s}/W_{K_{s}})}} : P(W_s/W_{K_s}) \to
 P(W_{rs}/W_{K_{rs}})\, .
 $$

If we denote by $G$ the Klein subgroup of $JW_{K_{rs}}[2]$ giving
the covers $c_n :  W_n \to W_{K_{rs}}$  for $n \in \{ r s, r^3 s,
r^2 \}$ for the case $\delta = \gamma_1 = \gamma_2 = 0$, then the
kernel of the isogeny has cardinality given as follows.

{\footnotesize{
\begin{equation*}
 \bigl|\ker \Nm a_{rs} \circ {a_s^*}_{|_{P(W_{s}/W_{K_{s}})}} \bigl| =
\begin{cases}
        2^{2g_T-5+2\delta}, &\text{ if } \delta >0
        \text{ and } \alpha = \gamma_1 = \gamma_2 = 0;   \\
        2^{2g_T-4+2\delta+\gamma_1+\gamma_2},
        & \text{ if } \delta >0 \text{ and } \\
        & \text{either } (\gamma_1 >0 \text{ and } \alpha = \gamma_2 = 0)   \\
        & \text{or } (\alpha > 0 \text{ and } \gamma_1 = \gamma_2 = 0) \\
        & \text{or } (\gamma_2 > 0 ,  \gamma_1 = 0 \text{ and any } \alpha); \\
        2^{2g_T-3+2\delta+\gamma_1+\gamma_2},
        & \text{ if } \gamma_1 > 0 \text{ and } \\
        & \text{either } (\delta = \alpha = 0 \text{ and any } \gamma_2)   \\
        & \text{or } (\delta\gamma_2 > 0 \text{ and any } \alpha) \\
        & \text{or } (\delta\alpha > 0 \text{ and } \gamma_2 = 0); \\
        2^{2g_T-2+\delta+\gamma_1+\gamma_2},
        & \text{ if } \delta = 0 \text{ and }  \\
        & \text{either } (\gamma_1\alpha > 0 \text{ and any } \gamma_2)   \\
        & \text{or } (\gamma_1 = \alpha = \gamma_2 = 0 \text{ and } \\
        & G \text{  isotropic}) \\
        & \text{or } (\gamma_1 = \alpha = 0 \text{ and } \gamma_2 >0); \\
        2^{2g_T-1+\gamma_2},
        & \text{ if } \delta = \gamma_1 = 0 \text{ and } \\
        & \text{either } (\gamma_2\alpha > 0)   \\
        & \text{or } (\gamma_2 = \alpha = 0 \text{ and } \\
        & G \text{  non isotropic}) \\
        & \text{or } (\gamma_2 = 0, \alpha > 0 \text{ and } \\
        & G \text{ isotropic}); \\
        2^{2g_T}, & \text{ if } \delta = \gamma_1 = \gamma_2 = 0, \\
        &   \alpha > 0 \text{ and } G \text{  non isotropic}.
\end{cases}
\end{equation*}
}}
\normalsize
\end{enumerate}
\end{theorem}

\begin{rem}
The conditions on the data in this case are the following.

 $$ \gamma_1 + \delta , \, \alpha + \delta, \, \alpha + \gamma_1 \equiv
 0 \ (2) \, .
 $$

Moreover, if $g = g_T = 0$ then $\gamma_1 + \delta$, $\alpha +
\delta$ and $\alpha + \gamma_1 \geq 2$.
\end{rem}

\begin{proof}
Consider the following two subdiagrams of (\ref{dihedral}):
\begin{equation}
\label{dia:fakeklein}
 \xymatrix{
     &  W_{r^2} \ar[dl]_{b_{r^2}} \ar[d]_{e}
      \ar[dr]^{c_{r^2}}
 \\
 W_{K_s} \ar[dr]_{d_{K_s}} & W_r  \ar[d]_{d_r}
    & W_{K_{rs}} \ar[dl]^{d_{K_{rs}}}
 \\
 & T
 } \ \
\ \
 \xymatrix{
    & &  W \ar[dll]_{a_{r^2 s}} \ar[dl]^{a_{s}}  \ar[d]^{a_{r^2}}
 \\
 W_{r^2 s} \ar[dr]_{b_{r^2 s}} & W_{s}  \ar[d]^{b_{s}} & W_{r^2}
 \ar[dl]^{b_{r^2}}
 \\
 & W_{K_s}
 }
\end{equation}

They correspond to actions of the Klein group, thus we can apply
Theorem~\ref{T:klein} to obtain isogenies
 $$\phi_{\mathcal K} :
  JT \times P(W_{r}/T) \times P(W_{K_s}/T) \times P(W_{K_{rs}}/T) \to
 JW_{r^2}
 $$
and
 $$\phi_{r^2} : P(W_s/W_{K_s}) \times P(W_{r^2 s}/W_{K_s}) \to
P(W/W_{r^2})
 $$
 defined by
 $$\phi_{\mathcal K} (t , w_1 , w_2 , w_3) =
 (d_{r} \circ e)^*(t) + e^{*}(w_1)
 + b_{r^2}^{*}(w_2) + c_{r^2}^{*}(w_3)
 $$
 and
 $$\phi_{r^2} (z_1 , z_2) = a_s^*(z_1) +
a_{r^2 s}^*(z_2)
 $$
respectively.

 Furthermore, we have that the cardinality of
their kernels is given by
\begin{equation}
\label{eq:dk} |\ker \phi_{\mathcal K}| =
\begin{cases}
        2^{8g_{T}-6}, &\text{if $\alpha = \delta = \gamma_1 = 0$ ;}   \\
        2^{8g_T-4+\alpha+\delta+\gamma_1}, &\text{if exactly two of $\alpha , \delta ,
        \gamma_1$ are  zero;} \\
        2^{8g_T-3+\alpha+\delta+\gamma_1}, &\text{if exactly one or none of $\alpha ,
        \delta , \gamma_1$ is zero;}
\end{cases}
\end{equation}
and
\begin{equation*}
|\ker \phi_{r^2}| =
\begin{cases}
        2^{16g_T-14+4\gamma_1} \, , &\text{if $\alpha = \delta =
        \gamma_2 = 0$ ;}   \\
        2^{16g_T-12+4\gamma_1+2\gamma_2+5\delta} \, , &\text{if }\alpha  = 0
            \text{ and } \gamma_2 + \delta >0   \\
        2^{16g_T-11+ 2\alpha + 4 \gamma_1 + 2 \gamma_2 + 5 \delta} \, ,
        &\text{ if } \alpha >0.
\end{cases}
\end{equation*}

Now since $r \, (r^2 s) = s \, r$,  the automorphism $r : JW \to
JW$ (induced by $r : W \to W$) induces an isomorphism $\tilde{r} :
JW_s \to JW_{r^2 s}$, therefore
 $$r_{|_{a_s^*(P(W_s/W_{K_s}))}} : a_s^*(P(W_s/W_{K_s}))
 \to a_{r^2s}^*(P(W_{r^2s}/W_{K_s}))
 $$
is an isomorphism.

In particular, it follows that $a_{r^2 s}^*(\tilde{r} (y)) = r
a_s^*(y)$, for any $y \in P(W_s/W_{K_s})$. By composing with
$\phi_{r^2}$ we obtain an isogeny
 $$ g : 2P(W_s/W_{K_s}) \to P(W/W_{r^2})
 $$
defined by
 $$g(y_1, y_2) = a_s^*(y_1) + r a_{s}^*(y_2)
 $$
whose kernel has cardinality given by
\begin{equation}
\label{eq:dg} |\ker g \, | =
\begin{cases}
        2^{16g_T-14+4\gamma_1} \, , &\text{if $\alpha = \delta =
        \gamma_2 = 0$ ;}   \\
        2^{16g_T-12+4\gamma_1+2\gamma_2+5\delta} \, , &\text{if }\alpha  = 0
            \text{ and } \gamma_2 + \delta >0   \\
        2^{16g_T-11+ 2\alpha + 4 \gamma_1 + 2 \gamma_2 + 5 \delta} \, ,
        &\text{ if } \alpha >0.
\end{cases}
\end{equation}

By Proposition~\ref{prop1} we also have an isogeny
 $$ \phi  : JW_{r^2} \times P(W_{r^2}/T) \to JW
 $$
given by
 $$\phi(y,x) = a_{r^2}^*(y) + x
 $$
whose kernel has cardinality
\begin{equation*}
\bigl| \ker \phi \bigl| = \dfrac{\bigl| JW_{r^2}[2] \bigl|}{\bigl|
\ker a_{r^2}^* \bigl|} = \dfrac{2^{8g_T -6 +\delta + \alpha_1 +
2\gamma_1}}{\bigl| \ker a_{r^2}^* \bigl|}.
\end{equation*}

By Remark~\ref{rem:lange} we know that $\bigl| \ker a_{r^2}^*
\bigl| = 2$, if $\delta = \gamma_2 = 0$, and that $\bigl| \ker
a_{r^2}^* \bigl| = 1$,  otherwise.

Therefore
\begin{equation}
\label{eq:dnat} \bigl| \ker \phi \bigl| =
\begin{cases}
    2^{8g_T -7 + \alpha + 2\gamma_1}, &\text{if $\delta = \gamma_2 = 0$;}   \\
    2^{8g_T -6 +\delta + \alpha + 2\gamma_1}, &\text{ if } \delta +
    \gamma_2 > 0.
\end{cases}
\end{equation}

Now note that
 $$\phi_{{\mathcal D}_4} = \phi \circ (\phi_{\mathcal K} , g)
 $$ and it follows that $\phi_{{\mathcal D}_4}$ is an isogeny
whose kernel has cardinality
 $$\bigl| \ker \phi_{{\mathcal D}_4} \bigl| = \bigl| \ker
 \phi \bigl| \cdot \bigl| \ker \phi_{\mathcal K} \bigl|
 \cdot \bigl| \ker g \bigl| \, .
 $$

Combining this equality with (\ref{eq:dnat}), (\ref{eq:dk}), and
(\ref{eq:dg}) we obtain (\ref{eq:d4}~ii).

Concerning the equivariance, note that the action of the
irreducible representation of order two of ${\mathcal D}_4$ is
given by $r(x,y) = (-y,x)$ and $s(x,y) = (x, -y)$.

For statement iii), we first prove that
 $\Nm a_{rs}(a_s^*(P(W_s/W_{K_s})))$ is contained in
$P(W_{rs}/W_{K_{rs}})$; then we show that
 $\Nm a_{rs} \circ a_s^*$ restricted to $P(W_s/W_{K_s})$ is an
isogeny and finally we compute the cardinality of its kernel.

Let us denote by
 $$ h = (\Nm a_{rs} \circ a_s^{*})|_{P(W_s/W_{K_s})}
 $$
the restriction of $\Nm a_{rs} \circ a_s^*$ to $P(W_s/W_{K_s})$.

Given $x \in P(W_s/W_{K_s})$ we know from the Klein case that
$a_s^*(x)$ is in $P(W/W_{r^2})$ and hence $\Nm a_{r^2}(a_s^*(x))
=0$.

Therefore
 \begin{multline*}
 \Nm c_{rs}(h(x)) = (\Nm c_{rs} \circ \Nm a_{rs})(a_s^*(x))) \\
 = (\Nm c_{r^2} \circ \Nm a_{r^2})(a_s^*(x))) = 0
 \end{multline*}
and it follows that $h(P(W_s/W_{K_s})) \subseteq
P(W_{rs}/W_{K_{rs}})$.

\medskip

To show that $h$ is an isogeny, denote by
 $$
 A = a_s^{*}(P(W_s/W_{K_s}))
 $$
and by
 $$
 B = a_{rs}^{*}(P(W_{rs}/W_{K_{rs}}))
 $$

 Then, by Corollary~\ref{coro:norm} iii), we know that
 $$A = \{  z \in JW^{\langle s \rangle} : z + r^2 s z = 0 \}^{\circ}
 $$
and that
 $$B = \{  w \in JW^{\langle r s \rangle} : w + r^2  w = 0
\}^{\circ} \, .
 $$

Moreover, the endomorphisms of $JW$ given by $1 + r s $ and $1 +
s$ induce morphisms which fit in the following commutative
diagram.

 {\tiny{
\begin{equation} \label{eq:mult2}
 \xymatrix{
      & A \ar[dr]_{\Nm a_{rs}} \ar[rr]^{1+rs} & & B \ar[dr]_{\Nm a_{s}} \ar[rr]^{1+s}& & A
 \\
 P(W_s/W_{K_s}) \ar[rr]_{h} \ar[ur]^{a_{s}^{*}} & & P(W_{rs}/W_{K_{rs}})
 \ar[ur]_{a_{rs}^{*}} \ar[rr]_{\Nm a_{s} \circ a_{rs}^*} & & P(W_s/W_{K_s})
 \ar[ur]_{a_{s}^{*}} }
\end{equation}}}
and a short computation shows that $(1 + s) \circ (1 + r s) = 2_A$
and $(1 + rs) \circ (1 + s) = 2_B$, respectively. Therefore $h$ is
an isogeny.

\medskip

We now compute the cardinality of the kernel $K$ of $h$ through
the next steps.
\begin{enumerate}
\item[1)]
Since
 $$| \ker \left( (1+rs) \circ a_s^{*}\right)_{|_{P(W_s/W_{K_s})}}
| = | \ker (a_{rs}^{*} \circ h) |  \, ,
 $$
we have that

\begin{equation}
\label{eq:kercard}
|K| = \frac{ |\ker \left( (1+rs) \circ a_s^{*}
\right)_{|_{P(W_s/W_{K_s})}} |}{| \ker
({a_{rs}^{*}}_{|_{P(W_{rs}/W_{K_{rs}})}}) |}
\end{equation}

Now from the rightmost diagram in (\ref{dia:fakeklein}) we compute
{\small{
 \begin{equation}
 \label{eq:ars}
  \bigl| \ker ({a_{rs}^{*}}_{|_{P(W_{rs}/W_{K_{rs}})}}) \bigl| =
\begin{cases}
    1, &\text{ if } \gamma_1 > 0;   \\
    2, &\text{ if } \gamma_1 = 0 \text{ and } \gamma_2 + \delta > 0; \\
    1, &\text{ if } \gamma_1 = \gamma_2 = \delta = 0 \text{ and $G$ is non-isotropic}; \\
    2, &\text{ if } \gamma_1 = \gamma_2 = \delta = 0 \text{ and if }
         G \text{ is isotropic}.
\end{cases}
\end{equation}}}

\normalsize

\noindent  if $G$ denotes the Klein subgroup of $JW_{K_{rs}}[2]$
giving the covers in the diagram for the case when they are all
unramified.

\item[2)] Since $a_s^{*}$ and $a_{rs}^{*}$ are isogenies and since $(1 + s)
\circ (1 + r s) = 2_A$ it follows that $\Nm a_{s} \circ
a_{rs}^* \circ h = 2_{P(W_s/W_{K_s})}$, so $\ker h \subseteq
P((W_s/W_{K_s})[2]$.

The same holds for
 $$ \Gamma = \ker ((1+r)\circ a_s^{*}|_{P(W_s/W_{K_s}}) \subseteq P(W_s/W_{K_s})[2]
 $$

In fact, we see immediately that
 $$ \Gamma = \{ x \in P(W_s/W_{K_s})[2] : a_s^{*}(x) \text{ is }
 D_4-\text{invariant} \}
 $$

\item[3)] We compute the cardinality of $\Gamma$ by decomposing
it into two parts: one coming from $JW_{W_{K_s}}$ and one coming
from the ramification of $b_s : W_s \to W_{K_s}$.

Therefore we define
 \begin{equation}
 \label{eq:ga1}
 \Gamma_1 = b_s^{*}(JW_{K_s}^{\langle \sigma \rangle}) \cap P(W_s/W_{K_s})[2]
 \end{equation}
where $\sigma : W_{K_s} \to W_{K_s}$ denotes the involution
induced by $rs$ and we will show that
 $$ \Gamma = \Gamma_1 \oplus \text{ part coming from the ramification of } b_s^{*}
 $$

In particular, if $b_s$ is unramified then it is clear that
$\Gamma = \Gamma_1$.

\item[4)] To compute $\Gamma_1$ we first apply Corollary
\ref{coro:double} to the cover \linebreak
 $d_{K_s} : W_{K_s} \to T$ to obtain
 $$ JW_{K_s}^{\langle \sigma \rangle} = {d_{K_s}}^{*}(JT) +
 P(W_{K_s}/T)[2]
 $$
and therefore
 \begin{equation}
 \label{eq:j2}
 JW_{K_s}^{\langle \sigma \rangle}[2] = ({d_{K_s}}^{*}(JT))[2] +
 P(W_{K_s}/T)[2]
 \end{equation}

But we also know (see \cite{mumprym} and Section
\ref{section:deg2}) that
 $$ P(W_{K_s}/T)[2] \subsetneq ({d_{K_s}}^{*}(JT))[2] \text{
 if } d_{K_s} \text{ is unramified}
 $$
and
 $$ ({d_{K_s}}^{*}(JT))[2] = {d_{K_s}}^{*}(JT[2]) \subseteq P(W_{K_s}/T)[2]
 \text{ if } d_{K_s} \text{ is ramified}
 $$

Combining this information with (\ref{eq:j2}) we obtain

\begin{equation}
\label{table:jac}
\begin{tabular}{|c|c|c|} \hline
   &  &   \\
  $JW_{K_s}^{\langle \sigma \rangle}[2]$
      & $\bigl| JW_{K_s}^{\langle \sigma \rangle}[2] \bigl|$ & \text{Case}  \\
   &  &   \\ \hline
   &  &   \\
  $({d_{K_s}}^{*}(JT))[2]$ & $2^{2g_T}$ & $\delta = \gamma_1 = 0$ \\
   &  &   \\ \hline
   &  &   \\
  $P(W_{K_s}/T)[2]$   & $2^{2g_T-2+\delta+\gamma_1}$ & $\delta + \gamma_1 > 0$ \\
   &  &   \\ \hline
\end{tabular}
\end{equation}

\vskip12pt

On the other hand we can now prove the following.

\begin{claim}
\label{claim:bs}
The following is always true.
 $$ {JW_{K_s}}^{\langle \sigma \rangle} \cap \ker b_s^{*} = \{ 0 \}
 $$
In particular,
 $$ P(W_{K_s}/T)[2] \cap \ker b_s^{*} = \{ 0 \} \, .
 $$
\end{claim}

Proof of the claim: if $b_s : W_s \to W_{K_s}$ is ramified, then
$b_s^{*} : JW_{K_s} \to JW_s$ is injective, by Proposition
\ref{prop:lange} and the claim follows.

If $b_s : W_s \to W_{K_s}$ is unramified then all maps on the
rightmost diagram of (\ref{dia:fakeklein}) are unramified.

Let $\ker {b_{s}}^{*} = \{ 0, \eta_{b_{s}} \} \subseteq
JW_{K_s}[2]$ the element defining the cover $b_s$. Also, let
$\sigma : W_{K_s} \to W_{K_s}$ denote the involution induced by
$rs : W \to W$.

>From the following commutative diagram $$
 \xymatrix{
  W \ar[r]^{rs} \ar[d]_{{a_{r^2s}}}   & W \ar[d]_{{a_{s}}} \\
  W_{r^2s} \ar[d]_{{b_{r^2 s}}} & W_s  \ar[d]_{{b_{s}}} \\
  W_{K_s} \ar[r]_{\sigma} & W_{K_s}
}
 $$
\noindent we see that $\sigma^{*}(\eta_{b_{s}})$ defines the cover
$b_{r^2 s} : W_{r^2 s} \to W_{K_s}$ and, in particular,
$\sigma^{*}(\eta_{b_{s}}) \neq \eta_{b_{s}}$; i.e., $\eta_{b_{s}}
\notin JW_{K_s}^{\langle \sigma \rangle}$, and the first part of
the claim follows in this case.

The last part follows from (\ref{eq:j2}), since it shows that
$P(W_{K_s}/T)[2] \subseteq {JW_{K_s}}^{\langle \sigma \rangle}$.
 \qed

\vskip12pt

We now continue the proof of the Theorem.

The claim just proved shows that $b_s^{*}({JW_{K_s}}^{\langle
\sigma \rangle})$ is isomorphic to ${JW_{K_s}}^{\langle \sigma
\rangle}$.

\vskip12pt

\item[5)] Computation of $\Gamma_1$, for the case $b_s$ ramified:
i.e., $\alpha + 2 \gamma_2 + \delta > 0$.

In this case it follows from \cite{mumprym} and Section
\ref{section:deg2} that $b_s^{*}$ injects $JW_{K_s}[2]$ into
$P(W_s/W_{K_s})[2]$ and from Claim \ref{claim:bs} that
 $$ b_s^{*}({JW_{K_s}}^{\langle \sigma \rangle})[2] =
 b_s^{*}({JW_{K_s}}^{\langle \sigma \rangle}[2]) \subseteq
 P(W_s/W_{K_s})[2]
 $$

Therefore
 $$ \Gamma_1 = b_s^{*}({JW_{K_s}}^{\langle \sigma \rangle}[2]) \cong
 {JW_{K_s}}^{\langle \sigma \rangle}[2]
 $$
and from (\ref{table:jac}) we obtain the following.

\begin{equation}
\label{table:jac2}
\begin{tabular}{|c|c|c|} \hline
     &  &   \\
 $\Gamma_1$ & $\bigl| \Gamma_1 \bigl|$ & \text{Case}  \\
   &  &   \\ \hline
   &  &   \\
 $b_s^{*}\left( ({d_{K_s}}^{*}(JT))[2] \right)$ & $2^{2g_T}$ & $\delta = \gamma_1 = 0$ and \\
   &  &   $\alpha + 2\gamma_2 >0$ \\
   &  &   \\ \hline
   &  &   \\
 $b_s^{*}\left( P(W_{K_s}/T)[2] \right)$&$2^{2g_T-2+\delta+\gamma_1}$&$\delta + \gamma_1 > 0$ and \\
   &  &   $\alpha + 2\gamma_2 + \delta >0$ \\
   &  &   \\ \hline
\end{tabular}
\end{equation}

\vskip12pt

\item[6)] Computation of $\Gamma_1$, for the case $b_s$ unramified:
i.e., $\alpha = \gamma_2 = \delta = 0$.

Note that in this case
 $$ \Gamma = \Gamma_1 \, .
 $$

It follows from \cite{mumprym} and Section \ref{section:deg2} that
if $\ker b_s^{*} = \{ 0, \eta_{b_s} \} \subseteq JW_{K_s}[2]$ then
 $$ \{ 0, \eta_{b_s}\}^{\perp}/\{ 0, \eta_{b_s}\} \cong
  b_s^{*}(\{ 0, \eta_{b_s}\}^{\perp}) = P(W_s/W_{K_s})[2]
 $$

Let us consider the map
 \begin{align*}
  \spp : P(W_{K_s}/T) \times P(W_s/W_{K_s}) & \to P(W_s/T) \\
 (x,y) & \to b_s^{*}(x) + y
 \end{align*}
whose kernel has order $\dfrac{1}{2} |P(W_{K_s}/T)[2]|$, from
Proposition \ref{prop:comp}.

It is clear that this kernel is isomorphic to
 $$ F= \{ x \in P(W_{K_s}/T)[2] : b_s^{*}(x) \in P(W_s/W_{K_s})[2] \}
 $$
via the map $F \to \ker \spp$ given by $x \to (x,-b_s^{*}(x))$.

We now subdivide into two cases, according to Table
\ref{table:jac}.

\item[6.1)] If $d_{K_s}$ is ramified, we know that
$JW_{K_s}^{\langle \sigma \rangle}[2] = P(W_{K_s}/T)[2]$.

Then (from Claim \ref{claim:bs})  we have
 $$ F \cong b_s^{*}(F) = b_s^{*}(JW_{K_s}^{\langle \sigma
 \rangle}[2]) \cap P(W_s/W_{K_s})[2] = \Gamma_1
 $$

Therefore we have proved that
 \begin{align}
 \label{eq:ga12} |\Gamma_1| & = 2^{2\dim P(W_s/T)-1} = 2^{2g_T-3+\gamma_1} \,
 , \\
 \notag & \ \ \ \ \ \text{ if } \alpha = \gamma_2=\delta = 0  \text{ and } \gamma_1 > 0
 \end{align}

\item[6.2)]  If $d_{K_s}$ is unramified, then $JW_{K_s}^{\langle \sigma
 \rangle}[2] = ({d_{K_s}}^{*}JT)[2]$.

Since $b_s^{*}(\{ 0, \eta_{b_s}\}^{\perp}) = P(W_s/W_{K_s})[2]$ we
can write
 $$ F = \{ 0, \eta_{b_s}\}^{\perp} \cap P(W_{K_s}/T)[2]
 \text{ of index two in }  P(W_{K_s}/T)[2] \, .
 $$

Therefore $P(W_{K_s}/T)[2]$ is not contained in $\{ 0,
\eta_{b_s}\}^{\perp}$; since $P(W_{K_s}/T)[2] \subseteq
({d_{K_s}}^{*}JT)[2]$ we obtain that
 $$ ({d_{K_s}}^{*}JT)[2] \not \subseteq \{ 0, \eta_{b_s}\}^{\perp}
 \, .
 $$

Moreover, $\{ 0, \eta_{b_s}\}^{\perp}$ is of index two in
$JW_{K_s}[2]$, and hence
 $$ ({d_{K_s}}^{*}JT)[2] \cap \{ 0, \eta_{b_s}\}^{\perp}
  \text{ is of index two in } ({d_{K_s}}^{*}JT)[2] \, .
 $$

But
 $$ b_s^{*}(({d_{K_s}}^{*}JT)[2] \cap \{ 0, \eta_{b_s}\}^{\perp})
 = b_s^{*}(({d_{K_s}}^{*}JT)[2]) \cap b_s^{*}(\{ 0,
 \eta_{b_s}\}^{\perp}) = \Gamma_1
 $$
where the first equality holds since $z = b_s^{*}(x) = b_s^{*}(y)$
with $x \in ({d_{K_s}}^{*}(JT))[2]$ and $y \in (\ker
b_s^{*})^{\perp}$ implies that $x = y + w$ with $w \in \ker
b_s^{*}$, which shows that $x \in (\ker b_s^{*})^{\perp}$ and
therefore $z \in b_s^{*}(({d_{K_s}}^{*}JT)[2] \cap \{ 0,
\eta_{b_s}\}^{\perp})$.

It follows that $|\Gamma_1| = \dfrac{1}{2} |({d_{K_s}}^{*}JT)[2]|
= \dfrac{1}{2} |JT[2]|$.

Therefore we have proved the following.
 \begin{equation}
 \label{eq:ga13}
  |\Gamma_1| = 2^{2g_T-1} \, ,
 \text{ if } \alpha = \gamma_2=\delta = \gamma_1 = 0
 \end{equation}

Putting together (\ref{eq:ga12}) and (\ref{eq:ga13}) we obtain
{\small{
\begin{equation}
\label{table:unram}
\begin{tabular}{|c|c|c|} \hline
    &  &   \\
 $\Gamma = \Gamma_1$ & $\bigl| \Gamma \bigl|$ & \text{Case}  \\
   &  &   \\ \hline
   &  &   \\
 $b_s^{*}(P(W_{K_s}/T)[2]) \cap P(W_s/W_{K_s})[2]$
   & $2^{2g_T-3+\gamma_1}$ & $\alpha = \gamma_2=\delta = 0$ \\
   & &  and $\gamma_1 > 0$ \\
   &  &   \\ \hline
   &  &   \\
 $b_s^{*}(({d_{K_s}}^{*}JT)[2] \cap \{ 0, \eta_{b_s}\}^{\perp})$
   &$2^{2g_T-1}$ & $\alpha = \gamma_2 = \delta = \gamma_1 = 0$  \\
   &  &   \\ \hline
\end{tabular}
\end{equation}
}}

\normalsize

\vskip12pt

\item[7)] We now complete the description of $\Gamma$ (in the case
$b_s$ ramified: $\alpha + 2\gamma_2 + \delta > 0$) by looking for
those elements coming from the ramification of $b_s$.

Denote by $\{ Q_1, \ldots , Q_{2\gamma_2} ,  \ldots ,
Q_{2\gamma_2+\delta+\alpha} \}$ the ramification points of $b_s$
ordered such that the corresponding branch points in $T$
 $$\{ d_{W_{K_s}}(b_s(Q_i)) = d_{W_{K_s}}(b_s(Q_{i+1})) \}_{i
 \in \{ 1,3, \ldots , 2\gamma_2 -1 \}} \text{ are of type }
 \gamma_2 \, ;
 $$
those of type $\delta$ are
 $$\{ d_{W_{K_s}}(b_s(Q_i)) \}_{2\gamma_2 + 1 \leq i \leq 2\gamma_2 +\delta} \, ; $$
and finally, those of type $\alpha$ are
 $$\{ d_{W_{K_s}}(b_s(Q_i)) \}_{2\gamma_2 + \delta +1 \leq i \leq 2\gamma_2 + \delta
 +\alpha} \, .
 $$

Now choose $m_i \in \pic^{0} W_{K_s}$ such that
 $$ m_i^{\otimes 2} = {\mathcal O}_{W_{K_s}}
  (b_s(Q_{i+1}) - b_s(Q_i))
 $$ and set
 $$ {\mathcal F}_i = {\mathcal O}_{W_{s}} (Q_i -Q_{i+1}) \otimes b_s^{*}(m_i) \, ,
 $$ for $1 \leq i \leq 2\gamma_2 + \delta + \alpha -1$.

Again, by Proposition~\ref{prop:case2} we obtain the description
{\small{
\begin{center}
\begin{tabular}{|c|c|} \hline
 & \\
$P(W_s/W_{K_s})[2]$ & Case \\
 & \\ \hline
 & \\
 $b_s^{*}(JW_{K_s}[2]) \oplus_{i=1}^{2\gamma_2 - 1}
 {\mathcal F}_i {\mathbb Z}/2{\mathbb Z}
 \oplus_{i=1}^{\delta + \alpha -1} {\mathcal F}_{2\gamma_2 + i}
 {\mathbb Z}/2{\mathbb Z}$
 &  $\gamma_2 \cdot (\delta + \alpha) > 0$ \\
  & \\ \hline
 & \\
 $b_s^{*}(JW_{K_s}[2]) \oplus_{i=1}^{2\gamma_2 - 2}
 {\mathcal F}_i {\mathbb Z}/2{\mathbb Z}$
 & $\delta = \alpha = 0$ \\
 & and $\gamma_2 > 0$ \\
   & \\ \hline
 & \\
 $b_s^{*}(JW_{K_s}[2])  \oplus_{i=1}^{\delta + \alpha -2} {\mathcal
  F}_{i} {\mathbb Z}/2{\mathbb Z}$
 &  $(\delta + \alpha) > 0$ \\
 & and $\gamma_2 = 0$ \\
   & \\ \hline
\end{tabular}
\end{center}}}

\normalsize

\vskip12pt

Noting that the  $a_s^{*}({\mathcal F}_i)$ are ${\mathcal
D}_4-$invariant precisely for $i \in \{ 1,3, \ldots , 2\gamma_2 -1
\}$ and for $ 2\gamma_2 +1 \leq i \leq 2\gamma_2 + \delta -1$, we
obtain that
 \small{$$ \Gamma =
 \begin{cases}
 \Gamma_1 \oplus_{i=1}^{\gamma_2} {\mathcal F}_{2i-1}
 {\mathbb Z}/2{\mathbb Z} \oplus_{i=1}^{\delta -1} {\mathcal F}_{2\gamma_2 + i}
 {\mathbb Z}/2{\mathbb Z} & \text{ if } \gamma_2 > 0 \text{ and } \delta  > 0; \\
 \Gamma_1 \oplus_{i=1}^{\gamma_2} {\mathcal F}_{2i-1} {\mathbb Z}/2{\mathbb Z}
 & \text{ if } \gamma_2 > 0, \alpha > 0 \text{ and } \delta = 0; \\
 \Gamma_1 \oplus_{i=1}^{\gamma_2-1} {\mathcal F}_{2i-1} {\mathbb Z}/2{\mathbb Z}
 & \text{ if } \gamma_2 > 0 \text{ and } \alpha = \delta = 0; \\
 \Gamma_1 \oplus_{i=1}^{\delta -1} {\mathcal F}_{2\gamma_2 + i}
 {\mathbb Z}/2{\mathbb Z} & \text{ if } \gamma_2 = 0 \text{ and }
 \alpha \cdot \delta  > 0; \\
 \Gamma_1  \oplus_{i=1}^{\delta -2} {\mathcal F}_{2\gamma_2 + i}
 {\mathbb Z}/2{\mathbb Z} & \text{ if } \gamma_2 = \alpha =0
 \text{ and } \delta  > 0; \\
 \Gamma_1  & \text{ if } \gamma_2 = \delta = 0 \text{ and } \alpha > 0 \, .
 \end{cases}
 $$
}

 Combining this expression with Table \ref{table:jac2}  we obtain
\begin{equation}
\label{eq:gammac}
 \bigl| \Gamma \bigl| =
\begin{cases}
2^{2g_T-4+2\delta+\gamma_1} & \text{ if } \alpha = \gamma_2 = 0
 \text{ and } \delta > 0; \\
 2^{2g_T-3+2\delta+\gamma_1+\gamma_2} & \text{ if }
 (\delta > 0  \text{ and } \alpha+\gamma_2 > 0) \\
 & \text{or } (\alpha=\delta = 0 \text{ and } \gamma_1\gamma_2 >0); \\
 2^{2g_T-2+\gamma_1+\gamma_2} & \text{ if }  \delta  = 0 \text{ and } \alpha\gamma_1 >
 0; \\
  2^{2g_T-1+\gamma_2} & \text{ if }  \delta  = \alpha = \gamma_1 = 0 \text{ and }
  \gamma_2 > 0; \\
   2^{2g_T+\gamma_2} & \text{ if }  \delta  = \gamma_1 = 0 \text{ and } \alpha
   > 0.
\end{cases}
 \end{equation}

\item[8)] The proof of Theorem \ref{T:dih} iii) is completed by
using (\ref{eq:kercard}) and putting together (\ref{eq:ars}),
(\ref{table:unram}) and (\ref{eq:gammac}).
\end{enumerate}
\end{proof}

\vskip12pt

\subsection{The bigonal construction}
\label{subsec:big}

Giving a curve $W$ with a ${\mathcal D}_4$ action and associated
Diagram~(\ref{dihedral}) is equivalent to giving a non-Galois
degree four cover  $f : X \to T$ which factorizes through two
covers of degree two $X \stackrel{\pi}{\to} Y \stackrel{g}{\to}
T$.

In the Diagram, $X$ corresponds to any one of the intermediate
\lq\lq first level\rq\rq quotients of $W$ except for $W_{r^2}$.
Without loss of generality, let $X$ correspond to $W_s$; then $Y$
corresponds to $W_{K_s}$.

In other words, we have associated to any nonorientable double
cover of a double cover, say given by
 $$
  W_s \stackrel{\pi}{\to} W_{K_s} \stackrel{g}{\to} T
 $$
another such double cover of a double cover,
given by
 $$
  X' = W_{rs} \to Y' = W_{K_{rs}} \to T \, .
 $$

That is, we have recovered the bigonal construction (see \cite{d})
for the nonorientable case, which together with Remark
\ref{rem:orient} gives the bigonal construction in general.

Furthermore, we have shown that $P(X/Y)$ is always isogenous to
$P(X'/Y')$ and we have described explicitly the kernel of the
isogeny in each case.

The case studied in \cite{p} is the bigonal construction applied
to a ramified double cover $K \to K_0$ of a hyperelliptic curve
$K_0$ to obtain the related cover $C \to C_0 \to {\mathbb P}^1$;
the main result there is that the two Prym $P(K/K_0)$ and
$P(C/C_0)$ in this special case are dual abelian varieties.

We will now obtain this result from our previous work on the
action of $D_4$.

Let $K_0$ be a hyperelliptic curve of genus $g$ and let $f : K_0
\to {\mathbb P}^1$ be the morphism given by the $g^1_2$ with
branch points $\{ a_1 , \ldots , a_{2g+2} \} \subseteq {\mathbb
P}^1$.

Consider a ramified cover of degree two $\pi : K \to K_0$, and
assume the ramification points $\{ p_1 , \ldots , p_{2h+2} \}$
satisfy the following condition
\begin{equation}
\label{cond:pantazis}
 \pi(p_i) + \pi(p_j) \notin g^1_2  \ \ \text{ for all } i,j \in
 \{ 1 , \ldots , 2h+2 \}  \ .
\end{equation}

Observe that $K$ has genus $2g+h$.

We claim that under the conditions just given, the Galois cover
associated to the fourfold cover $ F = f \circ \pi : K \to K_0 \to
{\mathbb P}^1$ has Galois group $D_4$, with $\delta = \gamma_2 =
0$ and $\alpha \cdot \gamma_1 > 0$ where $\alpha$, $\delta$,
$\gamma_1$ and $\gamma_2$ refer to the numbers in Theorem
\ref{T:dih}.

Proof of the claim: Let $b_i = F(p_i)$ for $1 \leq i \leq 2h+2$
and denote by $\Sigma = {\mathbb P}^1 - \{ a_1 , \ldots ,
a_{2g+2}, b_1 , \ldots , b_{2h+2} \}$ and by $\sigma_1 , \ldots ,
\sigma_{2g+2}, \tau_1 , \ldots , \tau_{2h+2}$ the corresponding
generators of $\Pi_1 (\Sigma)$.

Then the cover $F : K \to {\mathbb P}^1$ corresponds to a
representation
 $$
 F_1 : \Pi_1 (\Sigma) \to S_4
 $$

Let $G = F_1(\Pi_1(\Sigma))$. We know that  $G$ is a transitive
subgroup of $S_4$ and condition (\ref{cond:pantazis}) implies that
it is generated by transpositions $F_1(\tau_j)$ and products of
two disjoint transpositions $F_1(\sigma_k)$ for $1 \leq j \leq
2h+2$ and $1 \leq k \leq 2g+2$.

But furthermore, the associated diagram of cover maps
 $$\widetilde{K} \stackrel{\pi}{\to} \widetilde{K_0} \stackrel{f}{\to} \Sigma
 $$
where
 $$
 \widetilde{K} = K - \{F^{-1}(a_i) , F^{-1}(b_j) \}
 \text{ and }
 \widetilde{K_0} = K_0 - \{ f^{-1}(a_i) , f^{-1}(b_j) \}
 $$
corresponds to a chain of subgroups of index two
 $$
 S_3 \cap G \subseteq H \subseteq G
 $$
where $S_3$ is the subgroup of permutations of $S_4$ fixing the fourth
symbol say, and the correspondence is given by
 $$
 \Pi_1(\widetilde{K}) = F_1^{-1}(S_3 \cap G) \subseteq
\Pi_1(\widetilde{K_0}) = F_1^{-1}(H) \subseteq \Pi_1(\Sigma) =
F_1^{-1}(G) \, .
 $$

The existence of this chain of subgroups and the fact that $G$
contains transpositions and products of two disjoint
transpositions shows that $G$ must be isomorphic to $D_4$ and the
type for its generators proves that the numbers $\alpha$ ,
$\delta$ , $\gamma_1$ and $\gamma_2$ are as given.

In fact, recalling the notation of Theorem \ref{T:dih}, without
loss of generality we may assume that $K = W_s$ and $K_0 =
W_{K_s}$; it follows that then $C = W_{rs}$ and $C_0 =
W_{K_{rs}}$; furthermore, $\alpha = 2h+2$ and $\gamma_1 = 2g+2$
are positive and $\delta = \gamma_2 = 0$.

Under these conditions, observe that $a_s^{*} : JK \to JW$,
$a_{rs}^{*} : JC \to JW$ and $b_s^{*} : JK_0 \to JK$ are
injective.

Furthermore, part iii) of the Theorem shows that the isogeny $$
P_s = P(W_s/W_{K_s}) = P(K/K_0) \to P_{rs} = P(W_{rs}/W_{K_{rs}})
= P(C/C_0)
 $$
given by $H = \Nm a_{rs} \circ {a_s^*}_{|_{P_s}}$ has as kernel
$b_s^{*}(P(W_{K_s}/T)[2]) = \pi^{*}(JK_0[2])$, of cardinality
$2^{\gamma_1-2} = 2^{2g}$.

It also follows from Diagram (\ref{eq:mult2}) that the map
$\widetilde{H} \circ H : P_s \stackrel{H}{\to} P_{rs}
\stackrel{\widetilde{H}}{\to} P_s$ is multiplication by $2$ in
$P_s$, where $P_{rs} \stackrel{\widetilde{H}}{\to} P_s$ is given
by $\Nm a_{s} \circ {a_{rs}^*}_{|_{P_{rs}}}$ and therefore $\ker
(\widetilde{H} \circ H) = P_s [2]$.

But the polarization $\lambda_{P_s} : P_s \to \widehat{P_s}$
induced on $P_s$ by the principal polarization of $JK$ has as
kernel $\ker \lambda_{P_s} = \pi^{*}(JK_0) \cap P_s =
\pi^{*}(JK_0[2])$ by Proposition \ref{prop1}, and therefore $\ker
\lambda_{P_s} = \ker H$.

Similarly, $\ker \lambda_{P_{rs}} = \ker \widetilde{H} =
b_{rs}^{*}(P(W_{K_{rs}}/T)[2]) = b_{rs}^{*}(JC_0[2])$, of
cardinality $2^{\alpha-2} = 2^{2h}$.

Therefore,  there are respective isomorphisms $\Gamma :
\widehat{P_s} \to P_{rs}$ and $\Delta : \widehat{P_{rs}} \to
P_{s}$ making the following diagram commutative.
 $$ \xymatrix{
  P_s \ar[d]_{\lambda_{P_s}} \ar[r]^{H} &
    P_{rs} \ar[d]^{\lambda_{P_{rs}}} \ar[r]^{\widetilde{H}} & P_{s}     \\
  \widehat{P_s} \ar[ur]_{\Gamma}  &  \widehat{P_{rs}} \ar[ur]_{\Delta}
  }
 $$

If we denote by $\lambda$ the polarization on $\widehat{P_s}$
induced by $\lambda_{P_{rs}}$ via $\Gamma$, we may complete the
above to the following commutative diagram.
 $$ \xymatrix{
 P_s \ar[d]_{\lambda_{P_s}} \ar[r]^{H} &
  P_{rs} \ar[d]^{\lambda_{P_{rs}}} \ar[r]^{\widetilde{H}} & P_{s}         \\
 \widehat{P_s} \ar[ur]_{\Gamma} \ar[d]_{\lambda} &
  \widehat{P_{rs}} \ar[ur]_{\Delta} \ar[dl]^{\widehat{\Gamma}}        \\
  P_s
 }
$$

But then it follows that
 \begin{multline*}
 \ker (\lambda \circ \lambda_{P_s}) =
 \lambda_{P_s}^{-1}(\ker \lambda)  = H^{-1}(\ker \lambda_{P_{rs}}) \\
 = H^{-1}(\ker \widetilde{H})  = \ker (\widetilde{H} \circ H) =  P_s[2]
 \end{multline*}
and therefore
 $$
 \lambda = \lambda_{\widehat{P_s}}
 $$
showing that $(P_s , \lambda_{P_s})$ and $(P_{rs},
\lambda_{P_{rs}})$ are dual abelian varieties.

\section{The alternating case}

Let $W$ be a curve such that the alternating group on four letters
${\mathcal A}_4$ is contained in $\Aut(X)$.

Let  $\Delta = W/{\mathcal A}_4$ and $C = W /\langle \sigma
\rangle$ for $\sigma \in {\mathcal A}_4$ an element of order two.

Let $U = W/{\mathcal K}$ where  ${\mathcal K}$ denotes the Klein
subgroup of ${\mathcal A}_4$, and let $Y = W/\langle \tau \rangle$
where $\tau \in {\mathcal A}_4$ is an element of order three.

The corresponding cover maps will be denoted by $\gamma : W \to
\Delta$,  $ \varepsilon : W \to U$, $\nu : W \to C$, $c : C \to
U$, $\varphi : U \to \Delta$, $\psi : W \to Y$ and $h : Y \to
\Delta$.

Then we have the following diagram of curves and covers.

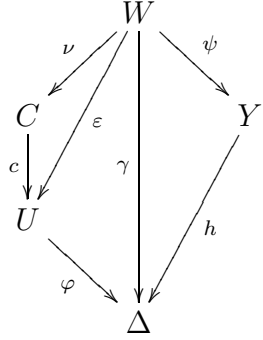
\begin{equation}
 \xymatrix{
  &              & W  \ar[ddl]^{\varepsilon} \ar[dl]_{\nu}  \ar[ddd]_{\gamma}  \ar[dr]^{\psi}\\
  & C \ar[d]_{c} &                   & Y \ar[ddl]^h  \\
  & U \ar[dr]_{\varphi} \\
  &              & \Delta
 }
 \label{dia:alt}
 \end{equation}

\begin{theorem}
\label{T:alt}
 If $W$ is a curve such that ${\mathcal A}_4
\subseteq \Aut(X)$ with associated Diagram~(\ref{dia:alt}), let
$\beta$ denote the number of fixed points of $\tau$, and let
$2\gamma_1$ denote the number of fixed points of $\sigma$.

Then:

\begin{enumerate}
\item[i)]
If $g$ denotes the genus of $\Delta$, the genera of the
intermediate covers and the cardinality of the corresponding
ramification loci are given in the following table. In particular,
the signature type of $\Delta$ is $(g ; \overbrace{3, \ldots
,3}^{\beta} , \overbrace{2, \ldots ,2}^{\gamma_1})$.

{\small{
\begin{center}
\begin{tabular}{|l|l|} \hline
\label{gt:a4} &  \\ genus & order of ramification $|B|$ \\ &  \\
\hline &      \\ $g_W = 12g - 11 + 4\beta + 3\gamma_1$  &   \\
       &                       \\ \hline
       &                       \\
$g_C = 6g - 5 + 2\beta + \gamma_1$ & $|B(W \to C)| =  2\gamma_1$\\
       &                       \\ \hline
       &                       \\
$g_Y = 4g - 3+ \beta + \gamma_1$ & $|B(W \to Y)| = 2\beta$ \\
       &                       \\ \hline
       &                       \\
$g_U = 3g - 2 + \beta$     & $|B(C \to U)|  = 2\gamma_1$ \\
       &                       \\ \hline
       &                       \\
                           &  $|B(U \to \Delta)|  = 2\beta$   \\
                           &  $|B(Y \to \Delta)|  = 2\beta + 2\gamma_1$   \\
                           &   \\ \hline
\end{tabular}
\end{center}
}}

\normalsize

\vskip12pt

\item[ii)] There is an ${\mathcal A}_4-$equivariant
isogeny
 $$ \phi_{{\mathcal A}_4} : J\Delta \times P(U/\Delta) \times 3P(C/U) \to JW
 $$
defined by
 $$\phi_{{\mathcal A}_4} (d, u , c_1 , c_2 , c_3 ) =
 \gamma^* (d) + \varepsilon^*(u) + \nu^*(c_1) + \tau \nu^*(c_2)
 + \tau^2 \nu^*(c_3)
 $$
where the action of ${\mathcal A}_4$ is: the trivial one on $ J\Delta $, the action of
the sum of the other two irreducible representations of degree one
of ${\mathcal A}_4$ on $ P(U/\Delta) $, and, on $3P(C/U)$, the
action of the irreducible representation of degree three of
${\mathcal A}_4$.

The cardinality of the kernel of $\phi_{{\mathcal A}_4}$ is given
by
\begin{equation*}
|\ker \phi_{{\mathcal A}_4}| =
\begin{cases}
        2^{24g-22} \, 3^{2g-1},
            &\text{if all covers are unramified;}   \\
        2^{24g-22+8\beta} \,  3^{2g},
            &\text{if $\gamma_1 = 0$ and $\beta$ is
            positive;} \\
        2^{24g-19+3\gamma_1} \,  3^{2g-1},
            &\text{if $\beta = 0$ and $\gamma_1$ is
            positive;}   \\
        2^{24g-19+8\beta+3\gamma_1} \,  3^{2g}, &\text{if
            $\beta \cdot \gamma_1$ is positive.}
\end{cases}
\end{equation*}

\vskip12pt

\item[iii)] There is a natural isogeny
 $$
 {\Nm \psi \circ \nu^{*}}_{|_{P(C/U)}} : P(C/U) \to P(Y/\Delta)
 $$
\noindent whose kernel is contained in $P(C/U)[2]$ and has cardinality given
as follows:

\begin{equation*}
|\ker ({\Nm \psi \circ \nu^{*}}_{|_{P(C/U)}})| =
\begin{cases}
        2^{4g-6+2\beta},
            &\text{if $\gamma_1 = 0$ and} \\
            & \text{if $P(U/\Delta)[2] \nsubseteq (\ker c^{*})^{\perp}$;}   \\
        2^{4g-5+2\beta+\gamma_1},
            &\text{otherwise.}
\end{cases}
\end{equation*}
\end{enumerate}
\end{theorem}

\begin{rem}
The conditions on the data in this case are the following: if
$g_{\Delta} = 0$ then $\beta \geq 2$ and $\beta + \gamma_1 \geq 3$
since $U$ and $Y$ must be connected covers.
\end{rem}

\begin{proof}

Statement i) is immediate.

For statement  ii), consider the commutative diagram
\begin{equation}
\label{dia:kalt}
 \xymatrix{
     &  W \ar[dl]_{\nu_2} \ar[d]_{\nu_3}  \ar[dr]^{\nu_4}
 \\
 C_2 \ar[dr]_{c_2} & C_3  \ar[d]_{c_3} & C_4 \ar[dl]^{c_4}
 \\
 & U
 }
\end{equation}
where $C_j = W/\langle (1 \, j)(k \, l) \rangle$ with
$\{ j,k,l\} = \{ 2,3,4 \}$.

Since it corresponds to an action of the Klein group, we may apply
Theorem~\ref{T:klein} to obtain an isogeny
 $$ \phi_{\mathcal K} : JU \times P(C_2/U) \times
 P(C_3/U) \times P(C_4/U) \to JW
 $$
given by
 $$
 \phi_{\mathcal K} (u_1, c_2, c_3, c_4) = \varepsilon^* (u_1) +
 \nu_2^*(c_2) + \nu_3^*(c_3) + \nu_4^*(c_4) \, .
 $$

Now we assume that $C = C_2$ and that $\tau$ is such that
 $$
 \tau (1 \, 2)(3 \, 4) = (1 \, 3)(2 \, 4)\tau \, .
 $$

Then, if we let $\tilde{\tau}$ denote the induced isomorphism from
$JC$ to $JC_3$, we obtain an isomorphism
 $$ n : 3P(C/U) \to P(C_2/U) \times P(C_3/U) \times P(C_4/U)
 $$
defined by
 $$ n(c_1 , c_2 , c_3) = (c_1,  \tilde{\tau} (c_2), \tilde{\tau}^2
 (c_3))\, .
 $$

We also have the natural isogeny
 \begin{align*}
 \phi_U : J\Delta \times P(U/\Delta) & \to JU \\
 (d,u) & \to \gamma^{*}(d) + u \, .
 \end{align*}

Hence we can write
 $$ \phi_{{\mathcal A}_4} (d, w , c_1 , c_2 , c_3)
   = \phi_{\mathcal K} (\phi_U (d,w), n(c_1 , c_2 , c_3))\, ,
 $$
and therefore $\phi_{{\mathcal A}_4}$ is an isogeny.

It also follows that its kernel has cardinality given by
 $$ \bigl| \ker
\phi_{{\mathcal A}_4} \bigl| = \bigl| \ker \phi_{\mathcal K}
\bigl| \cdot \bigl| \ker \phi_U \bigl| \, .
 $$

But
 $$
 \bigl| \ker \phi_{\mathcal K} \bigl| =
\begin{cases}
   2^{24g-22+8\beta}, &\text{if $\gamma_1 = 0$;} \\
   2^{24g-19+8\beta+3\gamma_1}, &\text{if $\gamma_1 > 0$,}
\end{cases}
 $$
by Theorem \ref{T:klein}, and
 $$
 \bigl| \ker \phi_U \bigl| =
 \begin{cases}
     3^{2g-1}, &\text{if $\beta = 0$;} \\
     3^{2g}, &\text{if $\beta > 0$}
 \end{cases}
 $$
by Remark \ref{rem:lange}.

Concerning iii), we first prove that $\Nm \psi (\nu^*(P(C/U)))$ is
contained in $P(Y/\Delta)$; then we show that  $\Nm \psi \circ
\nu^*$  restricted to $P(C/U)$ is an isogeny and finally we
compute the cardinality of its kernel.

Denote by $H$ the restriction of $\Nm \psi \circ \nu^*$ to
$P(C/U)$.

Given $x$ in $P(C/U)$ we obtain
 \begin{multline*}
 \Nm h (H(x)) = \Nm \varphi (\Nm c(\Nm \nu (\nu^* (x))))  \\
 =  2\Nm \varphi (\Nm c(x)) = 2\Nm \varphi (0) = 0
 \end{multline*}
and it follows that $H(P(C/U)) \subseteq P(Y/\Delta)$.

\medskip

We now show that $H$ is an isogeny. Let $A =\nu^{*}(P(C/U))$ and
$B = \psi^{*}(P(Y/\Delta))$.

By Corollary~\ref{coro:norm} iii) we have that
 $$ A = \{ z \in JW : z = \sigma(z) , z + \sigma'(z) = 0\}^{\circ}
 $$
where $\sigma'$ is any element of order two in
${\mathcal A}_4$ different from $\sigma$
 and that
 $$ B = \{ w \in JW : w = \tau (w) ,
 \displaystyle\sum_{k \in {\mathcal K}} k (w) = 0\}^{\circ} \, .
 $$

It is then clear that the endomorphisms of $JW$ given by $1 +\tau
+\tau^2$ and $1+\sigma$ induce respective isogenies $1 +\tau
+\tau^2 : A \to B$ and $1+\sigma : B  \to A$ such that
$(1+\sigma) \circ (1 +\tau +\tau^2) = 2_A$.

Now observe that  the commutative diagram
 \begin{equation}
\label{dia:isoa4}
 \xymatrix{
      & A \ar[dr]_{\Nm \psi} \ar[rr]^{1+\tau+\tau^2} & &
      B \ar[dr]_{\Nm \nu} \ar[rr]^{1+\sigma}  & & A  \\
  P(C/U) \ar[rr]_{H} \ar[ur]^{\nu^{*}} & & P(Y/\Delta) \ar[rr]_{\Nm \nu \circ \psi^{*}}
  \ar[ur]_{\psi^{*}} & & P(C/U) \ar[ur]_{\nu^{*}}
 }
\end{equation}
shows that $(\Nm \nu \circ \psi^{*}) \circ H = 2_{P(C/U)}$, since
$\nu^{*}$ is an isogeny.

In particular, $\ker H \subseteq P(C/U)[2]$.

Since it follows from $(1+\sigma) \circ (1 +\tau +\tau^2) = 2_A$
that $\ker (1 + \tau + \tau^2)_{|_{A}} \subseteq A[2]$ and since
$\deg \psi = 3$ it follows that any $z \in P(Y/\Delta) \cap \ker
\psi^{*}$ satisfies $2z = 0$ and $3z =0$, hence $z = 0$.

Hence $\ker \, \psi^{*}_{|_{P(Y/\Delta})} = \{ 0 \}$ and therefore
$\ker (1 + \tau + \tau^2)_{|_{A}} = \ker (\Nm \psi)_{|_{A}}$, so
the description of the kernel of the isogeny in iii) is as
follows.
 $$ \ker H = \{ z \in P(C/U)[2] : (1 + \tau  + \tau^2)
 (\nu^{*}(z)) = 0 \}.
 $$

We now compute the cardinality of this kernel.

Observe that all the cover maps in Diagram (\ref{dia:kalt}) have
ramification indices $2\gamma_1$, hence our calculation will be
divided into two cases: $\gamma_1 = 0$ and $\gamma_1 > 0$.

Since ${\mathcal K}$ is a normal subgroup of ${\mathcal A}_4$, the
action of $\tau $ on $W$ descends to an action on $U$, also
denoted by $\tau$; recall that $P(U/\Delta)$ is the connected
component of the identity of $\ker (1+\tau +\tau^2)$.

\begin{claim}
\label{claim:eps}
 If $\varepsilon : W \to U$ denotes the cover map
$c_2 \circ \nu_2$, then $\varepsilon \tau = \tau \varepsilon$.

We also have:

\begin{enumerate}
\item[i)] $P(U/\Delta)[2] = (\ker (1+\tau +\tau^2))[2]$,

\item[ii)] $\ker \varepsilon^{*} \subseteq (\ker (1+\tau +\tau^2))[2]$
\end{enumerate}

\end{claim}

Proof of the claim: The commutativity of $\varepsilon$ and $\tau$
is clear.

As for i), it follows from \cite[p. 61]{ri} that $\ker (1+\tau
+\tau^2) = P(U/\Delta) + \varphi^{*}J\Delta[3]$. Therefore, if $x
\in (\ker (1+\tau +\tau^2))[2]$ then $x = u +z$, with $u \in
P(U/\Delta)$ and $z \in \varphi^{*}J\Delta[3]$. But then $0 = 2x =
2u + 2z = 2u-z$, so $z \in P(U/\Delta)$ and hence $x \in
P(U/\Delta)$. This proves i).

To prove ii) observe that
 $$
 \ker \varepsilon^{*} =
 \begin{cases}
 0 & \text{if } \gamma_1 \neq 0 \\
 \{ 0, \eta_{c_2}, \eta_{c_3}, \eta_{c_4} = \eta_{c_2} +\eta_{c_3}\} & \text{if } \gamma_1 = 0
 \end{cases}
 $$
and that, in the second case, $\tau(\eta_{c_2}) = \eta_{c_3}$ and
$\tau(\eta_{c_3}) = \eta_{c_4}$.
 \qed

\vskip12pt

We continue the calculation of $\ker H$.

{\em Case I: $\gamma_1 = 0$.} In this case all the cover maps in
Diagram~(\ref{dia:kalt}) are \'etale, hence $\ker \varepsilon^{*}
= \{ 0, \eta_c , \eta_{c_3} , \eta_{c_4} \} \subseteq
P(U/\Delta)[2]$.

We also have $c^{*}(\{ 0 , \eta_c \}^{\perp}) = P(C/U)[2]$.

We further subdivide Case I into two subcases: Case I.a: $g (= g
_{\Delta}) = 0$ and Case I.b: $g > 0$.

{\em Case I.a: $g = 0$.} In this case $P(U/\Delta)[2] = JU[2]$ and
therefore $(1+\tau+\tau^2)m = 0$ for all $m \in JU[2]$ (by Claim
\ref{claim:eps}); in particular, for all $m \in \{ 0 , \eta_c
\}^{\perp}$.

But this implies that any $z \in P(C/U)[2]$ satisfies $ (1 + \tau
+ \tau^2) (\nu^{*}(z)) = 0$, since $P(C/U)[2] = c^{*}(\{ 0 ,
\eta_c \}^{\perp})$ and since $\varepsilon^{*} \tau = \tau
\varepsilon^{*}$.

We have thus proven that in this case $\ker H = P(C/U)[2]$ and
therefore
 $$ |\ker H |= 2^{2\beta-6} \ \text{ if } \gamma_1 = g = 0 \, .
 $$

{\em Case I.b: $g > 0$.} Now let $z \in \ker H$; that is, let $z$
be in $P(C/U)[2]$ such that $(1+\tau + \tau^2)(\nu^{*}(z)) = 0$.

Choose any $m \in \{ 0, \eta_c \}^{\perp}$ such that $c^{*}(m) =
z$. Then $m + \tau m + \tau^2 m$ belongs to $\ker \varepsilon^{*}
\subseteq P(U/\Delta)[2] $.

But $m + \tau m + \tau^2 m$ also belongs to $JU^{\langle \tau
\rangle} \cap P(U/\Delta) \subseteq P(U/\Delta)[3]$, by
Proposition \ref{prop:galois} i), and therefore $  m + \tau m +
\tau^2 m = 0$; that is, $m$ belongs to $P(U/\Delta)[2] \cap \{ 0,
\eta_c \}^{\perp}$.

But, since $ \{ 0, \eta_c \}^{\perp}$ is a subgroup of index two
of $JU[2]$, we have
 $$
 |P(U/\Delta)[2] \cap \{ 0, \eta_c \}^{\perp}|
=
\begin{cases}
       2^{2g_U - 2g_{\Delta}},
             & \text{if $P(U/\Delta)[2] \subseteq  \{ 0, \eta_c \}^{\perp}$;}   \\
       2^{2g_U - 2g_{\Delta}-1},
             & \text{otherwise;}
\end{cases}
 $$
\noindent and therefore we obtain
 $$ |\ker H| =
\begin{cases}
       2^{4g + 2\beta - 5},
  & \text{if $g > 0$, $\gamma_1 = 0$ and $P(U/\Delta)[2]\subseteq  \{ 0, \eta_c \}^{\perp}$;}  \\
       2^{4g + 2\beta - 6},
  & \text{if $g >0$, $\gamma_1 = 0$ and $P(U/\Delta)[2] \nsubseteq \{ 0, \eta_c \}^{\perp}$.}
\end{cases}
 $$

 \vskip12pt

 {\em Case II: $\gamma_1 > 0$.} In this case all of
the maps $\nu_j^{*}$ and $c_j^{*}$ are injective.

To compute the cardinality of the kernel we first describe
$P(C/U)[2]$.

We let $\{ Q_1 , \ldots , Q_{\gamma_1} , \sigma Q_1 , \ldots ,
\sigma Q_{\gamma_1} \}$ in $W$ denote the ramification points of
$\nu_3 : W \to C_3$, where $\sigma$ is the involution of $W$
giving the cover $\nu = \nu_2 : W \to C = C_2$.

Then the ramification points of $c : C \to U$ are $\{ P_1 =
\nu(Q_1) , P_2 = \nu(\tau(Q_1)), \ldots , P_{2\gamma_1-1} =
\nu(Q_{\gamma_1}) , P_{2\gamma_1} = \nu(\tau(Q_{\gamma_1})) \}$.

We now apply Corollary~\ref{coro:p2} to the cover $c : C \to U$
observing that ${\mathcal O}_U(\varepsilon Q_i - \varepsilon \tau
Q_i)$ belongs to $\ker (1+\tau+\tau^2) \subseteq JU$, hence for $i
\in \{ 1, \ldots , r \}$ we can choose $n_i \in \ker
(1+\tau+\tau^2)$ such that $n_i^{\otimes 2} = {\mathcal
O}_U(\varepsilon Q_i - \varepsilon \tau Q_i)$.

Now, setting ${\mathcal L}_i = {\mathcal O}_C (P_{2i}- P_{2i-1})
\otimes c^{*}(n_i) \in P(C/U)[2]$, we have that each ${\mathcal
L}_i$ clearly satisfies $(1+\tau+\tau^2)\nu^{*}{\mathcal L}_i =
0$, and since $ {\mathcal L}_1 \otimes {\mathcal L}_2 \ldots
\otimes {\mathcal L}_{\gamma_1} \in c^{*}JU[2] \cap \ker
(1+\tau+\tau^2)$, we obtain ${\mathcal L}_1 \otimes {\mathcal L}_2
\ldots \otimes {\mathcal L}_{\gamma_1} \in c^{*}P(U/\Delta)[2]$.

Letting ${\mathcal G}_1 , \ldots , {\mathcal G}_{\gamma_1-1}$ be
as in Corollary~\ref{coro:p2}, we obtain
 $$P(C/U)[2] = c^{*}JU[2]
\oplus_{i=1}^{\gamma_1-1} {\mathcal L}_i {\mathbb Z}/2{\mathbb Z}
\oplus_{i=1}^{\gamma_1-1} {\mathcal G}_i {\mathbb Z}/2{\mathbb Z}
\, .
 $$

To complete our computation, observe that $c^{*}JU[2] \cap  \ker H
= \{ z \in c^{*}JU[2] : (1+\tau+\tau^2)(\nu^{*}z) = 0 \}$ is
isomorphic (via $c^{*}$) to $JU[2] \cap \ker (1+\tau+\tau^2) =
P(U/\Delta)[2]$ and that no nontrivial combination of the
 ${\mathcal G}_i$ is in $\ker H$.

Therefore we have obtained the following result: if $\gamma_1
> 0$, then $ \ker H = c^{*}P(U/\Delta)[2]
 \oplus_{i=1}^{\gamma_1-1} {\mathcal L}_i {\mathbb Z}/2{\mathbb Z} $,
and  we have completed the proof.
\end{proof}

\subsection{The trigonal construction for the case ${\mathcal A}_4$}

As a corollary of Theorem \ref{T:alt} we obtain a particular case
of the trigonal construction as follows.

Let $h : Y \to {\mathbb P}^1$ be a tetragonal curve such that all
its branch points come from triple ramification points, with the
possible exception of one which is then of type $(2,2)$.

If we denote by $\gamma : W \to {\mathbb P}^1$ the corresponding
Galois extension, it has group ${\mathcal A}_4$ and we are in the
situation of Diagram (\ref{dia:alt}) with $\Delta = {\mathbb P}^1$
and $\gamma_1 = 0$ or $1$, respectively.

Moreover, $P(C/U)$ is a classical Prym.

We can  now prove the following.

\begin{coro}
Let $Y$ be a tetragonal curve such that all
its branch points come from triple ramification points, with the
possible exception of one which is then of type $(2,2)$.

Then the isogeny $H : P(C/U) \to JY$ from Theorem \ref{T:alt} iii)
induces an isomorphism between the principally polarized abelian
varieties $(JY, \lambda_{JY})$ and $(\wh{P(C/U)},
\lambda_{\wh{P(C/U)}})$.
\end{coro}

\begin{proof}
Let $P = P(C/U)$.

Under our hypothesis  the isogeny given in Theorem \ref{T:alt}
iii) by $H = \Nm \psi \circ \nu^{*} : P \to JY$ has as kernel
$P[2]$, which coincides with the kernel of $\lambda_P = 2 \lambda$
on $P$.

Therefore $H$ factorizes as follows, with $F$ an automorphism.
 $$
 \xymatrix{
  P \ar[d]_{\lambda_P = 2 \lambda} \ar[r]^(.6){H} & JY   \\
               \wh{P}  \ar[ru]_{F}^{\approx}           }
 $$

 We will now show that the isomorphism of complex tori
 $ F : \wh{P} \to JY$ is also an isomorphism of p.p.a.v.'s.

If we denote by $\lambda_1$ the polarization on $\wh P$ induced
via $F$ by $\lambda_{JY}$, we may complete the above diagram to
the following one.
 $$
 \xymatrix{
  P \ar[d]_{\lambda_{P}} \ar[rr]^{H} & &
              JY  \ar[d]^{\lambda_{JY}}  \\
  \wh{P} \ar[d]_{\lambda_1} \ar[urr]_{F}^{\thickapprox}& &
              \wh{JY} \ar[dll]^{\wh F}_{\thickapprox}   \\
  P
  }
 $$

It now follows from the commutativity of the above diagram that
$\lambda_1$ is principal and that $\ker (\lambda_1 \circ
\lambda_P) = P[2]$; therefore, $\lambda_1 = \lambda_{\wh P}$, as
claimed.
\end{proof}

\section{The symmetric case}

Let $W$ be a curve such that the symmetric group on four letters
${\mathcal S}_4$ is contained in $\Aut(W)$.

Let $T = W/{\mathcal S}_4$ and $\Delta = W/{\mathcal A}_4$.

For $\{ j, k, l \} = \{ 2, 3, 4 \}$ and $\sigma_j = (1 \, j)(k \,
l) \in {\mathcal S}_4$, let $C_j = W/\langle \sigma_j \rangle$ and
$U = W/{\mathcal K}$ where ${\mathcal K} = \{ 1,
\sigma_2,\sigma_3,\sigma_4 \}$ is the normal Klein subgroup of
${\mathcal S}_4$.

Also, let $Z_{kl} = W/\langle (k \, l) \rangle$, $S_j = W/\langle
\sigma_j , (k \, l) \rangle$ the quotient by a non normal Klein
subgroup, $V_j = W/\langle (1 \, k \, j \ l) \rangle$ the quotient
by a cyclic group of order four and $R_j = W/\langle (1 \, k \, j
\, l), (k \ l) \rangle$ the quotient by ${\mathcal D}_4$.

For $n \in \{ 1, 2, 3, 4 \}$ and $\{ j,k, l, n \} = \{ 1, 2, 3, 4
\}$, let $Y_n = W/\langle (j \ k \ l) \rangle$ be the quotient by
a cyclic group of order three, and let $X_n = W/\langle (j \ k \
l), (j \ k) \rangle$ denote the quotient of $W$ by the
corresponding ${\mathcal S}_3$.

Then we have a diagram of curves and covers as follows; the
sub-indices will be used as needed.
 
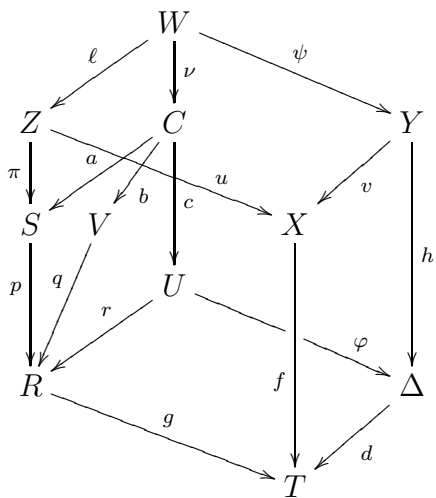
\begin{equation}
 \xymatrix@R=8pt@C=10pt
 { & & W \ar[ddrrrr]^{\psi} \ar[dd]^{\nu} \ar[ddll]_{\ell} \\
 &  \\
 Z \ar[dd]_{\pi} \ar[ddrrrr]^(.7){u} & & C \ar[ddd]^{c}
\ar[ddll]_{a} \ar[ddl]^(.6){b} & & & & Y \ar[ddll]^{v}
\ar[ddddd]^{h}  \\
   & &                 & &  \\
 S \ar[ddd]_(.4){p} & V \ar[lddd]_(.4){q} & & & X \ar[ddddd]_(.6){f}  \\
     & & U \ar[ddll]_(.4){r} \ar@{-}[drr]   \\
  & &         & &  \ar[drr]^{\varphi}  \\
 R  \ar[ddrrrr]^{g}   & & & & & & \Delta \ar[ddll]^{d}  \\
                     &                                  \\
                     & &                          &  & T
 }
 \label{dia:symm}
 \end{equation}

\begin{theorem}
\label{T:symm}
Let $W$ be a curve such that ${\mathcal S}_4
\subseteq \Aut(W)$ with associated Diagram~(\ref{dia:symm}).

We let $\gamma : W \to T = W/{\mathcal S}_4$ denote the quotient
map, and let $\tau$ denote any element of order three in
${\mathcal S}_4$.

Let $2\alpha$, $2\beta$,  $2\delta$, $4\gamma$ denote the number
of fixed points in $W$ of, respectively, any transposition, any
element of order three, any element of order four, the square of
any element of order four not fixed by the element of order four.

Then

\begin{enumerate}
\item[i)] If $g$ denotes the genus of $T$, the genera of the
intermediate covers and the cardinality of the corresponding
ramification loci are given in the following table. In particular,
the signature type of $T$ is $(g; \overbrace{2, \ldots
,2}^{\alpha}, \overbrace{3, \ldots ,3}^{\beta} , \overbrace{2,
\ldots ,2}^{\gamma} , \overbrace{4, \ldots ,4}^{\delta})$.

 {\small{
\hspace*{-1.5cm}
\begin{tabular}{|l|l|} \hline
\label{gt:s4} &  \\ genus & order of ramification $|B|$ \\ &  \\
\hline &      \\ $g_W = 24g-23+6\alpha+8\beta +6\gamma +9\delta$
&   \\
       &                       \\ \hline
       &                       \\
$g_C = 12g-11+3\alpha +4\beta+2\gamma+4\delta$ & $|B(W \to C)| =
                 4\gamma + 2\delta$  \\
&  \\ \hline &  \\ $g_Z = 12g-11 + \frac{5\alpha + 9\delta}{2}
+4\beta+3\gamma$ &
             $|B(W \to Z)| =  2\alpha $ \\
&   \\ \hline &   \\ $g_Y = 8g-7 + 2\alpha + 2 \beta + 2\gamma
+3\delta$ &  $|B(W \to Y)|
                 = 4\beta  $ \\
&  \\ \hline &  \\ $g_U = 6g- 5 + 3\frac{\alpha + \delta}{2} +
2\beta$   &
             $|B(C \to U)| = 4\gamma+2\delta$ \\
&  \\ \hline &   \\ $g_S = 6g-5 + \alpha + 2 \beta + \gamma
+2\delta$ &  $|B(Z \to S)|
                 = \alpha + 2\gamma + \delta  $ \\
      & $|B(C \to S)| = 2\alpha  $\\
&  \\ \hline &  \\ $g_X = 4g-3 + \frac{\alpha + 3\delta}{2} +
\beta + \gamma $ & $|B(Y \to X)| = 2\alpha  $ \\ &  $|B(Z \to X)|
= 2\alpha +2\beta$\\ &  \\ \hline & \\ $g_V = 6g-5 + 3\frac{\alpha
+ \delta}{2} + 2\beta + \gamma $ & $|B(C \to V)| = 2\delta  $ \\ &
\\ \hline & \\ &  $|B(U \to R)| = \alpha +\delta$  \\ $g_R = 3g-2
+ \frac{\alpha + \delta}{2} + \beta$ & $|B(S \to R)| = 2\gamma +
2\delta  $ \\ & $|B(V \to R)| = \alpha +2\gamma+\delta$\\ &  \\
\hline &  \\ $g_{\Delta} = 2g-1 + \frac{\alpha + \delta}{2} $ &
$|B(Y \to \Delta)| =  4\beta + 4\gamma +2\delta  $ \\ &  $|B(U \to
\Delta)| = 4\beta$\\ &  \\ \hline & \\ & $|B(R \to T)| = \alpha +
2\beta +\delta$\\ & $|B(X \to T)| = \alpha+ 2\beta +2\gamma
+3\delta$\\ & $|B(\Delta \to T)| = \alpha + \delta$\\ & \\ \hline
\end{tabular}
}}

\normalsize

\item[ii)] There is an ${\mathcal S}_4-$equivariant
isogeny
 $$ \phi_{{\mathcal S}_4} : JT \times P(\Delta/T) \times
 2P(R/T) \times 3P(S/R) \times 3P(V/R) \to JW
 $$
given by
 \begin{multline*}
 \phi_{{\mathcal S}_4} (t, d, r_1, r_2, s_1 , s_2 , s_3, v_1 , v_2 , v_3 )
 = \gamma^* (t) \\
 + (h \circ \psi)^*(d) + (p \circ \pi \circ \ell)^* (r_1) + \tau (p \circ \pi \circ \ell)^* (r_2) \\
 + (a \circ \nu )^* (s_1) + \tau (a \circ \nu )^* (s_2) + \tau^2 (a \circ \nu
 )^* (s_3) \\
 + (b \circ \nu)^* (v_1) + \tau (b \circ \nu)^* (v_2) + \tau^2 (b \circ \nu)^* (v_3)
 \end{multline*}
where the action of ${\mathcal S}_4$ is: the trivial one on $JT $,
the alternating action on $ P(\Delta/T) $, the action of the unique
irreducible representation of degree two of ${\mathcal S}_4$ on
$2P(R/T)$, the standard action (of degree three) of ${\mathcal S}_4$
on $3P(S/R)$, and, on $3P(V/R)$, the other irreducible action of
degree three of ${\mathcal S}_4$.

The cardinality of the kernel of $\phi_{{\mathcal S}_4}$ is given
by {\tiny{
\begin{equation*}
|\ker  \phi_{{\mathcal S}_4} \, | =
\begin{cases}
       2^{68g-65+22\beta} \cdot 3^{6g-3+\beta},
             & \text{if $\alpha = \gamma = \delta = 0$;}   \\
       2^{68g-61+15 \alpha + 22\beta} \cdot 3^{6g-3+\alpha +\beta},
             & \text{if $\gamma = \delta = 0$ and $\alpha >0$;}   \\
       2^{68g-59+22\beta+12\gamma} \cdot 3^{6g-3+\beta},
             & \text{if $\alpha = \delta = 0$ and $\gamma >0$;}   \\
       2^{68g-58+15 \alpha +22\beta+12\gamma +21\delta} \cdot
               3^{6g-3+\alpha +\beta+\delta},
             & \text{otherwise.}
\end{cases}
\end{equation*}
}}

\normalsize

\item[iii)] For any $C$ and $Y$  there is a natural isogeny
 $$
 {\Nm \psi \circ \nu^{*}}_{|_{P(C/U)}} : P(C/U) \to P(Y/\Delta)
 \, .
 $$

 Its kernel is contained in $P(C/U)[2]$ and has cardinality given as follows:
  {\small{
\begin{equation*}
|\ker  ({\Nm \psi \circ \nu^{*}}_{|_{P(C/U)}} ) \, | =
\begin{cases}
   2^{8g-10+4\beta+2\alpha} ,
    & \text{if $\gamma = \delta = 0$ and} \\
    &  \text{if $P(U/\Delta)[2] \nsubseteq (\ker c^{*})^{\perp}$;} \\
   2^{8g-9+4\beta+2\alpha +2\gamma+3\delta} ,
    & \text{otherwise.}
\end{cases}
\end{equation*}
}}

\normalsize

\item[iv)]
For any $Z$ and $C$ not covering the same $S$
there is a natural isogeny
 $$
 \Nm \nu \circ {\ell^{*}}_{|_{P(Z/S)}} : P(Z/S) \to P(C/U) \, .
 $$

Moreover, if for the case $\gamma = \delta = 0$ we denote by $G$
the Klein subgroup of $JU$ giving the covers $c_j : C_j \to U$ for
$j \in \{ 2,3,4 \}$, then the kernel of the isogeny  has
cardinality given as follows:

{\tiny{
 \begin{equation*} |\ker  (\Nm \nu \circ
{\ell^{*}}_{|_{P(Z/S)}} ) \, | =
\begin{cases}
       2^{6g_T-8+2\beta+4\delta},
             & \text{if $ \gamma  > 0$ and $\alpha = \gamma = 0$;}   \\
       2^{6g_T-7+\alpha+2\beta+3\gamma+4\delta},
             & \text{if either ($ \gamma > 0$ and $\alpha = \delta = 0$})   \\
             & \text{or ($\delta \alpha >0$ and
              $\gamma = 0$) or $ \delta \gamma  > 0$;} \\
       2^{6g_T-6+\alpha+2\beta+3\gamma},
             & \text{if either ($\alpha\gamma >0$ and $\delta = 0$)}  \\
             & \text{or ($\alpha = \delta = \gamma = 0$ and $G$ is isotropic);}   \\
       2^{6g_T-5+\alpha+2\beta},
             & \text{if either ($\alpha = \delta = \gamma = 0$ and} \\
             & \text{$G$ is non-isotropic)}  \\
             & \text{or ($\alpha >0$, $\delta = \gamma = 0$ and $G$ is isotropic);}   \\
       2^{6g_T-4+\alpha+2\beta},
             & \text{if $\alpha >0$, $\gamma = \delta = 0$ and} \\
             & \text{$G$ is non isotropic.}
\end{cases}
\end{equation*}
}}

\normalsize

\item[v)]
For any $S$ and $X$ covered by the same $Z$ there is a natural
isogeny (the trigonal construction)

 $$ {\Nm u \circ \pi^{*}}_{|_{P(S/R)}} : P(S/R) \to P(X/T)
 $$
whose kernel has cardinality as follows.

{\tiny{
\begin{equation*}
|\ker  (\Nm u \circ {\pi^{*}}_{|_{P(S/R)}} ) \, | =
\begin{cases}
       2^{2\dim P(S/R)},
             & \text{if $g=0$ and $ \gamma + \delta = 0$ or $1$;}   \\
       2^{4g-5+\alpha+2\beta-(\varepsilon+\zeta)},
             & \text{if $g > 0$ and $\gamma = \delta = 0$;}  \\
       2^{4g-4+\alpha+2\beta+\delta-\varepsilon},
             & \text{if $g > 0$ and $\gamma + \delta = 1$;}  \\
       2^{4g-5+\alpha+2\beta+\gamma+2\delta-\varepsilon},
             & \text{if $\gamma + \delta > 1$,}
\end{cases}
\end{equation*}
}} \normalsize
 where $\zeta = [P(R/T)[2]:P(R/T)[2]\cap \{ 0,
\eta_p \}^{\perp}] - 1$, $ \varepsilon = 0$ if $\alpha >0$,
$\varepsilon = 1$ if $\alpha = 0$.

\item[vi)] There is a natural isogeny
 \begin{align*}
 P(\Delta/T) \times P(V/R) & \longrightarrow P(Y/X) \\
 \intertext{given by}
 (d,v) & \longrightarrow h^{*}(d) + \Nm \psi \circ (b \circ
 \nu)^{*}(v) \, .
 \end{align*}

\item[vii)] There is a natural isogeny
 \begin{align*}
 P(R/T) \times P(Z/S) & \longrightarrow P(Z/X) \\
 \intertext{given by}
 (r,z) & \longrightarrow \pi^{*} \circ p^{*}(r) + z \, .
 \end{align*}

 $$
 P(Z/S) \times P(R/T) \to P(Z/X)
 $$
given by

\end{enumerate}
\end{theorem}

\begin{rem}
\label{rem:point}
 Equivalently, we could start with a degree four
cover  $f : X \to T$ whose Galois group is ${\mathcal S}_4$.

In this case, $\alpha$ is the number of simple ramification points
of $f$, $\beta$ is the number of ramification points of order
three, $\delta$ is the number of total ramification points, and
$2\gamma$ is the number of the remaining type of ramification
points.
\end{rem}

\begin{rem}
The conditions on the data are
 $$ \alpha + \delta \equiv 0 \ (2)
 $$
and if $g=0$ then $\alpha + \delta \geq 2$, $\beta \geq 1$ and
$\gamma+\delta \geq 1$.
\end{rem}

\begin{proof}
Using the Riemann-Hurwitz formula, statement i) is immediate.

For the proof of ii), we follow the idea behind the method of
\lq\lq little  groups\rq\rq \, and  apply it to the decomposition
of ${\mathcal S}_4$ given by the semi-direct product of the
abelian normal subgroup ${\mathcal K}$ and a subgroup ${\mathcal
S}_3$.

Parts  iii) and iv) are  immediate consequences of
Theorem~\ref{T:alt} iii) and  Theorem~\ref{T:dih} iii),
respectively.

For iii), just note that $\gamma_1$ for the case ${\mathcal A}_4$
corresponds to $2\gamma + \delta$ for the case ${\mathcal S}_4$
and $\beta$ for  the case ${\mathcal A}_4$ corresponds to $2\beta$
for the case ${\mathcal S}_4$.

As for case iv) note that $\alpha $ for the case $D_4$ coincides
with $\alpha$ for the case ${\mathcal S}_4$, as well as $\delta$,
that $\gamma_2$ for the case $D_4$ coincides with $\gamma$ for the
case ${\mathcal S}_4$ and that $\gamma_1$ for the case $D_4$
corresponds to $2\gamma + \delta$ for the case ${\mathcal S}_4$.

The isogenies in statements v) through vii) are suggested by
comparing statement ii) with the geometric decompositions of $JW$
obtained by respectively applying Theorems \ref{T:klein},
\ref{T:dih}, \ref{T:alt} and the results of \cite{s3} to the
actions of appropriate subgroups of ${\mathcal S}_4$ on $W$.

Concerning ii), by Proposition~\ref{prop1} we know that there is
an isogeny
 $$ \phi_U : JU \times P(W/U) \to JW \ , \, \phi_U (u,w) = (\nu \circ c)^*(u)+w
 $$
whose kernel has cardinality given by
 $$\bigl| \ker \phi_U \bigl| =
\begin{cases}
     2^{24g-22+6\alpha+8\beta}, &\text{if $\gamma=\delta = 0$;} \\
     2^{24g-20+6\alpha+8\beta+6\delta}, &\text{otherwise.}
\end{cases}
 $$

Note that this isogeny is ${\mathcal S}_4 -$equivariant with the
corresponding natural actions of ${\mathcal S}_4$ on each factor
on the left side.

We now decompose each such factor.

Since the action of ${\mathcal K} \subseteq {\mathcal S}_4$ is
trivial on $JU = J(W/{\mathcal K})$, there is a natural ${\mathcal
S}_4/{\mathcal K} = {\mathcal S}_3$ action on $JU$; under this
condition and from \cite{s3}, it follows that there is a natural
${\mathcal S}_3 -$equivariant isogeny
 $$ \phi_{{\mathcal S}_3} : JT \times P(\Delta/T) \times 2P(R/T) \to JU
 $$
defined by
 $$
 \phi_{{\mathcal S}_3}(t,d,r_1,r_2) = (r \circ g)^*(t)
 + \psi^*(d) + r^*(r_1) + \tilde{\tau}r^*(r_2)
 $$
where
 $\tilde{\tau}$ is the isomorphism induced on $JU$ by $\tau$ and
the action of ${\mathcal S}_3$ on the domain is given by: the
trivial representation of ${\mathcal S}_3$ on $JT$, the nontrivial
representation of degree one on $P(\Delta/T)$, and the unique
irreducible complex representation of degree two on $2P(R/T)$.

Note that each of these actions induces the corresponding
irreducible actions of ${\mathcal S}_4$ on each factor.

Also by \cite{s3},  the kernel of $\phi_{{\mathcal S}_3}$ has
cardinality given by
\begin{equation}
\label{eq:s3} |\ker  \phi_{{\mathcal S}_3} \, | =
\begin{cases}
       2^{2g-1} \cdot 3^{6g-3+\beta},
             & \text{if $\alpha = \delta = 0$;}   \\
       2^{2g} \cdot  3^{6g-3+ \alpha + \beta+\delta}, &\text{otherwise.}
\end{cases}
\end{equation}

\vspace{3mm}

Concerning the factor $P(W/U)$, there is a Klein-equivariant
isogeny $$ \nu^* + \tau\nu^* + {\tau}^2 \nu^* : 3P(C/U) \to P(W/U)
$$ whose kernel has cardinality given by
\begin{equation}
 \label{eq:kc}
 |\ker (  \nu^* + \tau \nu^* + {\tau}^2 \nu^*) \, |
  =
\begin{cases}
     2^{24g-24+6\alpha+8\beta},  & \text{if $\gamma = \delta = 0$,}   \\
     2^{24g-23+6\alpha+8\beta+6\gamma+9\delta} &\text{otherwise,}
\end{cases}
\end{equation}
as in the proof of Theorem~\ref{T:alt} and  by
Theorem~\ref{T:klein}.

Remark that this isogeny is, in fact, ${\mathcal S}_4 -$
equivariant.

To obtain the decomposition into irreducible representations, we
further decompose $P(C/U)$ as follows: Note that the following
piece of Diagram~(\ref{dia:symm}) corresponds to an action of a
Klein group on $C$ 
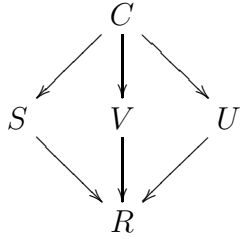
$$
 \xymatrix{
     &  C \ar[dl] \ar[d]  \ar[dr]  \\
 S \ar[dr] & V  \ar[d] & U \ar[dl] \\
 & R
 }
$$

By Theorem~\ref{T:klein} ii) we have a natural Klein-equivariant
isogeny defined by $$  a^* + b^*: P(S/R) \times P(V/R) \to P(C/U)
\, . $$

Furthermore, the cardinality of its kernel is given by
\begin{equation}
\label{eq:oklein} |\ker (a^* + b^*) | =
\begin{cases}
   2^{6g-6+2\beta}, &\text{if $\alpha =  \gamma = \delta = 0$;}   \\
   2^{6g-5+ \alpha + 2\beta + 2\gamma + 2\delta}, &\text{otherwise.}
\end{cases}
\end{equation}

Note that, then, the isogeny $$ 3P(S/R) \times 3P(V/R) \to P(W/U)
$$ given by the composition of the last two isogenies is
${\mathcal S}_4-$equivariant, where the action on $3P(S/R)$ is the
standard one and the action on $3P(V/R)$ is the other irreducible
action of degree three.

Combining the given isogenies we obtain part ii).

Concerning v), let $Z$ be the common cover of $S$ and $X$ via $\pi
:Z \to S$ and $u :Z \to X$ respectively.

We first prove that $\Nm u (\pi^*(P(S/R)))$ is contained in
$P(X/T)$; then we show that $\Nm u \circ \pi^*$ restricted to
$P(S/R)$ is an isogeny and finally we compute the cardinality of
its kernel.

Let $x \in P(S/R)$ and denote by $H$ the restriction of $\Nm u
\circ \pi^*$ to $P(S/R)$; we need to show that $H(x) \in P(X/T)$.
Since $P(S/R)$ is connected, it is enough to show that $\Nm f
(H(x)) = 0$. But
\begin{align*}
\Nm f (H(x)) & = \Nm g (\Nm p \, (\Nm \pi (\pi^* (x)))) \\
             & = \Nm g (\Nm p \, (2x)) \\
             & = 2\Nm g (0) \\
             & = 0
\end{align*}

\medskip

We now show that $H$ is an isogeny.

We may assume that $S = S_j = W/\langle (1 \ j), (k \ l) \rangle$
and that $X = W/\langle \tau , (1 \ j) \rangle$ where $\tau$
denotes an element of order three of $S_4$. Then $Z = Z_{1j}$ is a
common cover of $S$ and $X$, via $\pi : Z \to S$ and $u : Z \to
X$.

Let $A =\ell^{*}(\pi^{*}(P(S/R)))$ and let $B =
\ell^{*}(u^{*}(P(X/T)))$. Also, let $\sigma = (1 \ k)(j \ l)$
denote an element of ${\mathcal K} \subseteq S_4$ which induces
the involution $S \to S$ giving the cover $S \to R$.

By Corollary~\ref{coro:norm} iii) we have that
 $$
 A = \{ z \in JW^{\langle (1 \ j), (k \ l) \rangle} : z + \sigma z =
   0\}^{\circ}
 $$
and that
 $$
 B = \{ w \in JW^{\langle \tau , (1 \ j) \rangle} :
 \displaystyle\sum_{k \in {\mathcal K}} k (w) =
 0\}^{\circ} \, .
 $$

It is then clear that the endomorphisms of $JW$ given by $1 +\tau
+\tau^2$ and $1+ (k \ l)$ induce respective isogenies $1 +\tau
+\tau^2 : A \to B$ and $1+(k \ l) : B  \to A$ such that $(1+(k \
l)) \circ (1 +\tau +\tau^2) = 2_A$.

We also have the following commutative diagram
{\small
\begin{equation}
\label{dia:isotrig}
 \xymatrix{
      & A  \ar[rr]^{1+\tau+\tau^2} & &
      B  \ar[rr]^{1+(k \ l)}& & A  \\
  \pi^{*}(P(S/R)) \ar[drr]^{\Nm u} \ar[rr]^{u^{*}\circ \Nm u}
  \ar[ur]^{\ell^{*}} & & u^{*}(P(X/T))
  \ar[drr]^{\Nm \pi} \ar[rr]^{\pi^{*} \circ  \Nm \pi }
  \ar[ur]^{\ell^{*}} & &  \pi^{*}(P(S/R))
  \ar[ur]^{\ell^{*}} \\
  P(S/R) \ar[u]^{\pi^{*}}_{\cong} \ar[rr]_{H} & & P(X/T) \ar[u]^{u^{*}}_{\cong}
  \ar[rr]_{\Nm \pi \circ u^{*}} & & P(S/R) \ar[u]^{\pi^{*}}_{\cong}
 }
\end{equation}}
\normalsize
which shows that $H$ is an isogeny.

We are interested in computing the kernel of this isogeny.

For this, first note that it follows from Claim \ref{claim:bs} in
the proof of Theorem \ref{T:dih} iii) that $\pi^{*}$ restricted to
$P(S/R)$ is injective; therefore, the two external vertical arrows
of Diagram \ref{dia:isotrig} are isomorphisms. That the middle
vertical line is also an isomorphism follows from Lemma \ref{kerl}
(see also \cite[p. 136]{s3}).

We can now prove the following result, which will be fundamental
to complete the proof of the Theorem.
\begin{claim}
\label{claim:mult2}
The composition $$ \xymatrix@C=0.5cm{
  P(S/R) \ar[rr]^{H} & & P(X/T) \ar[rr]^{\Nm \pi \circ u^{*}} & & P(S/R)  }
  $$
is multiplication by $2$.
\end{claim}

Proof of the claim: Since the topmost line of Diagram
(\ref{dia:isotrig}) is multiplication by $2$ on $A$ and since
$\ell^{*}$ is an isogeny, it follows that the middle line is
multiplication by $2$ on $\pi^{*}(P(S/R))$. But then the claim
follows since we already know that the vertical arrows are
isomorphisms.
 \qed

A similar proof shows that we also have the following result.

\begin{lemma}
\label{lem:mult2}
The composition $$ \xymatrix@C=0.5cm{
  u^{*}(P(X/T)) \ar[rr]^{\pi^{*} \circ \Nm \pi} & &
   \pi^{*}(P(S/R)) \ar[rr]^{u^{*} \circ \Nm u} & & u^{*}(P(X/T))  }
  $$
is multiplication by $2$.
\end{lemma}

An immediate consequence of Claim \ref{claim:mult2} is that $\ker
H \subseteq P(S/R)[2]$. Now  it follows from Diagram
(\ref{dia:isotrig}) that $\ker H$ is contained in the following
set.

\begin{equation}
\label{eq:F}
 F = \{ x \in P(S/R)[2] : (1+\tau+\tau^2)\ell^{*}\pi^{*}(x) = 0 \}
\end{equation}

In fact, we can be more precise.
\begin{rem}
Recall that if $\alpha = 0$, it follows from the description given
in the table in part i) of Theorem \ref{T:symm} that the covers
$\ell :W \to Z$, $v : Y   \to X$ and $a : C \to S$ are unramified
of degree two, and therefore the kernels of the respective
$*$-induced maps on Jacobians are non trivial subgroups of order
two of $JZ[2]$, $JX[2]$ and $JS[2]$ respectively.

If moreover $\delta = 0$ then we must have $T \neq {\mathbb P}^1$
for $\Delta \to T$ to be connected.
\end{rem}

\begin{claim}
 \label{claim:nsbig}
If $\alpha = 0$ then there exists $s \in F$ such that $u^{*}(H(s))
= \eta_{\ell}$. In particular, $\eta_{\ell} \in u^{*}(P(X/T))$.
\end{claim}

\begin{proof}
Note that $u^{*}$ is always injective and that $u^{*}(\eta_v) =
\eta_{\ell}$; hence it is enough to show that there is $s \in
P(S/R)[2]$ such that $H(s) = \eta_v$ because if this holds then $0
= \ell^{*}(u^{*}(\eta_v)) = \ell^{*}(u^{*}(H(s)) =
(1+\tau+\tau^2)(\ell^{*}(\pi^{*}(s)))$ and therefore $s \in F$.

Furthermore, note that
 $$
 \Nm u(\eta_{\ell}) = \Nm u(u^{*}(\eta_v)) = 3 \eta_v = \eta_v
 $$
and therefore
 $$
 u^{*}(\Nm u(\eta_{\ell})) = \eta_{\ell}
 $$

We also observe that since $f^{*} : JT \to JX$ is always injective
(because it does not factor), we always have $f^{*}(JT[4])
\subseteq P(X/T)[4]$.

We consider two cases.
\begin{enumerate}
\item[i)] \textit{Case 1: $\delta > 0$:} In this case it follows
from Case III in the Appendix that $\eta_a \in P(S/R)[2]$.
Furthermore, it is clear that $\pi^{*}(\eta_a) = \eta_{\ell}$, by
commutativity of the following subdiagram and by the injectivity
of $\pi^{*}$.
 $$ \xymatrix{
     &  W \ar[dl]_{\ell} \ar[d]^{\nu}    \\
 Z \ar[d]_{\pi} & C \ar[dl]^{a}    \\
 S &
 }
 $$

But then
 $$ \eta_v = \Nm u (\eta_{\ell}) = \Nm u (\pi^{*}(\eta_a)) =
 H(\eta_a)
 $$
and we are done.
\item[ii)] \textit{Case 2: $\delta = 0$:} In this case
$\eta_{\delta}$ is a point of order two in $JT$ and therefore
$f^{*}(\eta_{\delta}) = \eta_v \in P(X/T)$.

But we also have
 $$
 \ell^{*}(\pi^{*} \circ \Nm \pi(\eta_{\ell})) = (1+(k \
 l))(\ell^{*}(\eta_{\ell})) = 0
 $$
where the first equality follows from Diagram \ref{dia:isotrig},
and hence $\pi^{*} \circ \Nm \pi(\eta_{\ell}) \in \ker \ell^{*} =
\{ 0, \eta_{\ell} \}$.

We will now show that $\pi^{*} \circ \Nm \pi(\eta_{\ell}) = 0$: if
not, then
 $$
 (\pi^{*} \circ \Nm \pi)(\eta_{\ell}) = \eta_{\ell}
 $$
and it follows that
 $$
 (u^{*} \circ \Nm u) \circ (\pi^{*} \circ \Nm \pi)(\eta_{\ell})
 = (u^{*} \circ \Nm u)(\eta_{\ell}) = \eta_{\ell} \, .
 $$

 But we also have that $(u^{*} \circ \Nm u) \circ (\pi^{*} \circ
\Nm \pi)$ is multiplication by $2$ on $u^{*}(P(X/T))$, from Lemma
\ref{lem:mult2}, which is a contradiction.

Having proven that $(\pi^{*} \circ \Nm \pi)(\eta_{\ell})= 0$, note
that then
 $$
 (u^{*} \circ \Nm u)^{-1}(\eta_{\ell}) \subseteq \ker ((\pi^{*}
 \circ \Nm \pi) \circ (u^{*} \circ \Nm u)) = (\pi^{*}(P(S/R))[2]
 $$
where the last equality follows from Claim \ref{claim:mult2}.

Hence there is $z \in (\pi^{*}(P(S/R))[2]$ such that $(u^{*} \circ
\Nm u)(z) = \eta_{\ell}$; since $\pi^{*}$ is an isomorphism on
$P(S/R)$, there exists $ s \in P(S/R)[2]$ such that $(u^{*} \circ
\Nm u)(\pi^{*}(s)) = \eta_{\ell}$.
\end{enumerate}
\end{proof}

We can now prove the following result, interesting on its own.

\begin{prop}
\label{prop:big}
$\ker H = F$ if and only if $\alpha \neq 0$.

Moreover, if $\alpha = 0$ then $[F : \ker H] = 2$.
\end{prop}

\begin{proof}
Assume there exists $s \in F$ with $s \notin \ker H$. Then
$u^{*}(H(s))$ is the non-zero element of $\ker \ell^{*}$ in
$u^{*}(P(X/T))$, because $u^{*} : JX \to JZ$ is injective as we
know from Lemma \ref{kerl} and by the commutativity of Diagram
\ref{dia:isotrig}; therefore $\alpha = 0$.

Furthermore, it follows that $H(s) \in P(X/T)[2]$ is the non-zero
element of $\ker v^{*}$ and hence the difference of two such $s$
is an element of $\ker H$; therefore $[F : \ker H] = 2$ in this
case. The result now follows from Claim \ref{claim:nsbig}.
\end{proof}

To compute the cardinality of $\ker H$ it now suffices to compute
the cardinality of $F$.

Towards this goal we now prove the following result.

\begin{claim}
\label{claim:prt}
With the notation of Diagram (\ref{dia:symm}))
we have the following description.
 $${r^{*}}^{-1}(P(U/\Delta)[2]) = P(R/T)[2] =  \{ x \in JR[2] :
(1+\tau+\tau^2)\ell^{*}\pi^{*}p^{*}x = 0 \}.
 $$
\end{claim}

Proof of the claim: We recall from the theory of $S_3-$actions
(see \cite{s3}) that  $r^{*} : P(R/T) \to P(U/\Delta)$ is
injective and also that $\ker \{ P(R/T) \times P(R/T)
\stackrel{r^{*}+\tau r^{*}}{\longrightarrow} P(U/\Delta) \} = \{
(x,y): r^{*}x =\tau r^{*}y, x \in P(R/T)[3] \}$.

Hence $r^{*}P(R/T)[2] \ \cap \ \tau r^{*}P(R/T)[2] = \{ 0 \}$ and
therefore
$$ r^{*}P(R/T)[2] + \tau r^{*}P(R/T)[2] = P(U/\Delta)[2]
$$
\noindent (by counting cardinalities).

The first equality in the claim follows now from the injectivity
of $r^{*}$ on $P(R/T)$. Noting that
$(1+\tau+\tau^2)\ell^{*}\pi^{*}p^{*} = (1+\tau+\tau^2)
\varepsilon^{*} r^{*} = \varepsilon^{*} (1+\tau+\tau^2) r^{*}$ and
using the first equality, we obtain that $P(R/T)[2] \subseteq \{ x
\in JR[2] : (1+\tau+\tau^2)\ell^{*}\pi^{*}p^{*}x = 0 \}$.

Conversely, given $x \in \{ x \in JR[2] :
(1+\tau+\tau^2)\ell^{*}\pi^{*}p^{*}x = 0 \}$, then
$(1+\tau+\tau^2) r^{*}(x) \in \ker \varepsilon^{*}$.

Since $\ker \varepsilon^{*} \subseteq \ker(1+\tau + \tau^2)$ by
Claim \ref{claim:eps} in the proof of Theorem \ref{T:alt}, it
follows that $(1+\tau + \tau^2)(1+\tau + \tau^2)r^{*}(x) = 0$.

But $(1+\tau + \tau^2)(1+\tau + \tau^2)r^{*}(x) = 3 (1+\tau +
\tau^2)r^{*}(x) = (1+\tau + \tau^2)r^{*}(3x) = (1+\tau +
\tau^2)r^{*}(x)$ since $3x =x$. Using the first equality again,
the claim is proved.
 \qed

An immediate corollary of this claim is the following result.
\begin{claim}
\label{claim:p}
If $p : S \to R$ is unramified, then $\eta_p$ is in $P(R/T)$.
\end{claim}

Proof of the claim: Since $p$ unramified is equivalent to $2\gamma
+2\delta =0$, in this case we have that all maps $c_i :C_i \to U$
are unramified.

Furthermore, it is clear that $r^{*} \eta_p = \eta_{c_i}$ for the
unique $i$ such that $C_i$ covers $S$. But then $r^{*} \eta_p$
belongs to $\ker \varepsilon^{*} \subseteq P(U/\Delta)[2]$ and
this claim follows from Claim \ref{claim:prt}.
 \qed

\begin{rem}
\label{rem:F}
Recall that we are interested in describing the set
$F$ defined in (\ref{eq:F}).

Note that $p^{*}(P(R/T)[2]) \cap P(S/R)[2]$ is contained in $F$,
by Claim \ref{claim:prt}.

The general philosophy to complete the description is based on
proving that the complementary part of $F$ arises from some
specific elements of $P(S/R)[2]$ which come from the ramification
of $p : S \to R$.
\end{rem}

In particular, if $p$ is unramified we should already have a
description of $F$.

Our next result shows that this is the case, even if $p$ has two
ramification points.

\begin{prop}
\label{prop:almclass}
 If $2\gamma + 2 \delta = 0$ or $2$, then
 $$  F = p^{*}(P(R/T)[2]) \cap P(S/R)[2]
 $$
\end{prop}

\begin{proof}
If $x$ is any element of $F$, then $x \in P(S/R)[2]$ and
$(1+\tau+\tau^2)\ell^{*}\pi^{*}(x) = 0$.

If $2\gamma + 2 \delta = 0$ we know that $P(S/R)[2] = p^{*}(\{ 0,
\eta_p \}^{\bot})$ and if $2\gamma + 2 \delta = 2$ then $P(S/R)[2]
= p^{*}(JR[2])$.

In both cases there exists $y \in JR[2]$ such that $p^{*}(y) = x$;
it now follows from Claim \ref{claim:prt} that $y \in P(R/T)[2]$
and the result is proved.
\end{proof}

We will analyze a further special case in the next section.

\subsection{The classical case of the trigonal construction}

Since the trigonal construction has been very useful in the theory
of Prym varieties of unramified double covers, we devote this
paragraph to it.

We say that we are in \textit{the classical case} if the curve $W$
with ${\mathcal S}_4$-action --as in Diagram (\ref{dia:symm})-- \,
is such that
\begin{enumerate}
\item[i)] $W/{\mathcal S}_4 = T = {\mathbb P}^1$, and such that
\item[ii)]  the canonical polarization on $JS$ induces
$\lambda_P = 2 \lambda$ twice a principal polarization on
$P(S/R)$.
\end{enumerate}

In this case $P(X/T) = JX$ is also a principally polarized abelian
variety.

Note that condition ii) occurs precisely when $p :S \to R$ is
unramified or when it has two ramification points; equivalently,
when $2\gamma + 2 \delta = 0$ or $2$.

Condition i) forces the double cover $d : \Delta \to {\mathbb
P}^1$ to have at least two points of ramification, since $\Delta$
must be connected; equivalently, $\alpha + \delta$ must be even
and greater or equal to two.

Observe that both conditions together exclude the possibility
$\alpha =0$ for the classical case, and also that they imply that
the triple cover $R \to T = {\mathbb P}^1$ has at least one simple
ramification point.

\begin{rem}
\label{rem:classical}

We could also say that the classical case corresponds to starting
with a tetragonal curve $X$ with at least one simple ramification
point and either no total ramification points nor points of type
$(2,2)$, or no total ramification point and one ramification point
of type $(2,2)$, or with one total ramification point and no
ramification point of type $(2,2)$.

Or, equivalently, to a double cover, either unramified or with two
ramification points, of a trigonal curve with at least one simple
ramification point: $S \to R \to {\mathbb P}^1$.

The trigonal construction (see \cite{recillas:jeg94}) shows that
these two situations are equivalent and, furthermore, that then
$P(S/R)$ and $JX$ are isomorphic as principally polarized abelian
varieties.

In both cases the corresponding Galois cover is given by the group
${\mathcal S}_4$ and we are in the situation of Theorem
\ref{T:symm} with $T = {\mathbb P}^1$, $\alpha >0$ and
$2\gamma+2\delta = 0$ or $2$; i.e., in the classical case.

We will now prove that we can also obtain from our methods that
the two principally polarized abelian varieties are isomorphic.
\end{rem}

We will first compute $\ker H$ for the classical case.
\begin{prop}
If $T = {\mathbb P}^1$ and either $2\gamma+2\delta = 0$ or $2$,
then the kernel of the morphism $H = \Nm u \circ \pi^{*} : P(S/R)
\to JX$ is $P(S/R)[2]$.
\end{prop}
\begin{proof}
Note first that if $T = {\mathbb P}^1$, it follows from Claim
\ref{claim:prt} that $JR[2] = \{ x \in JR[2]:
(1+\tau+\tau^2)\ell^{*}\pi^{*}p^{*}x = 0 \}$. Therefore
$p^{*}(JR[2]) \subseteq F$; applying Proposition
\ref{prop:almclass} we obtain $F = P(S/R)[2]$.
It follows from $\alpha \neq 0$ and  Proposition \ref{prop:big}
that
 $$ \ker H =  P(S/R)[2]
 $$
in the classical case.
\end{proof}
Now we can prove the following result.
\begin{theorem}
\label{thm:classical}

If $T = {\mathbb P}^1$ and either $2\gamma+2\delta = 0$ or $2$,
then the morphism $H = \Nm u \circ \pi^{*}|_{P(S/R)} : P(S/R) \to
JX$ induces an isomorphism between the principally polarized
abelian varieties $(JX,\lambda_{JX})$ and $(\wh{P(S/R)},
\lambda_{\wh{P(S/R)}})$.
\end{theorem}
\begin{proof}
If we denote by $\lambda_P$ the polarization on $P= P(S/R)$
induced by the natural principal polarization on $JS$, it follows
from our hypothesis that there exists a principal polarization
$\lambda$ on $P$ such that $\lambda_P = 2\lambda$.
Then note that, since $\ker \lambda_P = P[2] = \ker \left( P
\stackrel{\Nm u \circ \pi^{*} }{\longrightarrow} JX \right)$,
there exists an isomorphism $F : \wh P \to JX$ such that the
following diagram commutes.
$$ \xymatrix{
  P \ar[d]_{\lambda_P = 2\lambda} \ar[rr]^{\Nm u \circ \pi^{*}} & & JX        \\
  \wh P \ar[urr]_{F}^{\thickapprox}                     }
 $$
We now show that the isomorphism $F$ of tori is also  an
isomorphism of principally polarized abelian varieties.
If we denote by $\lambda_1$ the polarization on $\wh P$ induced
via $F$ by $\lambda_{JX}$, we may complete the above diagram to
the following one.
$$\xymatrix{
  P \ar[d]_{\lambda_{P}} \ar[rr]^{\Nm u \circ \pi^{*}} & &
              JX  \ar[d]^{\lambda_{JX}}  \\
  \wh{P} \ar[d]_{\lambda_1} \ar[urr]_{F}^{\thickapprox}& &
              \wh{JX} \ar[dll]^{\wh F}_{\thickapprox}   \\
  P
  }
$$
It now follows from the commutativity of the above diagram that
$\lambda_1$ is principal and that $\ker (\lambda_1 \circ
\lambda_P) = P[2]$; therefore, $\lambda_1 = \lambda_{\wh P}$, as
claimed.
\end{proof}

\subsection{Completion of the principally polarized $P(S/R)$ case}

In this section we  compute $\ker H$ for the remaining principally
polarized cases: $T \neq {\mathbb P}^1$ and either
$2\gamma+2\delta = 0$ or $2$.

\begin{prop}
Assume $T \neq {\mathbb P}^1$ and $2\gamma+2\delta = 0$ or $2$.
Then the cardinality of $\ker H$ is given as follows.
{\small{
\begin{equation*}
 \bigl| \ker H \bigr| =
\begin{cases}
       2^{4g-5+\alpha+2\beta},
             & \text{if $\gamma = \delta = 0$, $\alpha >0$} \\
             & \text{and $P(R/T)[2] \subseteq \{ 0, \eta_p \}^{\bot}$;}   \\
       2^{4g-6 + 2\beta},
             & \text{if $\gamma = \delta = 0 = \alpha$}   \\
             & \text{and $P(R/T)[2] \subseteq \{ 0, \eta_p \}^{\bot}$;}   \\
       2^{4g-6+\alpha+2\beta},
             & \text{if $\gamma = \delta = 0$, $\alpha >0$}   \\
             & \text{and $P(R/T)[2] \not\subseteq \{ 0, \eta_p \}^{\bot}$;}   \\
       2^{4g-7+2\beta},
             & \text{if $\gamma = \delta = 0 = \alpha$}   \\
             & \text{and $P(R/T)[2] \not\subseteq \{ 0, \eta_p \}^{\bot}$;} \\
       2^{4g-4+\alpha+2\beta+\delta},
             & \text{if $2\gamma+2\delta = 2$ and $\alpha>0$;}   \\
       2^{4g-5+2\beta+\delta},
             & \text{if $2\gamma+2\delta = 2$ and $\alpha=0$}   \\
\end{cases}
\end{equation*}
}}
\end{prop}

\begin{proof}
Recall that  $F = p^{*}(P(R/T)[2]) \cap P(S/R)[2]$ holds under our
hypothesis, from Proposition \ref{prop:almclass}.
If $2\gamma+2\delta = 2$, we know that $p^{*}$ is injective and
$P(S/R)[2] = p^{*}(JR[2])$; it follows that $F =
p^{*}(P(R/T)[2])$, with $\bigl| F \bigr| =
2^{4g-4+\alpha+2\beta+\delta}$. We may now apply Proposition
\ref{prop:big} to conclude that if $\alpha$ is positive then $\ker
H = F$ and if  $\alpha$ is zero then  $[F : \ker H] = 2$.
If $2\gamma+2\delta = 0$ we know that $\eta_p \in P(R/T)$, from Claim
\ref{claim:p}, and also that $P(S/R)[2]= p^{*}(\{ 0, \eta_p
\}^{\bot})$ and is isomorphic to $\{ 0, \eta_p \}^{\bot}/\{ 0, \eta_p
\}$.

Therefore we have to analyze two separate cases:

\begin{enumerate}
\item $P(R/T)[2] \subseteq \{ 0, \eta_p \}^{\bot}$, or
\item $P(R/T)[2] \not\subseteq \{ 0, \eta_p \}^{\bot}$, which
means $P(R/T)[2] \cap \{ 0, \eta_p \}^{\bot}$ is of ${\mathbb Z}/2{\mathbb
Z}$-codimension $1$ in $P(R/T)[2]$.
\end{enumerate}
In Case 1) we obtain $F = p^{*}(P(R/T)[2])$, which is isomorphic
to $P(R/T)[2]/\{ 0, \eta_p \}$, and therefore $\bigl| F \bigr| =
2^{4g-5+\alpha+2\beta+\delta}$. In Case 2) $F = p^{*}(P(R/T)[2]
\cap \{ 0, \eta_p \}^{\bot})$ is isomorphic to $P(R/T)[2] \cap \{
0, \eta_p \}^{\bot}/\{ 0, \eta_p \}$, and therefore $\bigl| F
\bigr| = 2^{4g-6+\alpha+2\beta+\delta}$. Now each of the two cases
splits into two more, depending on whether $\alpha > 0$ (with
$\bigl| \ker H \bigr| = \bigl| F \bigr|$) or $\alpha = 0$, with
$\bigl| \ker H \bigr| = \bigl| F \bigr|/2$).
\end{proof}

\subsection{The general case}
We may now assume that $2\gamma + 2 \delta >  2$; we will
continue the description of $F$ by constructing the elements
of $P(S/R)[2]$ coming from the ramification. Then we will
decide which of those lie in $F$.

First note that there are four types of points in $T$ over which
$f : X \to T$ ramifies: the $\alpha$, $\beta$, $\delta$ and
$\gamma$ points, corresponding to the images of simple, triple,
total or $(2,2)$ type of ramification points, respectively.

Their preimages via $g : R \to T$ are the places over which $p : S
\to R$ may ramify. Careful consideration of the group actions
involved shows that the $\alpha$ and the $\beta$ points do not
contribute and that the $\delta$ and $\gamma$ points do
contribute, in the following way.

{\em The $\gamma$ points:} If $\gamma_t \in T$ is a $\gamma$ point
we will denote
 $$
 g^*(\gamma_t) = {\gamma_1}{'} + {\gamma_2}{'} +\gamma_3{'} \ , \
 {\gamma_i}{'} \in R
 $$
where
 $$
 p^*(\gamma_1{'}) = 2\gamma_1, p^*(\gamma_2{'}) = 2\gamma_2 \text{ and }
 p^*(\gamma_3{'}) = \gamma_5 + \gamma_6 \ , \ \gamma_j \in S .
 $$

With respect to  $r : U \to R$, we choose $n \in JR$ such that
 $$
 n^{\otimes 2} = {\mathcal O}_{R}(\gamma_1{'}-\gamma_2{'}) \ \text{
 and } \ r^*(n) \in \ker (1+\tau+\tau^2)
 $$

Then
 $$ {\mathcal G} = {\mathcal O}_S(\gamma_2-\gamma_1) \otimes p^*(n)
 $$
is in $P(S/R)[2]$ and we also have
 $$
 (1+\tau+\tau^2)\ell^{*}\pi^{*}({\mathcal G}) = 0 \, ;
 $$
that is, we have constructed an element of $F$.

Therefore, if we enumerate the $\gamma$ points of $T$ as
$\gamma_t^1 , \ldots \gamma_t^{\gamma}$, we have that the
corresponding points $\gamma_1^1, \gamma_2^1, \ldots
\gamma_1^{\gamma}, \gamma_2^{\gamma}$ in $S$ are ramification
points of $p : S \to R$, and as above we construct, for each $i$
in $\{ 1, \ldots , \gamma \}$, elements in $F$ given as follows
 $$ {\mathcal G}_{2i-1} = {\mathcal O}_S(\gamma_2^i-\gamma_1^i) \otimes
 p^*(n_i)
 $$
\noindent with $n_i \in JR$, $n_i^{\otimes 2} = {\mathcal
O}_{R}({\gamma_1^{i}}{'}-{\gamma_2^i}{'})$  and
$(1+\tau+\tau^2)r^*(n_i) = 0$.

Similarly, we construct elements of $P(S/R)[2]$ as follows.
 $$ {\mathcal G}_{2i} = {\mathcal O}_S(\gamma_1^{i+1}-\gamma_2^i) \otimes p^*(m_i)
 $$
\noindent with $m_i \in JR$, $m_i^{\otimes 2} = {\mathcal
O}_{R}({\gamma_2^{i}}{'}-{\gamma_1^i}{'})$.

\vskip12pt

{\em The $\delta$ points:}

If $\delta_t \in T$ is a $\delta$ point we will denote
 $$
 g^*(\delta_t) = {\delta_1}{'} + 2{\delta_{23}}{'} \ , \
 {\delta_i}{'} \in R
 $$
and
 $$
 p^*(\delta_1{'}) = 2\delta_1, p^*(\delta_{23}{'}) = 2\delta_{23} \ , \
 {\delta_j} \in S \, .
 $$

Next choose $n \in JR$ such that
 $$
 n^{\otimes 2} = {\mathcal O}_{R}(\delta_{23}{'}-\delta_1{'}) \  \text{ with }
   \ r^*(n) \in \ker (1+\tau+\tau^2) \, .
 $$

If we define ${\mathcal D}$ as follows
 $$ {\mathcal D} = {\mathcal O}_S(\delta_1-\delta_{23}) \otimes p^*(n)
 $$
then it is in $P(S/R)[2]$ and we also have
 $$
 (1+\tau+\tau^2)\ell^{*}\pi^{*}({\mathcal D}) = 0 \, ;
 $$
that is, we have constructed an element of $F$.

Therefore, if we enumerate the $\delta$ points of $T$ as
$\delta_t^1 , \ldots \delta_t^{\delta}$, we have that the
corresponding points $\delta_1^1, \delta_{23}^1, \ldots
\delta_1^{\delta}, \delta_{23}^{\delta}$ in $S$ are ramification
points of $p : S \to R$, and as above we construct, for each $i$
in $\{ 1, \ldots , \delta \}$, elements in $F$ given as follows
 $$
 {\mathcal D}_{2i-1} = {\mathcal O}_S(\delta_1^i-\delta_{23}^i)
 \otimes p^*(n_i)
 $$
\noindent with $n_i \in JR$, $n_i^{\otimes 2} = {\mathcal
O}_{R}({\delta_{23}^{i}}{'}-{\delta_1^i}{'})$  and
$(1+\tau+\tau^2)r^*(n_i) = 0$.

Similarly, we construct elements of $P(S/R)[2]$ as follows.
 $$
 {\mathcal D}_{2i} = {\mathcal O}_S(\delta_{23}^{i+1}-\delta_1^i)
 \otimes p^*(m_i)
 $$
\noindent with  $m_i \in JR$, $m_i^{\otimes 2}
 = {\mathcal O}_{R}({\delta_1^{i}}{'}-{\delta_{23}^i}{'})$.

Finally, if $\gamma \delta > 0$ we consider one more sheaf which
links both cases:
 $$
 {\mathcal L} = {\mathcal O}_S(\delta_{1}^{1}-\delta_2^{\delta}) \otimes p^*(m)
 $$
\noindent with  $m \in JR$ and
 $$
 m^{\otimes 2} = {\mathcal O}_{R}({\delta_2^{\gamma}}{'}-{\delta_{1}^1}{'}) \, .
 $$

\vskip12pt

Now we can apply Corollary \ref{coro:p2} to give a description of
$P(S/R)[2]$ in this case: $2\gamma + 2\delta > 2$.

\small{
\begin{equation}
\hspace*{-0.5cm}
\begin{tabular}{|c|c|} \hline
    &     \\
 $P(S/R)[2]$ &  \text{Case}  \\
   &     \\ \hline
   &     \\
 $p^*JR[2] \oplus_{j=1}^{2\gamma - 1} {\mathcal G}_j {\mathbb Z}/2{\mathbb Z}
           \oplus {\mathcal L} {\mathbb Z}/2{\mathbb Z}
           \oplus_{j=1}^{2\delta - 2} {\mathcal D}_j {\mathbb Z}/2{\mathbb
           Z}$
   & $\delta \gamma > 0$ and $2\gamma+2\delta >2$\\
   &     \\ \hline
   &     \\
 $p^*JR[2] \oplus_{j=1}^{2\gamma - 2} {\mathcal G}_j {\mathbb Z}/2{\mathbb
 Z}$
   & $\delta = 0$ and $\gamma \geq 2$ \\
   &     \\ \hline
   &     \\
 $p^*JR[2] \oplus_{j=1}^{2\delta - 2} {\mathcal D}_j {\mathbb Z}/2{\mathbb Z}$
   & $\delta \geq 2$ and $\gamma = 0$ \\
   &     \\ \hline
\end{tabular}
\end{equation}
}

We are now ready to describe $F$ for the case $2\gamma + 2\delta
> 2$.

\begin{prop}

With the above notation, $F$ is described as follows.

\small{
\begin{equation} \hspace*{-0.5cm}
\begin{tabular}{|c|c|} \hline
   &     \\
 $F$ &  \text{Case}  \\
   &     \\ \hline
   &     \\
 $p^*(P(R/T)[2]) \oplus_{j=1}^{\gamma} {\mathcal G}_{2j-1} {\mathbb Z}/2{\mathbb Z}
            \oplus_{j=1}^{\delta - 1} {\mathcal D}_{2j-1} {\mathbb Z}/2{\mathbb Z}$
   & $\delta \gamma > 0$ and $2\gamma+2\delta >2$\\
   &     \\ \hline
   &     \\
 $p^*(P(R/T)[2]) \oplus_{j=1}^{\gamma - 1} {\mathcal G}_{2j-1} {\mathbb Z}/2{\mathbb  Z}$
   & $\delta = 0$ and $\gamma \geq 2$ \\
   &     \\ \hline
   &     \\
 $p^*(P(R/T)[2]) \oplus_{j=1}^{\delta - 1} {\mathcal D}_{2j-1} {\mathbb Z}/2{\mathbb Z}$
   & $\delta \geq 2$ and $\gamma = 0$ \\
   &     \\ \hline
\end{tabular}
\end{equation}
}
\end{prop}

\vskip12pt

\begin{proof}

Since we are assuming $2\gamma + 2\delta > 2$, $p^{*} : JR \to JS$
is injective; hence $p^{*}(P(R/T)[2]) \subseteq P(S/R)[2]$ and
therefore the factor of $F$ not coming from the ramification is
$p^{*}(P(R/T)[2])$.

For $\delta\gamma>0$, let $\mathcal{F}= \{ {\mathcal G}_{2j-1} ,
{\mathcal D}_{2i-1} \}$ with $ 1 \leq j \leq \gamma $ and $1\leq i
\leq \delta$; if $\gamma \geq 2$ and $\delta = 0$, let
$\mathcal{F}= \{ {\mathcal G}_{2j-1} \}$ with $1 \leq j \leq
\gamma $; and if $\delta \geq 2$ and $\gamma = 0$, let
$\mathcal{F}= \{ {\mathcal D}_{2i-1} \}$ with $1 \leq i \leq
\delta $.

Note that, in each case, the elements of the collection
$\mathcal{F}$ span the sheaves ${\mathcal S}$ which come from the
ramification and such that
$(1+\tau+\tau^2)\ell^{*}\pi^{*}({\mathcal S}) = 0$; i.e., those
elements of $F$ coming from the ramification.

Also, there is exactly one relation among them (c.f. Remark
\ref{rem:rel}).
\end{proof}

In this way we have obtained that $\bigl| F \bigr| =
2^{2(g_R-g_T)+\gamma+\delta-1}=2^{4g-5+\alpha+2\beta+\gamma+2\delta}$
whenever $2\gamma + 2\delta > 2$.

Now we can compute $\bigl| \ker H \bigr|$ for this case.

\begin{prop}
If $2\gamma + 2\delta > 2$, then the cardinality of $\ker H$ is
given as follows.

{\small{
\begin{equation*}
 \bigl| \ker H \bigr| =
\begin{cases}
       2^{4g-5+\alpha+2\beta+\gamma+2\delta},
             & \text{if $2\gamma+2\delta > 2$ and $\alpha >0$;} \\
       2^{4g-6+\alpha+2\beta+\gamma+2\delta},
             & \text{if $2\gamma+2\delta >2$ and  $\alpha=0$.}   \\
\end{cases}
\end{equation*}
}}
\end{prop}

\begin{proof}

We know from Proposition \ref{prop:big} that if $\alpha$ is
positive, then $\bigl| \ker H \bigr| = \bigl| F \bigr|$ and that
if $\alpha = 0$, then $\bigl| \ker H \bigr| = \bigl| F \bigr|/2$.
\end{proof}

We have thus completed the proof of v) in Theorem \ref{T:symm}.

As for vi), let us fix the notation: $V = V_3$, $R = R_3$, $Y =
Y_4$ and $X = X_4$ and let
 $$
 T_0 JW = \mathbf{U} \oplus \mathbf{U'} \oplus \mathbf{V_2}
  \oplus \mathbf{V_3} \oplus \mathbf{V_3'}
 $$
denote the isotypical decomposition of the tangent space to $JW$
at the origin, where $\mathbf{U} = {\mathcal U}^{n_0}$,
$\mathbf{U'} = {\mathcal U'}^{n_1}$, $\mathbf{V_2} = {\mathcal
V}_2^{n_2}$, $\mathbf{V_3} = {\mathcal V}_3^{n_3}$, $\mathbf{V'_3}
= {{\mathcal V}'_3}^{n_4}$ with ${\mathcal U}$, ${\mathcal U'}$,
${\mathcal V}_2$, ${\mathcal V}_3$, ${{\mathcal V}'_3} = {\mathcal
V}_3 \otimes {\mathcal U'}$ the complex irreducible
representations of ${\mathcal S}_4$ of respective degrees $1$,
$1$, $2$, $3$ and $3$.

Following \cite{sa}, we compare actions to obtain the isogeny, as
follows.

A short computation shows that
 \begin{align*}
 (d(b \circ \nu)^{*})_0 (T_0 P(V/R)) & =
 {\mathbf{V'_3}}^{\langle (1 \, 3) (2 \,4) \rangle} \, , \\
 (d(\psi \circ h)^{*})_0 (T_0 P(\Delta/T)) & =
 {\mathbf{U'}} \ \text{ and }\\
 (d \psi^{*})_0 (T_0 P(Y/X)) & =
 {\mathbf{U'}} \oplus {\mathbf{V'_3}}^{\langle (1 \, 2 \, 3)  \rangle}  \, . \\
 \end{align*}

  Therefore it follows from the second and third equalities that
   $$
 T_0 P(Y/X) = (d \, h^{*})_0 (T_0 P(\Delta/T)) \oplus (d
 \psi^{*})_0^{-1} {\mathbf{V'_3}}^{\langle (1 \, 2 \, 3)  \rangle}
   $$

   But we can also prove that
   $$
   (1+\tau+\tau^2){{\mathcal V}'_3}^{\langle (1 \, 3)(2 \, 4) \rangle}
    = {{\mathcal V}'_3}^{\langle (1 \, 2 \, 3) \rangle} \, ,
   $$
hence
   $$
    (1+\tau+\tau^2)(d(b \circ \nu)^{*})_0 (T_0 P(V/R))  =
    (1+\tau+\tau^2){{\mathbf V}'_3}^{\langle (1 \, 3)(2 \, 4) \rangle}
    = {{\mathbf V}'_3}^{\langle (1 \, 2 \, 3) \rangle}
   $$

  If we now observe that on $JW$ we have $1+\tau+\tau^2 = \psi^{*}
  \circ \Nm \psi$, we obtain
   $$
   T_0 P(Y/X) = (d \, h^{*})_0 (T_0 P(\Delta/T)) \oplus
   d(\Nm\psi \circ (b \circ \nu)^{*})_0 (T_0 P(V/R))
   $$
 proving vi).

As for vii), let us fix $Z = Z_{13}$, $S= S_3$, $R = R_3$ and
$X=X_4$.

Since $|P(Z/S) \cap \pi^{*}(JS)| < \infty$ then
 \begin{equation}
 \label{eq:fin}
 |P(Z/S) \cap \pi^{*} \circ p^{*}(P(R/T))| < \infty \, .
 \end{equation}

On the other hand, some computations show that
 \begin{align*}
 T_0 (\ell^{*}\circ \pi^{*} \circ p^{*} P(R/T)) & =
  {{\mathbf V}_2}^{\langle (1 \, 3) \rangle}\\
 T_0 (\ell^{*} \circ u^{*} JX) & =
 {{\mathbf V}_3}^{\langle (1 \, 2),(1 \, 3) \rangle} \oplus
 {{\mathbf U}}  \\
  \intertext{and}
 T_0 \ell^{*}JZ & = {{\mathbf V}_3}^{\langle (1 \, 3) \rangle} \cap
 T_0(\ell^{*}P(Z/X))\oplus {{\mathbf V}'_3}^{\langle (1 \, 3)
 \rangle} \\
 & \hskip36pt \oplus {{\mathbf V}_2}^{\langle (1 \, 3) \rangle} \oplus
 {{\mathbf V}_3}^{\langle (1 \, 2),(1 \, 3) \rangle} \oplus
 {{\mathbf U}} \, .\\
  \end{align*}

These equalities imply that $\pi^{*} \circ p^{*} P(R/T) \subseteq
P(Z/X)$. But
 $$
 \dim P(R/T)+ \dim P(Z/S) = \dim P(Z/X)
 $$
 which together with (\ref{eq:fin}) complete the proof of vii) and
 of Theorem \ref{T:symm}.
\end{proof}

\section{Other applications}

\subsection{Examples}

Throughout the paper we did put some obvious restrictions to the
ramification data in some formulae.

Here we will actually construct curves with given ramification
data and given $G-$action, where $G$ is one of the groups
associated to non-Galois fourfold covers, as considered in this
paper: ${\mathcal S}_4$, ${\mathcal A}_4$ or ${\mathcal D}_4$.

For this we recall the general construction: consider an $n-$fold
cover between complex curves
 $$ f : X \to T \, .
 $$

If we denote by $B(f) = \{ P_1, \ldots , P_{\omega}\} \subseteq T$
the branch locus of $f$, we have an induced homomorphism
 $$ f^{\#} : \Pi_1(T - B(f)) \to {\mathcal S}_n
 $$

By the Monodromy theorem and the Riemann extension theorem, we
know that the covers $f$ (up to isomorphism) are classified by the
 homomorphisms $f^{\#}$ with transitive image (up to inner automorphisms).

Moreover, one knows that ${f^{\#}}^{-1}({\mathcal S}_{n-1})
\approx \Pi_1(X - f^{-1}(B(f)))$ and that $ \ker f^{\#} \approx
\Pi_1(W-\gamma^{-1}(B(f)))$, where $\gamma : W \to T$ is the
corresponding Galois extension of $f$, with group $G = \text{Im }
f^{\#} \subseteq {\mathcal S}_n$.

Also recall that
 \begin{multline*}
 \Pi_1 =
 \Pi_1 (T - B(f)) = \langle \alpha_1, \beta_1, \ldots , \alpha_g,
 \beta_g, \sigma_1, \ldots , \sigma_{\omega} : \\
 \alpha_1 \beta_1 \alpha_1^{-1} \beta_1^{-1} \ldots
 \alpha_g \beta_g \alpha_g^{-1} \beta_g^{-1} =
 \sigma_1 \ldots \sigma_{\omega}
 \rangle
 \end{multline*}
where $g = $ genus of $T$, $\alpha_1, \beta_1, \ldots , \alpha_g,
\beta_g $ are canonical generators for $\Pi_1(T)$ and each
$\sigma_j$ is represented by a trajectory going from the base
point to near $P_j$, around it once in the appropriate direction,
and back.

So, to construct $W$ with $G = {\mathcal S}_4$ action and given
values of $\alpha$, $\beta$, $\gamma$ and $\delta$, where $\omega
= \alpha + \beta +\gamma + \delta$, is equivalent to the
construction of a four-fold cover of curves $f : X \to T$ with
those ramification values, and therefore also equivalent to the
construction of a \textit{surjective} homomorphism
 $$
 f^{\#} : \Pi_1 \to {\mathcal S}_4
 $$
such that
 \begin{align}
 \label{cond}
 \begin{cases}
 f^{\#}(\sigma_1), \ldots , f^{\#}(\sigma_{\alpha})
 & \text{ are transpositions;  } \\
 f^{\#}(\sigma_{\alpha+1}), \ldots ,f^{\#}(\sigma_{\alpha+\beta})
 & \text{ are three-cycles;  } \\
 f^{\#}(\sigma_{\alpha+\beta+1}), \ldots ,f^{\#}(\sigma_{\alpha+\beta+\gamma})
 & \text{ are products of two } \\
 & \text{ disjoint transpositions;  } \\
 f^{\#}(\sigma_{\alpha+\beta+\gamma+1}), \ldots ,f^{\#}(\sigma_{\alpha+\beta+\gamma+\delta})
 & \text{ are four-cycles;  } \\
 \end{cases}
 \end{align}

For the other cases $G = {\mathcal A}_4$ or $G = {\mathcal D}_4$
the condition that $f^{\#}$ be surjective changes to $\text{Im }
f^{\#} = G$ is a transitive group of ${\mathcal S}_4$, whereas
(\ref{cond}) stays the same.

\subsection{Jacobians of curves isogenous to a product of Jacobians}

The equivariant isogeny $\phi_{{\mathcal S}_4}$ of Theorem
\ref{T:symm} induces an isogeny
 $$ \phi : JU \times 3P(S/R) \times 3P(V/R) \to JW
 $$
Hence we obtain curves $W$ with ${\mathcal S}_4$ action whose
Jacobian $JW$ is isogenous to a product of Jacobians if $g_R = 0$.

Since $g_R = 3g_T -2 + \displaystyle{\frac{\alpha+\delta}{2}}
+\beta$, this is equivalent to $g_T= 0$ and $\alpha + \delta
+2\beta = 4$.

We now describe the cases for which such covers actually exist and
the dimension of the corresponding moduli.

Recall from the previous section that we are looking for
surjective homomorphisms
 \begin{equation}
 \label{eq:homsp}
 f^{\#} : \Pi_1 = \langle  \sigma_1, \ldots , \sigma_{\omega} : \,
  \sigma_1 \ldots \sigma_{\omega} = 1  \rangle \to {\mathcal S}_4
 \end{equation}
that satisfy (\ref{cond}) with $\alpha + \delta +2\beta = 4$ and
such that
 \begin{equation}
 \label{eq:homg0}
 f^{\#}(\sigma_1) \cdots f^{\#}(\sigma_{\alpha+\beta+\gamma+\delta})= 1
 \end{equation}

Many a priori possibilities for $ \alpha$, $\beta$, $\gamma$ and
$\delta$ are excluded by our conditions; for instance we already
know that $g_T = 0$ implies $\alpha + \delta > 0$ and even.

A final observation before we actually give all possible cases is
that under our conditions we obtain $g_U = g_{\Delta}$, hence in
fact we are considering the isogeny of Theorem \ref{T:symm} ii) in
the special cases
 $$ \phi_{{\mathcal S}_4} : J\Delta \times 3JS \times 3JV \to JW
 $$

Let us observe that when $T = {\mathbb P}^1$ then $P(S/R)$ is
isogenous to $JX$, so apparently we only need to impose the
condition that $P(V/R)$ be isogenous to a Jacobian. The only way
we know to do this at the moment is that $\dim P(V/R) = 1$, but
this is equivalent to
 $$
 \alpha + \beta +\gamma+\delta = 4
 $$
so all possible cases are included in our next result.

\begin{theorem}
The Jacobian of a curve W with ${\mathcal S}_4$ action is
isogenous to a product of Jacobians, via the isogeny
$\phi_{{\mathcal S}_4}$ of Theorem \ref{T:symm} ii) in the
 form
 $$ J\Delta \times 3JS \times 3JV \to JW
 $$
if the ramification data satisfies the following.

 $$
\begin{tabular}{|c|c|c|c|c|} \hline
        &          &         &          &      \\
  Case & $\alpha$ & $\beta$ & $\gamma$ & $\delta$ \\
       &          &         &          &      \\  \hline
  $I$ & $4$ & $0$ & $\geq 1$ & $0$ \\   \hline
  $II$ & $2$ & $1$ & $\geq 1$ & $0$ \\   \hline
  $III$ & $2$ & $0$ & $\geq 0$ & $2$ \\   \hline
  $IV$ & $1$ & $0$ & $\geq 0$ & $3$ \\    \hline
  $V$ & $1$ & $1$ & $\geq 0$ & $1$ \\   \hline
  $VI$ & $0$ & $0$ & $\geq 0$ & $4$ \\    \hline
  $VII$ & $0$ & $1$ & $\geq 0$ & $2$ \\   \hline
  $VIII$ & $3$ & $0$ & $\geq 0$ & $1$ \\ \hline
\end{tabular}
 $$

\vskip12pt

Furthermore, in the following table we list the genera of the
corresponding curves, the degree of the isogeny $\phi_{{\mathcal
S}_4}$ and the number of respective moduli.

 $$
 \begin{tabular}{|c|c|c|c|c|c|c|}  \hline
    &  &  &  &  &  &  \\
   Case & $g_{\Delta}$ & $g_S$ & $g_V$ & $g_W$ & $\deg \phi$ & moduli  \\
    &  &  &  &  &  &  \\  \hline
   $I$     & $1$ & $\gamma-1$ & $\gamma+1$ & $6\gamma+1$  & $2^{12\gamma+2} 3$ & $\gamma+1$ \\   \hline
   $II$    & $0$ & $\gamma-1$ & $\gamma$   & $6\gamma-3$  & $2^{12\gamma-6}$   & $\gamma$ \\   \hline
   $III$   & $1$ & $\gamma+1$ & $\gamma+1$ & $6\gamma+7$  & $2^{12\gamma+14}3$ & $\gamma+1$ \\   \hline
   $IV$    & $1$ & $\gamma+2$ & $\gamma+1$ & $6\gamma+10$ & $2^{12\gamma+20}3$ & $\gamma+1$ \\   \hline
   $V$     & $0$ & $\gamma$   & $\gamma$   & $6\gamma$    & $2^{12\gamma}$     & $\gamma$ \\   \hline
   $VI$    & $1$ & $\gamma+3$ & $\gamma+1$ & $6\gamma+13$ & $2^{12\gamma+26}3$ & $\gamma+1$ \\   \hline
   $VII$   & $0$ & $\gamma+1$ & $\gamma$   & $6\gamma+3$  & $2^{12\gamma+6}$   & $\gamma$ \\   \hline
   $VIII$  & $1$ & $\gamma$   & $\gamma+1$ & $6\gamma+4$  & $2^{12\gamma+8}3$  & $\gamma+1$ \\ \hline
 \end{tabular}
 $$

\end{theorem}

\begin{proof}
The proof of this result is a direct application of Theorem
\ref{T:symm}, but we must show that the given covers exist and why
other cases are excluded.

In particular, since $\alpha + \delta$ must be positive and even,
it is either $4$ or $2$, in which cases $\beta$ must be $0$ or $1$
respectively. From here we obtain eight cases for the values of
$\alpha$, $\beta$ and $\delta$ as in the first table; the next
step is to show that in the first two cases the value $\gamma = 0$
is excluded and that all the other possibilities actually exist.

This is done by either constructing surjective homomorphisms
 $$
 f^{\#} : \Pi_1 = \langle  \sigma_1, \ldots , \sigma_{\omega} : \\
  \sigma_1 \ldots \sigma_{\omega} = 1  \rangle \to {\mathcal S}_4
 $$
as in (\ref{eq:homsp}) with $\omega = \alpha +
\beta+\gamma+\delta$, where $\alpha +2\beta+\delta = 4$, and
satisfying (\ref{cond}) and (\ref{eq:homg0}), or showing that they
cannot exist, depending on the values of the first table.

First of all it is clear that ${\mathcal S}_4$ cannot be generated
by either $4$ transpositions with trivial product nor by two
transpositions with product a $3$-cycle; therefore, if $\alpha=4$
and $\beta = \delta = 0$, or if $\alpha = 2$, $\beta = 1$ and
$\delta = 0$, then $\gamma >0$.

The proof is completed by writing down specific homomorphisms for
the other cases.

We will illustrate with a couple of examples, as follows.

\textit{Case I: $\alpha = 4$, $\beta = \delta =0$ and $\gamma \geq
1$.} Then $f^{\#}$ is given by
 \begin{align*}
 \sigma_1 & \to (1 \, 2) \\
 \sigma_2 & \to (2 \, 3) \\
 \sigma_3 & \to (3 \, 4) \\
 \sigma_4 & \to (1 \, 3) \\
 \intertext{and,  for $\gamma \equiv 1 \ (2)$}
 \sigma_5 & \to (1 \, 4)(2 \, 3) \\
 \sigma_6 , \ldots , \sigma_{4+\gamma} & \to (1 \, 2)(3 \, 4), \text{ only if } \gamma > 1 \\
 \intertext{and,  for $\gamma \equiv 0 \ (2)$}
 \sigma_5 & \to (1 \, 2)(3 \, 4) \\
 \sigma_6 & \to (1 \, 3)(2 \, 4) \\
 \sigma_7 , \ldots , \sigma_{4+\gamma} & \to (1 \, 2)(3 \, 4), \text{ only if } \gamma  >  2.
 \end{align*}
\end{proof}

\subsection{Rigid Jacobians with ${\mathcal S}_4$ actions}

This time we have to look for surjective homomorphisms as follows.

   \begin{equation}
 \label{eq:homrig}
 f^{\#} : \Pi_1({\mathbb P}^1-\{ P_1,P_2,P_3\})  \to {\mathcal S}_4
 \end{equation}
that satisfy (\ref{cond}) with $\alpha +\beta + \gamma + \delta  =
3$ and such that
 \begin{equation}
 \label{eq:homrig3}
 f^{\#}(\sigma_1) \cdots f^{\#}(\sigma_{\alpha+\beta+\gamma+\delta})= 1
 \end{equation}

Again $g = g_T = 0$ implies $\alpha + \delta$ even and positive,
and hence $\alpha + \delta = 2$; therefore $\beta +\gamma = 1$.

One verifies then that the only cases that do appear are the
following two special cases from the previous section, the other
possibilities not being realizable.

\textit{Case V with $\gamma = 0$:}i.e., $\alpha = \beta = \delta =
1$.

In this case $g_W = 0$; that is  we obtain a rational function
 $$ \gamma : {\mathbb P}^1 \to {\mathbb P}^1
 $$
of degree $24$, which corresponds to the quotient map by the
action of ${\mathcal S}_4$ (see \cite{klein}).

\textit{Case VII with $\gamma = 0$:}i.e., $\alpha = 0$, $\beta =
1$, $\delta = 2$.

In this case $g_W = 3$ and $S = E$ is an elliptic curve with an
isogeny
 $$
 E \times E \times E \to JW
 $$
of degree $2^6$.

\subsection{One dimensional families of curves with ${\mathcal S}_4$ action}

Here we have to look for surjective homomorphisms
 $$ f^{\#} : \Pi_1  \to {\mathcal S}_4
 $$
where either
 \begin{align*}
 (I) \ \Pi_1 & = \Pi_1({\mathbb P}^1-\{ P_1,P_2,P_3,P_4\}) \\
 \intertext{ respectively }
 (II) \ \Pi_1 & =  \Pi_1(E-\{ P \})
 \end{align*}
where $E$ is an elliptic curve, and such that $f^{\#}$ satisfies
(\ref{cond}) with
 \small{
 \begin{align*}
 (I) \ \alpha +\beta + \gamma + \delta & = 4 & \text{ with }
        \ f^{\#}(\sigma_1) \cdots f^{\#}(\sigma_{4}) & = 1\\
  \intertext{ respectively }
 (II) \ \alpha +\beta + \gamma + \delta & = 1 & \text{ with }
        \ f^{\#}(\alpha_1) f^{\#}(\beta_1) f^{\#}(\alpha_1^{-1})
      f^{\#}(\beta_1^{-1})  &  = f^{\#}(\sigma_{1})
 \end{align*}
}

The next result gives all possible cases.

\begin{theorem}
The one--parameter families of curves $W$ with ${\mathcal S}_4$
action correspond exactly to those in the following table, where
$g_T = 0$ except in the last case, where $g_T = 1$.

\begin{equation}
\begin{tabular}{|c|c|c|c|c|c|c|c|c|c|} \hline
    &  &  &  &  &  &  &  &  &    \\
 Case & $\alpha$ & $\beta$ & $\gamma$ & $\delta$ & $g_R$ & $g_S$ & $g_V$ & $g_X$ & $g_W$  \\
  $I$ & $3$ & $0$ & $0$ & $1$ & $0$ & $0$ & $1$ & $0$ & $4$  \\
  $II$ & $2$ & $2$ & $0$ & $0$ & $1$ & $1$ & $2$ & $0$ & $5$   \\
  $III$ & $2$ & $0$ & $0$ & $2$ & $0$ & $1$ & $1$ & $1$ & $7$   \\
  $IV$ & $2$ & $1$ & $1$ & $0$ & $0$ & $0$ & $1$ & $0$ & $3$   \\
  $V$ & $1$ & $0$ & $0$ & $3$ & $0$ & $2$ & $1$ & $2$ & $10$   \\
  $VI$ & $1$ & $2$ & $0$ & $1$ & $1$ & $2$ & $2$ & $1$ & $8$   \\
  $VII$ & $1$ & $1$ & $1$ & $1$ & $0$ & $1$ & $1$ & $1$ & $6$   \\
  $VIII$ & $0$ & $2$ & $0$ & $2$ & $1$ & $3$ & $2$ & $2$ & $11$   \\
  $IX$ & $0$ & $1$ & $1$ & $2$ & $0$ & $2$ & $1$ & $2$ & $9$   \\
  $X$ & $0$ & $0$ & $0$ & $4$ & $0$ & $2$ & $1$ & $3$ & $13$   \\
  $XI$ & $0$ & $1$ & $0$ & $0$ & $2$ &$3$ & $3$ & $2$ & $9$  \\ \hline
\end{tabular}
\end{equation}

Furthermore, in each case we have respective isogenies to $JW$ as
follows.

\tiny{
\begin{equation}
\begin{tabular}{|c|c|c|} \hline
   &  &  \\
  Case & Isogeny & Degree of Isogeny \\
   &  &  \\
  $I$ & $3E_0 \times E_1$  & $2^8 3$ \\
  $II$ & $3E_3 \times 2E_4$ & $2^{16} 3$ \\
  $III$ & $3E_2 \times 3E_0 \times E_1$ & $2^{14} 3$ \\
  $IV$ & $3E_0$ & $2^6$ \\
  $V$ & $3JS \times 3E_0 \times E_1$ &  \\
  $VI$ & $3E_5 \times 3E_3 \times 2E_4$ & $2^{22}$ \\
  $VII$ & $3E_2 \times 3E_0 $ & $2^{12}$ \\
  $VIII$ & $3JX \times 3E_3 \times 2E_4 $ &  \\
  $IX$ & $3JS \times 3E_0 $ &  \\
  $X$ & $3JS \times 3E_0 \times E_1$  &  \\
  $XI$ & $E \times 2E_6 \times 3E_5 \times 3E_3  $ & $2^{25} 3^4$\\
   & with $E_0 = V$, $E_1 = \Delta$, $E_2 = S$, $E_3 = P(V/R)$,& \\
   &  $E_4 = R$, $E_5 = P(S/R)$, $E_6 = P(R/T=E)$ &  \\ \hline
\end{tabular}
\end{equation}
}
\end{theorem}

\subsection{Prym varieties of  genus seven double covers of
genus three curves, isogenous to a product of elliptic curves}

In \cite{B-C-V} it is shown that the Prym varieties of genus $7$
double covers of genus three curves, branched at $4$ points, are
dense in the the moduli space of abelian fourfolds of type
$(1,2,2,2)$.

Here we describe a one parameter family of such Prym varieties
which are, moreover, isogenous to the product of elliptic curves.

Our example corresponds to Case III of the previous Section, which
can also be described as follows.

Let $E$ be an elliptic curve and consider $E : \hookrightarrow
{\mathbb P}^3$ its projective normal embedding as a degree four
space curve.  Take two hyper-osculating points $P_1$ and $P_2$ of
$E$; i.e., there exist planes $H_1$ and $H_2$ in ${\mathbb P}^3$
such that $H_i \cdot E = 4P_i$ for $i= 1,2$.

Now consider the degree four meromorphic function $f : E \to
{\mathbb P}^1$ obtained by projection from the line $L = H_1 \cap
H_2$; since the construction depends upon one parameter, we have
that the ramification data for $f$ is given by $\alpha = \delta =
2$ and $\beta = \gamma = 0$ (the construction actually depends
only on the corresponding cyclic subgroup of order $4$ of $E[4]$).

If we denote by $V = E_0$, $\Delta = E_1$, $S = E_2$ and $U = E_7$
the corresponding elliptic curves in this case and with a
superscript the genus of the other curve, we obtain the following
version of Diagram (\ref{dia:symm}).

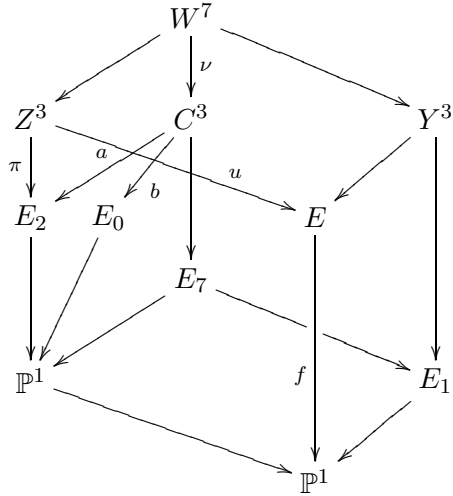
\begin{equation*}
 \xymatrix@R=8pt@C=10pt
 { & & W^7   \ar[ddrrrr] \ar[dd]^{\nu} \ar[ddll] \\
   &  \\
 Z^3 \ar[dd]_{\pi} \ar[ddrrrr]^(.7){u} &               & C^3 \ar[ddd]
\ar[ddll]_{a} \ar[ddl]^(.6){b} & &   & & Y^3 \ar[ddll]
\ar[ddddd]\\
   & &                 & &  \\
 E_2 \ar[ddd]                        & E_0 \ar[lddd] &
                               & & E \ar[ddddd]_(.6){f}  \\
     & & E_7 \ar[ddll] \ar@{-}[drr]   \\
  & &         & &  \ar[drr]  \\
 {\mathbb P}^1  \ar[ddrrrr]   & & & & & & E_1 \ar[ddll]  \\
                     &                                  \\
                     & &                          &  & {\mathbb P}^1
 }
 \end{equation*}

By the general trigonal construction we know that
 $$ \ker(\Nm u \circ \pi^{*} : E_2 \to E) = E_2[2]
 $$
hence $E_2$ and $E$ are isomorphic.

>From the theory of the Klein group action we know that the isogeny
 $$ a^{*}+b^{*} : E_2 \times E_0 \to P(C/E_7)
 $$
is of degree $2^2$ (Case IV in the Appendix).

With the notation of Theorem \ref{T:alt} and again from a Klein
action we have that the isogeny
 $$ \nu_2^{*} + \nu_3^{*} : P(C_2/U) \times P(C_3/U) \to P(W/C_4)
 $$
is of degree $2^4$.

Therefore we have an isogeny
 $$ E_2 \times E_0 \times E_2 \times E_0 \to P(W/C)
 $$
of degree $2^6$.

\section*{Appendix}

In this section we complete the proof of Theorem~\ref{T:klein} by
computing the following quantities  for $\{ j, k, l \}  = \{
\sigma  \tau, \sigma, \tau \}$ and $P_k = P(X_k/T)$.

  \begin{align*}
 \bigl| \ker (\phi_{j}) \bigl| & = \deg a_{k}^*\bigl|_{P_{k}}
 \cdot \deg a_{l}^*\bigl|_{P_{l}} \cdot \bigl| a_{k}^*(P_{k}[2])
 \cap a_{l}^*(P_{l}[2]) \bigl| \\
 \intertext{and}
 \bigl| \ker (\spp_{j}) \bigl| & =|P_{j}[2]| \dfrac{| b_{j}^* JT
 \cap \ker a_{j}^*|}{|\ker a_{j}^*|}
  \end{align*}
where
  \begin{align*}
 \phi_{j} : P_k \times P_l & \to P(X/X_{j}) \, ,& \phi_{j} (x_1, x_2) & = a_{k}^*(x_1) +
 a_{l}^*(x_2) \\
 \intertext{ and }
 \spp_{j} :  P_j \times P(X/X_{j}) & \to P(X/T) \, ,& \spp_j (y, x) & = a_{j}^*(y) + x \, .
 \end{align*}

We will compute by analyzing the different possible cases for the
ramifications of the covers appearing in the following Diagram.

\begin{equation}
 \xymatrix{
 & &  X \ar[dll]_{a_{\sigma}}^{2s} \ar[d]_{a_{\tau}}^{2t}  \ar[drr]^{a_{\sigma \tau}}_{2r} \\
 X_{\sigma} \ar[drr]_{b_{\sigma}}^{t+r} & & X_{\tau}  \ar[d]_{b_{\tau}}^{s+r}
   & &  X_{\sigma \tau} \ar[dll]^{b_{\sigma \tau}}_{s+t}   \\
 & & T
 }  \ .
 \label{klein2}
 \end{equation}

 The possibilities are as follows.
{\em Case I:} All covers in Diagram~(\ref{klein2}) are unramified:
$r=s=t=0$.
{\em Case II:} Exactly two  of the top covers in
Diagram~(\ref{klein2}) are unramified.
{\em Case III:} Exactly two of the top covers in
Diagram~(\ref{klein2}) are ramified.
{\em Case IV:} All covers in Diagram~(\ref{klein2}) are ramified:
$rst \neq 0$.

\bigskip

{\em Case I : $r=s=t=0$.}  In this case, all induced morphisms
between corresponding Jacobians are non-injective.

Let $H_{j} = \ker b_{j}^* = \{ 0, \eta_{b_{j}} \} \subseteq
JT[2]$. Then $H_{l} = \ker b_{l}^* = \{ 0, \eta_{b_j} + \eta_{b_k}
\}$ and
 $$
 \ker \gamma^* = H_{j} + H_{k} = \{ 0, \eta_{b_{\sigma}} ,
 \eta_{b_{\tau}} , \eta_{b_{\sigma}} + \eta_{b_{\tau}} \} \, .
 $$

By \cite{mumprym} we have, for each $j$,  induced isomorphisms
$H_{j}^{\perp}/H_{j} \to b_{j}^*(H_{j}^{\perp}) = P_{j}[2]$.

Therefore
 $$ |P_{\sigma}[2]| = |P_{\tau}[2]| = |P_{\sigma \tau}[2]| =
 2^{2g_T -2} \ \, \text{ if } \, r = s = t = 0.
 $$

We now distinguish two subcases, according to the value of the
Weil pairing $(\eta_{b_{\sigma}} , \eta_{b_{\tau}})$ of
$\eta_{b_{\sigma}} $ and $\eta_{b_{\tau}}$.

\vskip12pt

{\em Case a):} Assume $(\eta_{b_{\sigma}} , \eta_{b_{\tau}}) = 0$;
or, equivalently, $\ker \gamma^* \subseteq H_{j}^{\perp} \cap
H_{k}^{\perp}$ for some  pair (equivalently, each pair) $j,k$.

Then $\ker a_{l}^* = b_{l}^*(\ker \gamma^*) \subseteq P_{l}[2]$
for each $l$ and it follows that for each $l \in  \{ \sigma \tau,
\sigma, \tau \}$ we have
 $$ \deg a_{l}^*\bigl|_{P_{l}} = 2 \, , \ \,  \text{ if } \, r = s = t =
 0 \text{ and } (\eta_{b_{\sigma}} , \eta_{b_{\tau}}) = 0 \, .
 $$

On the other hand,
 $$
 |a_{j}^*(P_{j}[2]) \cap a_{k}^*(P_{k}[2])| = \left|
 \dfrac{H_{j}^{\perp}}{H_{j} + H_{k}} \cap
 \dfrac{H_{k}^{\perp}}{H_{j} + H_{k}} \right|
 $$
but our assumption for this case ($\ker \gamma^* = H_{j} + H_{k}
\subseteq H_{j}^{\perp} \cap H_{k}^{\perp}$)  implies
that
 $$
 \dfrac{H_{j}^{\perp}}{H_{j} + H_{k}} \cap
 \dfrac{H_{k}^{\perp}}{H_{j} + H_{k}} =
 \dfrac{(H_{j} + H_{k})^{\perp}}{H_{j} + H_{k}}
 $$
from where it follows that for each pair $j,k$ we have
 $$
 |a_{j}^*P_{j}[2] \cap a_{k}^*P_{k}[2]| = |(H_{j} +
 H_{k})^{\perp}/(H_{j} + H_{k})| = 2^{2g_T -4} \, .
 $$

Hence we have obtained that
 $$ \bigl| \ker (\phi_{j}) \bigl| = \deg
a_{k}^*\bigl|_{P_{k}} \cdot \deg a_{l}^*\bigl|_{P_{l}} \cdot
\bigl| a_{k}^*(P_{k}[2]) \cap a_{l}^*(P_{l}[2]) \bigl| =
2^{2g_T-2}
 $$
in this case.

We also have $\ker a_{j}^* = b_{j}^*(\ker \gamma^*) \subseteq
b_{j}^*JT$ and therefore
 $$
 \ker a_{j}^* \cap b_{j}^*JT = \ker a_{j}^*  \ \text{ for each } j
 \, .
 $$

It follows that
 $$\bigl| \ker (\spp_{j}) \bigl| =|P_{j}[2]|
 \dfrac{| b_{j}^* JT \cap \ker a_{j}^*|}{|\ker a_{j}^*|} =
 |P_{j}[2]| = 2^{2g_T-2} \ .
 $$

\vskip12pt

{\em Case b):} Assume $(\eta_{b_{\sigma}} , \eta_{b_{\tau}}) \neq
0$; or, equivalently, $(H_{j} + H_{k}) \cap H_{j}^{\perp} = H_{j}$
for some pair (equivalently, any pair) $j, k$.

Then it follows that $a_{l}^*\bigl|_{P_{l}}$ is injective for
every $l$.

Now $H_{l}^{\perp}$ is of index two in $JT[2]$ and hence
 $$
 \gamma^* (H_{l}^{\perp}) = \gamma^* JT[2] \, ;
 $$
therefore,
 $$
 |a_{l}^*P_{l}[2] \cap a_{k}^*P_{k}[2]| = |\gamma^*
 (H_{k}^{\perp}) \cap \gamma^* (H_{k}^{\perp})| = |\gamma^* JT[2]|
 = 2^{2g_T-2} \, .
 $$

So again we obtain
 $$ \bigl| \ker (\phi_{j}) \bigl| =
2^{2g_T-2} \ \, \text{ if } \, r = s = t = 0
 $$
 and
 $$\bigl| \ker (\spp_{j}) \bigl| = 2^{2g_T-2} \ \, \text{ if } \, r = s = t = 0
 $$
which complete the proof for the unramified case.

\vskip12pt

{\em Case II:} Without loss of generality, in this paragraph we
assume $r=t=0$ and $s>0$, which imply $s$ is even.

In this case $a_{\sigma}^*$, $b_{\tau}^*$ and $b_{\sigma \tau}^*$
are the only injective induced homomorphisms between Jacobians.

Let $\ker \gamma^* (= \ker b_{\sigma}^*) = \{ 0, \eta_{b_{\sigma}}
\} \subseteq JT[2]$. Then
 $$
 \ker a_{j}^* = \{ 0, b_{j}^*(\eta_{b_{\sigma}}) \}
 \subseteq b_{j}^*(JT[2]) \subseteq P_{j}[2] \text{ for } j \neq \sigma
 $$
where the last inclusion follows from \cite{mumprym} since $b_{j}^*$ is injective.

It follows that $\deg a_{j}^*\bigl|_{P_{j}} = 2$ for $j \neq
\sigma$.

It also follows that $|\ker \spp_j| = |P_j[2]|$ for $j \neq
\sigma$, whereas $|\ker \spp_{\sigma}| = |P_{\sigma}[2]|$ follows
from the injectivity of $a_{\sigma}^*$. Therefore in this case we
have
 $$ |\ker \spp_k| =
 \begin{cases}
  2^{2g_T-2+s} & \text{if } k \neq \sigma \text{ and } r=t=0, s>0;\\
  2^{2g_T-2} & \text{if } k = \sigma \text{ and } r=t=0, s>0\\
 \end{cases}
 $$

Now from $P_{\sigma}[2] = b_{\sigma}^*(\{ 0,
\eta_{b_{\sigma}}\}^{\perp})$ we obtain
$a_{\sigma}^*(P_{\sigma}[2]) = {\gamma}^*(\{ 0,
\eta_{b_{\sigma}}\}^{\perp}) \subseteq {\gamma}^*(JT[2])$.

But $b_{j}^*(JT[2]) \subseteq P_{j}[2]$ for $j \neq \sigma$, and
therefore ${\gamma}^*(JT[2]) \subseteq a_{j}^*(P_{j}[2])$, from
where $a_{\sigma}^*(P_{\sigma}[2]) \cap a_{j}^*(P_{j}[2]) =
a_{\sigma}^*(P_{\sigma}[2])$ and $|a_{\sigma}^*(P_{\sigma}[2])
\cap a_{j}^*(P_{j}[2])|= 2^{2g_T-2}$.

In order to compute $|a_{\sigma \tau}^*(P_{\sigma \tau}[2]) \cap
a_{\tau}^*(P_{\tau}[2])|$ we use Proposition~\ref{prop:case2} as
follows: let $S_1, \ldots , S_s , \tau S_1 , \ldots , \tau S_s$
denote the ramification points of $a_{\sigma}$ in $X$.

Then $a_{\tau}(S_i) = a_{\tau}(\tau (S_i))$ for $i \in \{ 1,
\ldots , s \}$ are the ramification points of $b_{\tau}$ in
$X_{\tau}$ and $a_{\sigma \tau}(S_i) = a_{\sigma \tau}(\sigma \tau
(S_i))$ for $i \in \{ 1, \ldots , s \}$ are the ramification
points of $b_{\sigma \tau }$ in $X_{\sigma \tau}$.

For $i \in \{ 2, \ldots , s \}$ we let $u_i$ in $\pic^{0}(T)$ be
such that
 $$
 u_i^{\otimes 2} = {\mathcal O}_T(\gamma(S_1) - \gamma(S_i))
 $$ and  we define
 $$
 {\mathcal S}_i^{j} = {\mathcal O}_{X_{j}}
 (a_{j}(S_i) - a_{j}(S_1)) \otimes b_{j}^*(u_i) \text{ for }
 j \in \{ \tau , \sigma \tau \} \text{ and } 2 \leq i \leq s \, ;
 $$
then Proposition~\ref{prop:case2} implies that
 $$ P_j [2] = b_{j}^*(JT[2]) \oplus_{i=2}^{s-1}  {\mathcal S}_i^{j}
\ {\mathbb Z}/2{\mathbb Z} \ .
 $$

But $\ker a_{j}^* \subseteq b_{j}^*(JT[2])$ and hence
 $$ a_{j}^*(P_j [2]) = \gamma^*(JT[2]) \oplus_{i=2}^{s-1}  a_{j}^*{\mathcal S}_i^{j}
\ {\mathbb Z}/2{\mathbb Z} \ .
 $$

Now, since $\sigma \tau(S_k) = \tau(S_k)$ for all $k$, we obtain
 \begin{align*}
 a_{\tau}^*({\mathcal S}_i^{\tau}) & = {\mathcal O}_{X}
 (S_i + \tau(S_i) - S_1 - \tau(S_1)) \otimes a_{\tau}^* b_{\tau}^*(u_i)
 \\
 & = {\mathcal O}_{X} (S_i + \sigma \tau(S_i) - S_1 - \sigma \tau(S_1))
 \otimes a_{\sigma \tau}^* b_{\sigma \tau}^*(u_i) \\
 & = a_{\sigma \tau}^*({\mathcal S}_i^{\sigma \tau})
 \end{align*}
and therefore
 $$ |a_{\tau}^*P_{\tau} [2] \cap a_{\sigma \tau}^*P_{\sigma \tau}
 [2]| = |\gamma^*(JT[2]) \oplus_{i=2}^{s-1}
  a_{\tau}^*{\mathcal S}_i^{\tau} \ {\mathbb Z}/2{\mathbb Z}| =
  2^{2g_T+s-3} \ .
 $$

In this way we have obtained that
 $$ |\ker \phi_k | =
    \begin{cases}
    2^{2g_T-1}  & \text{if } k \neq  \sigma \text{ and } r=t=0, s>0; \\
    2^{2g_T-1+s}  &  \text{if } k = \sigma \text{ and } r=t=0, s>0
    \end{cases}
 $$
which concludes case II.

 \vskip12pt

{\em Case III:} Without loss of generality, in this paragraph we
assume $s=0$ and $rt>0$, which imply $r$ and $t$ are even.

In this case $a_{\sigma}^*$ is the only non-injective induced
homomorphism between Jacobians.

In particular $\eta_{a_{\sigma}} \notin b_{\sigma}^*(JT)$. This
proves that $\ker \spp_{\sigma}$ is a subgroup of index two of
$P_{\sigma}[2]$, and hence
 $$
 |\ker \spp_{k}| =
 \begin{cases}
 2^{2g_T-3+r+t} & \text{if } k = \sigma \text{ and } rt > 0, s =
 0; \\
 2^{2g_T-2+r}  & \text{if } k = \tau \text{ and } rt > 0, s =  0;
 \\
 2^{2g_T-2+t}  & \text{if } k = \sigma \tau \text{ and } rt > 0, s =  0
 \end{cases}
 $$

In order to compute $\deg a_{\sigma}^{*} \bigl|_{P_{\sigma}}$, we
now show that $\eta_{a_{\sigma}} \in P_{\sigma}[2]$: if
$\tilde{\tau} : X_{\sigma} \to X_{\sigma}$ denotes the involution
induced by $\tau : X \to X$, then
$a_{\sigma}^{*}(\tilde{\tau}(\eta_{a_{\sigma}})) =
\tau(a_{\sigma}^{*}(\eta_{a_{\sigma}})) = 0$; hence
${\gamma}^{*}(\Nm b_{\sigma} (\eta_{a_{\sigma}})) =
a_{\sigma}^{*}(\eta_{a_{\sigma}} + \tilde{\tau}
\eta_{a_{\sigma}})= 0$.

But ${\gamma}^{*}$ is injective and therefore $\Nm
b_{\sigma}(\eta_{a_{\sigma}}) =0$. Since $\Nm b_{\sigma}$ and
$b_{\sigma}^{*}$ are dual morphisms (see \cite{mumprym}) and since
$\ker b_{\sigma}^{*} = \{ 0 \}$ has only one connected component,
$\ker \Nm b_{\sigma}$ must have only one connected component;
i.e., $\ker \Nm b_{\sigma} = P_{\sigma}$.

The claim follows, and we obtain $\deg a_{\sigma}^{*}
\bigl|_{P_{\sigma}} = 2$.

In order to compute $|a_{\sigma}^*(P_{\sigma}[2]) \cap
a_{l}^*(P_{l}[2])|$ we use Corollary~\ref{coro:p2} as follows: let
$T_1, \ldots , T_t , \sigma(T_1) , \ldots , \sigma(T_t)$ and
$R_1, \ldots , R_r , \sigma(R_1), \ldots , \sigma(R_r)$ denote the
ramification points of $a_{\tau}$ and  $a_{\sigma \tau} $ in $X$
respectively.

Then $\{ a_{\sigma}(R_1), \ldots , a_{\sigma}(R_r) ,
a_{\sigma}(T_1), \ldots , a_{\sigma}(T_t) \}$ are the ramification
points of $b_{\sigma}$ in $X_{\sigma}$, $\{ a_{\tau}(R_1), \ldots
, a_{\tau}(R_r)\}$ are the ramification points of $b_{\tau}$ in
$X_{\tau}$ and $\{ a_{\sigma \tau}(T_1), \ldots , a_{\sigma \tau}
(T_t)\}$ are the ramification points of $b_{\sigma \tau }$ in
$X_{\sigma \tau}$.

For $i \in \{ 2, \ldots , r+t \}$ we let $\lambda_i$ in
$\pic^{0}(T)$ be such that
 \begin{enumerate}
\item $\lambda_i^{\otimes 2} =
 {\mathcal O}_T(\gamma (R_i) - \gamma (R_{i+1}))$, for $1 \leq i \leq r-1$;
\item $\lambda_r^{\otimes 2} =
 {\mathcal O}_T(\gamma (R_r) - \gamma (T_1))$, for $ i = r$;
\item $\lambda_{r+i}^{\otimes 2} =
 {\mathcal O}_T(\gamma (T_i) - \gamma (T_{i+1}))$, for $1 \leq i \leq t-1$.
 \end{enumerate}

\noindent and we define

\begin{enumerate}
\item $ {\mathcal F}_i^{\sigma} =
 {\mathcal O}_{X_{\sigma}} (a_{\sigma} (R_{i+1}) - a_{\sigma}
 (R_i)) \otimes b_{\sigma}^*(\lambda_i) \in P_{\sigma}[2]$ and \\
 $ {\mathcal F}_i^{\tau} =
 {\mathcal O}_{X_{\tau}} (a_{\tau} (R_{i+1}) - a_{\tau} (R_i))
 \otimes b_{\tau}^*(\lambda_i) \in P_{\tau}[2]$, for $1 \leq i \leq r-1$;
\item $ {\mathcal F}_r^{\sigma} =
 {\mathcal O}_{X_{\sigma}} (a_{\sigma} (T_1) - a_{\sigma} (R_r))
 \otimes b_{\sigma}^*(\lambda_r) \in P_{\sigma}[2]$, for $ i = r$;
\item $ {\mathcal F}_{r+i}^{\sigma} =
 {\mathcal O}_{X_{\sigma}} (a_{\sigma} (T_{i+1}) - a_{\sigma} (T_i))
 \otimes b_{\sigma}^*(\lambda_{r+i}) \in P_{\sigma}[2]$ and \\
 $ {\mathcal F}_{r+i}^{\sigma \tau} =
 {\mathcal O}_{X_{\sigma \tau}} (a_{\sigma \tau} (T_{i+1}) - a_{\sigma \tau} (T_i))
 \otimes b_{\sigma \tau}^*(\lambda_{r+i}) \in P_{\sigma \tau}[2]$, for $1 \leq i \leq t-1$.
 \end{enumerate}

Then Corollary~\ref{coro:p2} implies that
 \begin{align}
 \label{eq:psigma}
  P_{\sigma}[2] & = b_{\sigma}^*(JT[2]) \oplus_{i=1}^{r+t-2} {\mathcal F}_i^{\sigma}
\ {\mathbb Z}/2{\mathbb Z} \ , \\
 \label{eq:ptau}
 P_{\tau}[2] & = b_{\tau}^*(JT[2]) \oplus_{i=1}^{r-2} {\mathcal F}_i^{\tau}
\ {\mathbb Z}/2{\mathbb Z} \ , \\
  \intertext{and}
 \label{eq:psigmatau}
 P_{\sigma \tau}[2] & = b_{\sigma \tau}^*(JT[2]) \oplus_{i=1}^{t-2}
 {\mathcal F}_{r+i}^{\sigma \tau} \ {\mathbb Z}/2{\mathbb Z} \ .
 \end{align}

 Observe that $a_{\sigma}^* {\mathcal F}_i^{\sigma} =
 a_{\tau}^* {\mathcal F}_i^{\tau}$, which will be denoted
 by $ {\mathcal F}_i$,  for $1 \leq i \leq r-1$
 and that $a_{\sigma}^* {\mathcal F}_{r+i}^{\sigma} = a_{\sigma \tau}^* {\mathcal
 F}_{r+i}^{\sigma \tau}$, which will be denoted by ${\mathcal
 F}_{r+i}$, for $i \in \{ 1, \ldots, t-1 \}$.

Since $a_{\tau}^*$ and $a_{\sigma \tau}^*$ are injective, applying
them to (\ref{eq:ptau}) and (\ref{eq:psigmatau}) respectively we
obtain
 \begin{align}
 \label{eq:aptau}
 a_{\tau}^* (P_{\tau}[2]) & = \gamma^*(JT[2]) \oplus_{i=1}^{r-2}
 {\mathcal F}_i \ {\mathbb Z}/2{\mathbb Z} \ , \\
 \intertext{and}
 \label{eq:apsigmatau}
 a_{\sigma \tau}^*(P_{\sigma \tau}[2]) & = \gamma^*(JT[2]) \oplus_{i=1}^{t-2}
 {\mathcal F}_i \ {\mathbb Z}/2{\mathbb Z} \ .
 \end{align}

Since $\ker a_{\sigma}^* = \{ 0, \eta_{a_{\sigma}} \} \subseteq
P_{\sigma}[2]$, in order to describe $a_{\sigma}^*(P_{\sigma}[2])$
we first have to express $\eta_{a_{\sigma}}$ in terms of the given
generators for $P_{\sigma}[2]$: by Lemma~\ref{lemma:case2}
 we know that
  $${\mathcal F}_1^{\tau} \otimes \ldots \otimes
  {\mathcal F}_{r-1}^{\tau} = b_{\tau}^* (m) \text{ for some } m \in JT[2]
  $$
and that
  $${\mathcal F}_2^{\sigma} \otimes \ldots \otimes {\mathcal F}_{r-1}^{\sigma}
  \in P_{\sigma}[2] \text{ is \underline{not} in } b_{\sigma}^*JT[2] \, .
  $$

This implies that
 $${\mathcal F}_1^{\sigma} \otimes \ldots \otimes
 {\mathcal F}_{r-1}^{\sigma} \otimes b_{\sigma}^* (m^{-1}) \in P_{\sigma}[2]-
 \{ 0 \}
 $$
but
 \begin{align*}
  a_{\sigma}^*({\mathcal F}_1^{\sigma} \otimes \ldots \otimes
 {\mathcal F}_{r-1}^{\sigma} \otimes b_{\sigma}^*(m^{-1})) & =
 {\mathcal F}_1 \otimes \ldots \otimes
 {\mathcal F}_{r-1} \otimes \gamma^*(m^{-1}) \\
 & =  a_{\tau}^*({\mathcal F}_1^{\tau} \otimes \ldots \otimes
 {\mathcal F}_{r-1}^{\tau} \otimes b_{\tau}^*(m^{-1})) \\
 & = a_{\tau}^*(0) \\
 & = 0
 \end{align*}
and therefore
 $$ \eta_{a_{\sigma}} = {\mathcal F}_1^{\sigma} \otimes \ldots \otimes
 {\mathcal F}_{r-1}^{\sigma} \otimes b_{\sigma}^*(m^{-1}) \ .
 $$

We may now apply  $a_{\sigma}^*$ to (\ref{eq:psigma}) to obtain
  \begin{equation}
  \label{apsigma}
 a_{\sigma}^*(P_{\sigma}[2]) = \gamma^*(JT[2])
  \bigoplus_{ \substack{ i=1 \\ i \neq r-1   }}^{r+t-2}
  {\mathcal F}_i \ {\mathbb Z}/2{\mathbb Z} \ .
  \end{equation}

Equality (\ref{apsigma})  shows that there is no dependence
relation between $\mathcal F_1 , \ldots , {\mathcal F}_{r-2}$ and
 ${\mathcal F}_{r+1} , \ldots , {\mathcal F}_{r+t-2}$
 and comparing  we obtain
 \begin{align*}
 a_{\tau}^*(P_{\tau}[2]) \cap a_{\sigma \tau}^*(P_{\sigma
 \tau}[2]) & = \gamma^*(JT[2]) \ , \\
 a_{\tau}^*(P_{\tau}[2]) \cap a_{\sigma}^*(P_{\sigma}[2])
 & = \gamma^*(JT[2]) \oplus_{i=1}^{r-2}
  {\mathcal F}_{i} \ {\mathbb Z}/2{\mathbb Z} \\
 \intertext{and}
 a_{\sigma}^*(P_{\sigma}[2]) \cap a_{\sigma \tau}^*(P_{\sigma
 \tau}[2]) & = \gamma^*(JT[2]) \oplus_{i=1}^{t-2}
  {\mathcal F}_{r+i} \ {\mathbb Z}/2{\mathbb Z} \ .
 \end{align*}

Thus we have proven
 $$ |\ker \phi_k | =
    \begin{cases}
    2^{2g_T-1+t}  & \text{if } k =  \tau \text{ and } rt>0, s=0;\\
    2^{2g_T-1+r}  & \text{if }  k = \sigma \tau \text{ and } rt>0, s=0; \\
    2^{2g_T}  & \text{if }  k = \sigma \text{ and } rt>0, s=0
    \end{cases}
 $$
which concludes case III.

\vskip12pt

{\em Case IV:} All covers in Diagram~(\ref{klein2}) are ramified:
$rst \neq 0$.

In this case all induced homomorphisms between Jacobians are
injective and therefore
 $$ |\ker \spp_k| = |P_k[2]| =
   \begin{cases}
   2^{2g_T-2+s+r}  & \text{if } k = \tau \text{ and } rst > 0;  \\
   2^{2g_T-2+s+t}  & \text{if }  k = \sigma \tau  \text{ and } rst > 0; \\
   2^{2g_T-2+t+r}  & \text{if }  k = \sigma  \text{ and } rst > 0
   \end{cases}
 $$

We also have $\deg a_{j}^{*} \bigl|_{P_{j}} = 1$ for every $j$;
hence all that is left is to compute
 $$
 |\ker \phi_j| = \bigl| a_{k}^*(P_{k}[2]) \cap a_{l}^*(P_{l}[2]) \bigl| \, .
 $$

By the symmetries involved, it is enough to compute one of them:
we will compute
 $$
 |\ker \phi_{\sigma \tau}| = \bigl| a_{\sigma}^*(P_{\sigma}[2]) \cap a_{\tau}^*(P_{\tau}[2]) \bigl| \, .
 $$

The computation for this case is subdivided into two cases,
according to whether $r+t = 2$ or $r+t>2$.

\vskip12pt

{\em Case a): $r = t = 1$.} In this case $P_{\sigma}[2] =
b_{\sigma}^{*}(JT[2])$ and therefore $\gamma^{*}(JT[2]) =
a_{\sigma}^{*}(P_{\sigma}[2])$.

On the other hand, $b_{\tau}^{*}(JT[2]) \subseteq P_{\tau}[2]$ and
therefore $a_{\sigma}^{*}(P_{\sigma}[2])
 \subseteq a_{\tau}^{*}(P_{\tau}[2])$.

Thus
 $$
 a_{\sigma}^{*}(P_{\sigma}[2]) \cap a_{\tau}^{*}(P_{\tau}[2])
 = \gamma^{*}(JT[2])
 $$
and in this case we obtain
 $$ |\ker \phi_{\sigma \tau}| = 2^{2g_T} \ \text{ if } \, r = t = 1 , s > 0.
 $$

\vskip12pt

 {\em Case b): $r + t > 2$.} Without loss of generality,
we will assume $r>1$.

Let $\{ S_1 , \ldots, S_s , \tau(S_1) , \ldots , \tau(S_s) \}$,
$\{T_1 , \ldots , T_t , \sigma (T_1), \ldots , \sigma (T_t)\}$ and
$\{ R_1 , \ldots , R_r , \sigma (R_1), \ldots , \sigma (R_r)\}$
denote the ramification points of $a_{\sigma}$, $a_{\tau}$ and
$a_{\sigma \tau}$ respectively.

Then the ramification points for $b_{\sigma}$, $b_{\tau}$ and
$b_{\sigma \tau}$ are as follows.

\vskip12pt

\begin{center}
\begin{tabular}{|c|c|} \hline
Map & Ramification points \\ \hline
 $b_{\sigma}$ & $\{ a_{\sigma} (R_1),
 \ldots , a_{\sigma} (R_r), a_{\sigma} (T_1), \ldots , a_{\sigma}
 (T_t) \}$ \\ \hline
 $b_{\tau}$   & $\{ a_{\tau}(R_1) , \ldots ,
 a_{\tau}(R_r), a_{\tau}(S_1) , \ldots , a_{\tau}(S_s) \}$ \\ \hline
 $b_{\sigma \tau}$ & $\{ a_{\sigma \tau}(S_1) , \ldots, a_{\sigma \tau}(S_s),
 a_{\sigma \tau}(T_1) , \ldots , a_{\sigma \tau}(T_t) \}$ \\ \hline
\end{tabular}
\end{center}

\vskip12pt

If we choose  ${\mathcal F}_i^{\sigma} \in
 P_{\sigma}[2]$ for $i \in \{ 1, \ldots , r+t-1 \}$ as in the proof of
 Case III, then (\ref{eq:psigma}) holds.

 Furthermore, we choose
 ${\mathcal F}_j^{\tau} \in P_{\tau}[2]$ for $1 \leq j \leq r-1$ and
 ${\mathcal F}_{r+j}^{\sigma \tau} \in P_{\sigma \tau}[2]$ for
 $1 \leq j \leq t-1$ as in the proof of Case III.

 Also,  for $i \in \{ 1, \ldots , s-1 \}$, choose $m_i$ in $\pic^{0}(T)$
 such that
  $$
  m_{i}^{\otimes 2} = {\mathcal O}_T(\gamma(S_i) - \gamma(S_{i+1}))
  $$
and define
 \begin{align*} {\mathcal G}_{i}^{\sigma \tau} & =
 {\mathcal O}_{X_{\sigma \tau}} (a_{\sigma \tau}(S_{i+1}) - a_{\sigma \tau}(S_i))
 \otimes b_{\sigma \tau}^*(m_{r+i}) \in P_{\sigma \tau}[2] \\
 \intertext{and}
 {\mathcal G}_{r+i}^{\tau} & = {\mathcal O}_{X_{\tau}} (a_{\tau}(S_{i+1}) - a_{\tau}(S_i))
 \otimes b_{\tau}^*(m_{i}) \in P_{\tau}[2] \, .
 \end{align*}

\vskip12pt

With these generators and Corollary~\ref{coro:p2} we obtain the
following descriptions.

\begin{align}
 \label{eq:ptaun}
 P_{\tau}[2] & = b_{\tau}^*(JT[2]) \oplus_{i=1}^{r-1} {\mathcal F}_i^{\tau}
 \ {\mathbb Z}/2{\mathbb Z} \oplus_{i=1}^{s-1}
 {\mathcal G}_{r+i}^{\tau} \ {\mathbb Z}/2{\mathbb Z} \ , \\
\intertext{and}
 \label{eq:psigmataun}
 P_{\sigma \tau}[2] & = b_{\sigma \tau}^*(JT[2]) \oplus_{i=1}^{s-1}
 {\mathcal G}_{i}^{\sigma \tau} \ {\mathbb Z}/2{\mathbb Z} \oplus_{i=1}^{t-1}
 {\mathcal F}_{r+i}^{\sigma \tau} \ {\mathbb Z}/2{\mathbb Z} \ .
 \end{align}

 Note that $a_{\sigma \tau}^{*}({\mathcal G}_i^{\sigma \tau}) =
a_{\tau}^{*}({\mathcal G}_{r+i}^{\tau})$ for $i \in \{ 1, \ldots ,
s-1 \}$; the common value will be denoted by ${\mathcal G}_i$.
Applying $a_{\sigma}^{*}$, $a_{\tau}^{*}$ and $a_{\sigma
\tau}^{*}$ to (\ref{eq:psigma}), (\ref{eq:ptaun}) and
(\ref{eq:psigmataun}) respectively, we obtain the following.

 \begin{align}
 \label{eq:apsig2}
 a_{\sigma}^*(P_{\sigma}[2]) & = \gamma^*(JT[2]) \oplus_{i=1}^{r+t-1}
  {\mathcal F}_i \ {\mathbb Z}/2{\mathbb Z} \ , \\
 \label{eq:aptau22}
 a_{\tau}^*(P_{\tau}[2]) & = \gamma^*(JT[2]) \oplus_{i=1}^{r-1}
  {\mathcal F}_i \ {\mathbb Z}/2{\mathbb Z} \oplus_{i=1}^{s-1}
  {\mathcal G}_i \ {\mathbb Z}/2{\mathbb Z} \ , \\
 \intertext{and}
 \label{eq:apsigtau2}
 a_{\sigma \tau}^*(P_{\sigma \tau}[2]) & = \gamma^*(JT[2]) \oplus_{i=1}^{s-1}
  {\mathcal G}_i \ {\mathbb Z}/2{\mathbb Z} \oplus_{i=1}^{t-1}
  {\mathcal F}_{r+i} \ {\mathbb Z}/2{\mathbb Z} \ ,
 \end{align}

It follows from equalities (\ref{eq:apsig2}), (\ref{eq:aptau22})
and (\ref{eq:apsigtau2}) that

\begin{align*}
a_{\sigma}^* P_{\sigma}[2] \cap a_{\tau}^* P_{\tau}[2] & =
 \gamma^*JT[2] \oplus_{i=1}^{r-1} {\mathcal F}_i
 \ {\mathbb Z}/2{\mathbb Z} \\
a_{\tau}^* P_{\tau}[2] \cap  a_{\sigma \tau}^* P_{\sigma \tau}[2]
& =
 \gamma^*JT[2] \oplus_{i=1}^{s-1} {\mathcal G}_i
 \ {\mathbb Z}/2{\mathbb Z} \\
a_{\sigma}^* P_{\sigma}[2] \cap a_{\sigma \tau}^* P_{\sigma
\tau}[2]
 & = \gamma^*JT[2] \oplus_{i=1}^{t-1} {\mathcal F}_{r+i}
 \ {\mathbb Z}/2{\mathbb Z}
\end{align*}

Therefore, we  have proven that in this case we have
 $$ \ker \phi_{k} =
     \begin{cases}
    2^{2g_T-1+t}  & \text{if } k =  \tau \text{ and } r>1, st > 0;\\
    2^{2g_T-1+r}  & \text{if }  k = \sigma \tau \text{ and } r>1, st > 0;\\
    2^{2g_T-1+s}  & \text{if }  k = \sigma \text{ and } r>1, st > 0;
    \end{cases}
 $$

Thus the proof of Theorem~\ref{T:klein} is now complete.

\bibliographystyle{amsalpha}

\begin{thebibliography}{ceawo}

\bibitem[B-C-V]{B-C-V} F. Bardelli, C. Ciliberto and A. Verra,
 \textit{Curves of minimal genus on a general abelian variety}
  Compositio Math. \textbf{96} (1995), no. 2, 115--147.

\bibitem[D1]{d} R. Donagi,
\textit{The Fibers of the Prym Map},  Curves, Jacobians, and abelian
varieties (Amherst, MA, 1990), pp. 55--125. Contemp. Math.,
\textbf{136}, Amer. Math. Soc., Providence, RI, 1992.

\bibitem[D2]{d2} R. Donagi,
\textit{Decomposition Of Spectral Covers}, Journ\'ees de
G\'eom\'trie Alg\'ebrique d'Orsay (Orsay, 1992). Ast\'erisque
\textbf{218} (1993), 145--175.

\bibitem[G]{gg} E. G\'omez--Gonz\'alez,
 \textit{Prym varieties of curves with an automorphism of
  prime order},
  Aportaciones Mat. \textbf{13} (1998), 103--116.

\bibitem[Ka]{kanev} \textit{Spectral Curves, simple Lie Algebras, and
Prym-Tjurin varieties}, Theta functions -- Bowdoin 1987, Part
1,627--645, Proc. Sympos. Pure Math. \textbf{49} Part 1, Amer.
Math. Soc., 1989.

\bibitem[K]{klein} F. Klein,
 \textit{Lectures on the Icosahedron}.  Dover Publications, New
 York, 1913.

\bibitem[L-B]{lange:cav92} H. Lange and  C. Birkenhake,
\textit{Complex Abelian Varieties}, Grundlehren der Mathematischen
Wissenschaften \textbf{302}. Springer-Verlag, Berlin, 1992.

\bibitem[M]{merindol} J-Y. M\'erindol, \textit{Vari\'et\'es de
Prym d'un rev\^{e}tement galoisien}, J. Reine Angew. Math
\textbf{461} (1995), 49--61.

\bibitem[M1]{mumav} D. Mumford, \textit{Abelian Varieties},
Tata Institute of Fundamental Research, Bombay: Oxford University
Press, 1970.

\bibitem[M2]{mumprym} D. Mumford, \textit{Prym Varieties I},
Contributions to analysis (a collection of papers dedicated to Lipman
Bers), pp. 325--350. Academic Press, New York, 1974.

\bibitem[P]{p} S. Pantazis,
\textit{Prym Varieties and the Geodesic Flow on $SO(n)$},
Mathematische Annalen  \textbf{273} (1986), 297--315.

\bibitem[Rec1]{r} S. Recillas,
\textit{Jacobians of curves with  $g^1_4$'s are  the Prym's of
trigonal curves}, Bol. Soc. Mat. Mexicana (2) \textbf{19} (1974), no.
1, 9--13. \textbf{19} (1974), 9--13.

\bibitem[Rec2]{recillas:jeg94} S. Recillas,
\textit{La Jacobiana de la extensi\'on de Galois de una
  curva trigonal},
Aportaciones Matem\'aticas de la Soc Mat. Mexicana \textbf{14}
(1994), 159--167.

\bibitem[Rec--Ro]{s3} S. Recillas and R. E. Rodr\'\i guez,
\textit{Jacobians and Representations of $S_3$},
 Aportaciones Mat. Inv. \textbf{13} (1998), 117--140.

\bibitem[R]{ri} J. Ries,
\textit{The Prym variety for a cyclic unramified cover of a
hyperelliptic Riemann surface}, J. Reine Angew. Math. \textbf{340}
(1983), 59--69.

\bibitem[SA]{sa} A. S\'anchez-Arg\'aez,
 \textit{Acciones del grupo ${\mathcal A}_5$ en variedades
jacobianas}, Aportaciones Mat. Com. \textbf{25} (1999), 99--108.

\bibitem[W]{w} W. Wirtinger,
\textit{Untersuchungen \"uber Theta Funktionen}. Teubner, Berlin (1895).
\end{thebibliography}

\end{document}